\documentclass[10pt]{amsart}
\usepackage{amssymb}
\usepackage{latexsym}
\usepackage{amscd}
\usepackage{amsmath}
\usepackage{epsfig}
\usepackage{eepic}

\include{psfig}

\newtheorem{theorem}{Theorem}[section]
\newtheorem{introtheorem}{Theorem}
\newtheorem{corollary}[theorem]{Corollary}
\newtheorem{lemma}[theorem]{Lemma}
\newtheorem{definition}[theorem]{Definition}
\newtheorem{proposition}[theorem]{Proposition}
\newtheorem{remark}[theorem]{Remark}
\newtheorem{notation}[theorem]{Notation}
\def\A{\mathbb{A}}
\def\B{\mathbb{B}}
\def\CC{\mathbb{C}}
\def\C{\mathbb{C}}

\def\H{\mathbb{H}}
\def\hepsilon{\H_\epsilon}
\def\R{\mathbb{R}}
\def\Q{\mathbb{Q}}
\def\Qa{\mathbb{Q}^a}
\def\N{\mathbb{N}}
\def\L{\mathbb{L}}
\def\Z{\mathbb{Z}}
\def\D{\mathbb{D}}
\def\S{{\mathbb{S}}}

\def\CDC{\CC \setminus \overline{\D}}
\def\CP2{{\mathbb{CP}^2}}

\def\be{\bar{\mathrm{e}}}
\def\ba{\bar{a}}
\def\bv{\bar{v}}
\def\bb{\bar{\beta}}
\def\bT{{\mathbf{T}}}
\def\bM{{\mathbf{M}}}

\def\cC{{\mathcal{C}}}
\def\cD{{\mathcal{D}}}
\def\cE{{\mathcal{E}}}
\def\cEvphi{\mathcal{E}^\vphi}
\def\cEpsi{\mathcal{E}^\psi}
\def\cF{{\mathcal{F}}}
\def\cH{{\mathcal{H}}}
\def\cL{{\mathcal{L}}}
\def\cM{{\mathcal{M}}}

\def\cO{{\mathcal{O}}}
\def\cP{{\mathcal{P}}}
\def\cS{{\mathcal{S}}}

\def\z{\zeta}
\def\vphi{\varphi}

\def\znot{{\z_0}}

\def\zone{\z_1}

\def\hA{{\hat{A}}}

\def\hnu{{\hat{\nu}}}

\def\ha{{\hat{a}}}
\def\hv{{\hat{v}}}

\def\hp{{\hat{p}}}

\def\hr{{\hat{r}}}
\def\hz{\hat{\zeta}}
\def\hzone{\hat{\zeta}_1}

\def\zp{\z^\prime}
\def\nup{\nu^\prime}
\def\mup{{\mu^\prime}}

\def\tL{\widetilde{\mathbb{L}}}
\def\tS{\widetilde{\mathbb{S}}}

\def\tvarphi{\tilde{\varphi}}
\def\tz{\tilde{\zeta}}
\def\talpha{\tilde{\alpha}}

\def\itin{\mathit{it}}
\def\fav{{f_{a,v}}}

\def\laut{\Qa((t))}
\def\lautm{\Qa((t^{1/m}))}
\def\pui{\Qa \langle\langle t \rangle\rangle }
\def\sol{\S}

\def\order{o}
\def\lval{|}
\def\rval{|_o}
\def\valgL{\lval \L^\star \rval}
\def\valgS{\lval \S^\star \rval}
\def\mdl{\operatorname{mod}}
\def\mult{\operatorname{mult}}
\def\chordal{\operatorname{chordal}}

\def\im{\operatorname{Im}}
\def\mathe{\mathrm{e}}

\def\fatou{F}
\def\julia{J}
\def\filled{K}
\def\diam{\operatorname{diam}}
\def\db{D} %dynamical ball
\def\dbb{D^{\beta}} %dynamical ball
\def\dbnu{D^\nu} %dynamical ball
\def\dbhnu{D^{\hat{\nu}}} %dynamical ball
\def\dbnup{D^{\nu^\prime}} %dynamical ball
\def\pb{\cD} %parameter ball
\def\pbpsi{\cD^\psi} %parameter ball
\def\da{A} %dynamicall annulus
\def\danu{A^\nu} %dynamicall annulus
\def\dn{E}
\def\dN{\mathbf{Ends}}

\def\cubicl{\cP_3 (\L)}
\def\hcubicc{\widehat{\cP}_3(\C)}
\def\cubicc{\cP_3(\C)}
\def\ab{{\alpha,\beta}}
\def\b{{\beta}}

\def\varphiab{{\varphi_{\alpha,\beta}}}
\def\vphiab{{\varphi_{\alpha,\beta}}}
\def\vphib{{\varphi_{\beta}}}
\def\psianu{\psi_{\alpha,\nu}}
\def\psinu{\psi_\nu}
\def\psihnu{\psi_{\hat{\nu}}}
\def\psinup{\psi_{\nu^\prime}}
\def\coneccubicl{\cC_3 (\L)}

\def\hconeccubicc{\widehat{\cC}_3 (\C)}
\def\filledab{\filled (\varphiab)}
\def\shiftcubicl{\cS_3 (\L)}

\def\hshiftcubicc{\widehat{\cS}_3 (\C)}

\def\hescapepc{\widehat{\cE}^+ (\C)}

\def\hescapepmc{\widehat{\cE}^\pm (\C)}

\def\escapepml{\cE^\pm (\L)}

\def\per{\operatorname{Per}}
\def\omep{\omega^+}

\def\bot{\phi}
\def\basinstar{\mathcal{B}_\nabla}
\def\curve{\mathcal{W}}
\def\bav{\overline{\mathrm{av}}}
\def\av{\mathrm{av}}

\def\bottcher{{B\"{o}ttcher }}

\begin{document}

\title{Puiseux series polynomial dynamics and iteration of complex cubic polynomials}

\author{Jan Kiwi}
\address{Facultad de Matem\'aticas,
Pontificia Universidad Cat\'olica,
Casilla 306, Correo 22, Santiago,
Chile.}
\thanks{Supported by ``Proyecto Fondecyt \#1020711''}
\email{jkiwi@puc.cl}
\date{September 15, 2004}
%\keywords{Puiseux series, Julia set}
%\subjclass{37F}
\begin{abstract}
We study polynomials with coefficients in a field $\L$ as dynamical systems where
$\L$ is any  algebraically closed and complete ultrametric field with
dense valuation group and characteristic zero residual field.
 We give a complete description of the dynamical and parameter space of
cubic polynomials. In particular we characterize cubic polynomials with compact Julia sets.
Also, we prove that any infraconnected  connected component of a filled Julia set (of a cubic polynomial)
is either a point or eventually periodic. 
 A smallest field $\S$ with the above properties is, up to isomorphism, the completion of the field of formal Puiseux series
with coefficients in an algebraic closure of $\Q$. 
 We show that some elements of $\S$ naturally correspond to the Fourier series of analytic almost
periodic functions (in the sense of Bohr) which parametrize (near infinity) the quasiconformal classes of 
non-renormalizable complex cubic polynomials. 
 Our techniques are based on the ideas introduced by Branner and Hubbard to study complex cubic polynomials.
\end{abstract}

\maketitle
%\tableofcontents
\markboth{ Puiseux series polynomial dynamics}{ \normalfont\scriptsize Jan Kiwi}
\section{Introduction}

The aim of this paper is to study  the dynamics of polynomials over
a complete algebraically closed ultrametric field $\L$ with a dense
valuation group and characteristic zero residual field.  In
particular, such a field $\L$ is an extension of $\Q$ and induces the
trivial valuation in $\Q$.  Our interest arises from  the extensive research on the dynamics
of rational functions over $\C$ and the recent one over $\C_p$.  Fields such as $\L$ seem to be a
natural dynamical space to explore the interplay between
non-archimedean and complex dynamics. The focus of this paper is
on cubic polynomials. We will show that the techniques developed
by Branner and Hubbard~\cite{branner-88,branner-92} to study complex
cubic polynomials merge with some basic ideas from $p$-adic dynamics
to give a complete picture of the dynamical behavior and the parameter
space of cubic polynomials in $\L$.  Moreover, we will show that
the dynamics of a family of cubic polynomials acting on $\S$, a
smallest ultrametric field with the above properties, is intimately
related to the structure of the parameter space of complex cubic
polynomials near infinity.  In particular, we show that some elements
of $\S$ naturally correspond to the Fourier series of analytic almost
periodic functions (in the sense of Bohr) which parametrize (near
infinity) the quasiconformal classes of non-renormalizable complex
cubic polynomials.

\medskip
Let us now describe the context and statements of our main results in more detail. 
For this let us fix a field $\L$ with the above properties. 
For the sake of simplicity we restrict  to polynomial dynamics.
Given a polynomial $\vphi$ with coefficients in one of the fields $\C$, $\C_p$ or $\L$,
the set of non-escaping points is the filled Julia set $\filled (\vphi)$ and 
its boundary $\julia (\vphi)$ is called the Julia set of $\vphi$.  The complement of the Julia set is the Fatou set 
$F(\vphi)$ (see~\cite{milnor-99} and Chapter 6 in~\cite{rivera-00}). 

\smallskip
In complex polynomial dynamics it is useful to study the connected components of the filled Julia set.
Non-archimedean fields are totally disconnected and, following Rivera~\cite{rivera-00},  the analogue discussion requires to replace 
the definition of connected  component by the weaker notion of infraconnected component (see Subsection~\ref{infra-ss}). 
In~\cite{benedetto-02},  Benedetto gives examples of  polynomials $\vphi$  in $\C_p$ 
which have  non-trivial wandering infraconnected components in $\filled (\vphi)$. That is, 
an infraconnected component $C$ which is not a singleton and such that 
the iterated images of $C$ under $\vphi$ are pairwise disjoint.
These examples rely on the presence of a fixed infraconnected component of the filled Julia set
with ``inseparable reduction'' (compare with~\cite{fernandez-04}). 
A non-trivial infraconnected component of $\filled (\vphi)$ is automatically a ball contained in the Fatou set 
(see Subsection~\ref{infra-ss}).
Thus a main difference with complex dynamics appears, since it follows from  a Theorem by Sullivan  (e.g., see~\cite{milnor-99})
that  a wandering 
connected component of a complex filled Julia set is contained in the Julia set.
The examples of $p$--adic polynomials with non-trivial wandering infraconnected components show that the analogue
statement of Sullivan's Theorem is false in $p$-adic dynamics.
Since the residual field of $\L$ has characteristic zero 
there are no  inseparable components 
 and we show that the situation is similar to that in $\C$. In fact, to prove the theorem below
we adapt  Branner and Hubbard's techniques used in~\cite{branner-92} to show  that the  filled Julia set of any complex cubic
polynomial does not have  non-trivial wandering components.

\begin{introtheorem}
  \label{wandering-it}
  For any cubic polynomial $\varphi \in \L [\z]$, every  infraconnected component of the filled Julia set
is either a singleton or  eventually periodic. 
\end{introtheorem}

The field $\L$ is not locally compact. Nevertheless some polynomial Julia sets in $\L$ are non-empty and compact. 
In particular, we show that
if $\vphi \in \L[\z]$ is a polynomial of degree $d \geq 2$  with all its critical points escaping, then the
Julia set $\julia (\vphi)$ is a non-empty compact set and the dynamics over it is topologically
conjugated to the one--sided shift on $d$ symbols (Theorem~\ref{cantor-th}). 
This is the analogue of a classical result in complex  dynamics (e.g., see Theorem 9.9 in~\cite{blanchard-84}). 
In $p$-adic dynamics, a similar statement  is false.

The question of characterizing compactness in $p$-adic dynamics has been addressed by Bezivin 
in~\cite{bezivin-04} where it is shown
that an obstruction for a Julia set to be compact is the existence of non-repelling cycles. 
For cubic polynomials in $\L$ we show that the absence of non-repelling cycles is in fact equivalent to the compactness
of the Julia set. Also, we show where in the parameter space $\cubicl$ of cubic
 polynomials those with non-empty and compact Julia set may be found.
More precisely, we work in the parameter space  $\cubicl$ formed by the polynomials of the form:
$$\vphiab (\z) = \z^3 - 3 \alpha^2 \z + \beta$$
where $(\alpha, \beta ) \in \L^2$.  Thus parameter space  is naturally identified with $\L^2$. As suggested by 
the previous paragraph, we say that the shift locus
$\shiftcubicl$ is the subset of $\cubicl$ formed by the polynomials with all their critical points escaping.
Now our characterization of compact Julia sets reads as follows:

\begin{introtheorem}
  \label{compact-it}
  Let $\varphi \in \cubicl$ be a cubic polynomial. Then the following are equivalent:
 
(i) The Julia set $\julia (\varphi)$ is a compact non-empty set.

(ii) $\julia (\varphi)= \filled (\varphi)$.

(iii) All the cycles of $\varphi$ are repelling.

(iv) $\varphi$ is in the closure of the shift locus $\shiftcubicl$.
\end{introtheorem}

The previous results rely on the detailed study of both the dynamical and parameter space of cubic polynomials
contained in sections~\ref{dynamical-s} and~\ref{parameter-s}. One  interesting consequence of our description of
parameter space we show that the subset $\cH$ of $\cubicl$ formed by all the polynomials whose Julia set is critical point 
free is open and dense (see Corollary~\ref{hyperbolic-c}). According to Benedetto~\cite{benedetto-01}, 
polynomials in $\cH$ exhibit some sort of hyperbolicity over their Julia set. Another consequence of our description
is the existence of cubic polynomials with coefficients in $\laut$ (the formal Laurent series with coefficients in an algebraic closure 
of $\Q$) which have a recurrent and non-periodic critical point (Corollary~\ref{recurrent-c}).

\medskip
The situation for quadratic polynomials in $\L$  is rather trivial. If the critical point of a quadratic polynomial is non-escaping, then
the filled Julia set is a closed ball and all the cycles are non-repelling. 
Otherwise, the filled Julia set is a Cantor set and all the cycles are repelling.
The situation for polynomials of degrees greater than $3$ is more subtle. Nevertheless, for polynomials of any degree ($\geq 2$)
with coeffiecients in a  smallest field
$\S$,  we conjecture
that the statement of the theorems above still hold. 

\bigskip
Let us now outline the results which establish a connection between dynamics over $\S$ and $\C$.
In complex cubic polynomial dynamics we work in the parameter space $\cubicc$ formed by all the polynomials
of the form:
$$g_{a,b} (z) = z^3 - 3 a^2 z + b$$
where $(a,b) \in \C^2$. 
This parameter space is naturally identified with $\C^2$. Note that the critical points of $g_{a,b}$ are $\pm a$.
If the critical point  $+a$ lies in $\filled (g_{a,b})$ we denote by $C_{a,b}(+a)$ the connected component of 
$\filled (g_{a,b})$ that contains $+a$.
We are interested in the set
$$A_\C := \{ (a,b) \in \C^2 \,\,/\, +a \in \filled (g_{a,b}) \mbox{ and } C_{a,b} (+a) \mbox{ is not periodic} \}.$$
Equivalently, $A_\C$ consists of  the polynomials such that  the critical point $-a$ escapes and
that  are not renormalizable about  the critical point $+a$. 

\smallskip
The field of formal Puiseux series $\pui$ is an algebraic closure of the field formal Laurent series $\laut$ with coefficients
in the algebraic closure $\Qa$ of $\Q$. We always regard $\Qa$ as a subset of $\C$. The field $\S$ is, up to isomorphism,
 the completion 
of $\pui$ with respect to an appropriate valuation (see Subsection~\ref{example-ss}). So for the rest of this paper
$\S$ will denote the completion of $\pui$. The elements of $\S$ may be identified
with series of the form $$ \zeta = \sum_{\lambda \in \Q} a_\lambda t^\lambda$$
where $a_\lambda \in \Qa$ and 
the set $\{ \lambda \, / \, a_\lambda \neq 0 \}$
is discrete and bounded below in $\R$.  An automorphism $\sigma$ of $\S$ over $\laut$ will play an special role in our work.
More precisely, we let $\sigma$ be the unique automorphism of $\S$ such that $\sigma (t^{1/m}) = \mathe^{2 \pi i/m} t^{1/m}$
for all $m \in \N$. 

\smallskip
Consider the family of cubic polynomials in $\S$:
$$\vphib (\z) = \z^3 - 3 t^{-2} \z + \beta$$
where $\beta \in \S$. 
Here the critical points are $\pm t^{-1}$ and 
if $t^{-1} \in \filled (\vphib)$ we denote by $IC_\beta (t^{-1})$ the infraconnected component of 
$\filled (\vphib)$ that contains $t^{-1}$. Now the analogue of $A_\C$ is 
$$A_\S := \{ \beta \in \S \,\,/\, t^{-1} \in \filled (\vphib) \mbox{ and } IC_\beta (t^{-1})  \mbox{ is not periodic} \}.$$

To state the correspondence between $A_\S$ and $A_\C$ we also need to introduce,
for $\epsilon >0$, the infinite strip 
$$\hepsilon := \{ T \in \C \,\,/\, \operatorname{Im} T > \frac{- \log \epsilon}{2 \pi} \}.$$

\begin{introtheorem}
  \label{complex-it}
  There exists $\epsilon >0$ such that for all $\beta = \sum_{\lambda \in \Lambda} a_\lambda t^\lambda$
in $A_\S$ the series $$ \sum_{\lambda \in \Lambda} a_\lambda  \mathe^{2 \pi i T \lambda}$$
is the Fourier series of an analytic almost periodic function $b_\beta : \hepsilon \rightarrow \C$.
Moreover, 
$$  \begin{array}{rccc}
  \tilde{\Phi} : & \hepsilon \times A_\S & \rightarrow &   A_\C \cap \{ | a | > 1/\epsilon \} \\
                      & (T, \beta)                    & \mapsto      & (\mathe^{-2 \pi i T}, b_\beta (T))
  \end{array}$$
is a well defined surjective map which is continuous in $(T,\beta)$ and holomorphic in $T$.
Furthermore, $\tilde{\Phi}$ projects to a homeomorphism:
$$  
  {\Phi} :  \hepsilon \times A_\S  / \sim   \rightarrow   A_\C \cap \{ | a | > 1/\epsilon \} 
$$
where $\sim$ is the smallest equivalence relation that identifies $(T-1, \beta)$ with $(T, \sigma (\beta))$.
\end{introtheorem}

For a short summary regarding almost periodic functions see Subsection~\ref{almost-ss}.

\smallskip
As an immediate consequence we recover (near infinity) a result by Branner and Hubbard which says that the local structure of 
$A_\C$ is that of a totally disconnected set cross a disk.
However, in the proof of the previous theorem we use two important ideas introduced by Branner and Hubbard:
the wringing construction and the tableaux. Thus, the theorem above and some direct consequences cannot be regarded
as independent from Branner and Hubbard's work. From the complex dynamics viewpoint, 
the novelty is the natural  parametrization  for $A_\C$ near infinity which is a manifestation of the interplay between the
dynamics over $\S$ and $\C$.

\smallskip
The proof of Theorem~\ref{complex-it} relies in our description of the parameter space $\cubicl$ achieved in Section~\ref{parameter-s}
as well as some complex dynamics techniques.

\bigskip
Let us now outline the structure of the paper:

\smallskip
Section~\ref{prelim-s} consists of some preliminaries. 
After giving a short discussion about the smallest field $\S$ we summarize the basic properties of
the action of polynomials on $\L$. Then we introduce ``affine partitions'' of a closed ball (which in the language of
~\cite{escassut-03} are the ``classes'' of a closed ball) and show that polynomials act on affine partitions. 
We continue with some dynamical aspects of polynomials in $\L$ such as their Fatou and Julia sets, and infraconnected components
of their filled Julia set. Simultaneously we discuss the basic combinatorial structure of the dynamical space of polynomials
in $\L$ given by balls and annuli of level $n$.

\smallskip
Section~\ref{shift-s} is devoted to the proof of Theorem~\ref{cantor-th} which describes the Julia set of polynomials with all their critical points escaping. 

\smallskip
Section~\ref{dynamical-s} contains a detailed study of the geometry of the filled Julia set of cubic polynomial 
with one critical point non-escaping and the other one in the basin of infinity. 
This study is based on Branner and Hubbard's ideas for organizing the relevant combinatorial information by introducing 
marked grids and tableaux. This section concludes with the proof of a stronger version of Theorem~\ref{wandering-it},
a corollary which establishes the equivalence of (i) through (iii) as stated in Theorem~\ref{compact-it} and Proposition~\ref{top-ent-p} regarding the topological entropy of cubic polynomials.

\smallskip
In Section~\ref{parameter-s} we give a detailed description of the parameter space $\cubicl$. At the end of this section
the reader may find  the proof of Theorem~\ref{compact-it}. 

\smallskip
In Section~\ref{complex-s} we prove Theorem~\ref{complex-it}. 
Here the key to pass from $\S$ to $\C$ are the Puiseux series of the ends of periodic curves in $\cubicc$.

\bigskip
\noindent
{\bf Acknowledgements.} I would like to thank Juan Rivera Letelier for
introducing me into non-archimedean dynamics. The influence of
conversations with him are scattered all over this work. 
I am grateful to John Milnor since an important motivation for this paper was my interest on understanding
the Puiseux-Laurent series of some ends of periodic curves computed by
him in a short note that he kindly gave to me some years ago.  
I thank  Manuel Elgueta and Alejandro Ramirez for useful conversations regarding almost periodic functions.  
I would like to thank the Royal Society for funding my trip to England
during April 2004. The conversations held with Adam Epstein and Mary
Rees during this trip were very helpful to organize the exposition of these results.

\section{Preliminaries}
\label{prelim-s}
Throughout the paper $\L$ will denote an algebraically closed field endowed with a non--archimedean valuation $\lval \cdot \rval$  such that $\L$ is complete, the residual field $\tL$ has characteristic zero and the valuation group 
$\valgL$ is dense in $(0, +\infty)$. 

  \begin{remark}
    {\em The residual field $\tL$ has characteristic zero if and only if $\L \supset \Q$ and $\lval x \rval = 1$ for all 
$x \in \Q^* = \Q \setminus \{ 0 \}$. }
  \end{remark}

\subsection{Example}
\label{example-ss}
Let  $\laut$ be the field of formal Laurent series in $t$ with coefficients in $\Qa \subset \C$ where $\Qa$ is the algebraic closure of $\Q$.
Given  a non-zero Laurent series $$\zeta = \sum_{j \geq j_0} a_j t^j \in \laut$$
define the {\bf order} of $\z$ by 
$$\order(\zeta) := \min \{ j \, / \, a_j \neq 0 \}$$
and consider the non-archimedean valuation in $\laut$ given by
$$\lval \zeta \rval := \be^{\order(\zeta)} \mbox{ where } \be = \mathe^{-1}.$$

The {\bf field of formal Puiseux series} with coefficients in $\Qa$,  denoted  $\pui$,
is  the algebraic closure of $\laut$ (e.g., see page 17 in~\cite{casas-00}).
The elements of $\pui$ may be identified with
the Laurent series in $t^{1/m}$ for some $m \in \N$. That is, 
for any  $\zeta \in \pui$ there exists 
$m \in \N$ such that  $$\zeta = \sum_{j \geq j_0} a_j t^{j/m} \in \lautm.$$
The unique extension of $\lval \cdot \rval$ from $\laut$ to $\pui$ is given by
$$\lval \zeta \rval = \mathrm{e}^{-\order(\zeta)}$$
where
$$\order(\zeta) = \frac{\min \{j \, / \, a_j  \neq 0 \}}{m}$$
provided that $\z \neq 0$.
The valuation group of $\pui$ is $\mathe^\Q$.

\medskip
We denote by  $\S$ the completion of $\pui$. 
The elements of $\S$ may be identified with the series
$$\zeta = \sum_{\lambda \in \Q} a_\lambda t^\lambda$$
where $a_\lambda \in \Qa$ and 
the set $\{ \lambda \, / \, a_\lambda \neq 0 \}$
is discrete and bounded below in $\R$.  
Moreover, $\lval \z \rval = \be^{\order(\z)}$ where 
$\order(\z) = \min \{ \lambda \in \Q \,\,/ \, a_\lambda \neq 0 \}$ if $\z \neq 0$. 
Therefore, $\valgS$ is also $\mathrm{e}^\Q$.
Since $\S$ is the completion of an algebraically closed field 
we have that $\S$ is also algebraically closed (e.g., see~\cite{cassels-85}). 

\medskip
The field $\L$ contains a copy of $\S$. In fact, given $\zeta \in \L$ such that $0 < \lval \z \rval <1$ it is not difficult to show that 
$$    \begin{array}{rccl}
\iota : & \laut & \rightarrow & \L \\
                 &  \sum_{i \geq i_0} a_i t^i        & \mapsto      & \sum_{i \geq i_0} a_i \z^i   
    \end{array}$$
extends to a  monomorphism $\iota : \S \rightarrow \L$. Therefore $\S$ is the smallest algebraically closed complete
field such that the residual field $\tS= \Qa$ has characteristic zero and $\valgS = \mathe^\Q$ is dense in $(0,+\infty)$.

\subsection{Polynomial maps in $\L$}

In this subsection we summarize some basic properties of polynomial maps in $\L$. 
Although most of these properties also hold for the larger class of holomorphic maps we only state them for polynomials in order to keep the exposition as simple as possible.

\medskip
For $r \in \valgL$ and $\znot \in \L$ we say that 
$$B^+_r(\znot) := \{ \z \in \L \, / \, \lval \z - \znot \rval \leq r \}$$
is a {\bf closed ball} and 
$$B_r(\znot) := \{ \z \in \L \, / \, \lval \z   - \znot  \rval  < r \}$$
is an {\bf open ball}.
If $r \notin \valgL$, then $B^+_r(\znot) =B_r(\znot)$ is an {\bf irrational ball}.
The reader should be aware that, in despite these names, topologically speaking every ball is  open and closed.

\smallskip
Consider $\varphi(\z) \in \L[\z]$
and $\z_0 \in \L$. The largest integer $d_0$ such that
$(\z- \z_0)^{d_0}$ divides $\varphi(\z)- \varphi(\z_0)$ 
is called {\bf the degree of $\varphi$ at $\znot$} and
denoted by $\deg_{\znot} (\varphi) $.
If the degree of $\varphi$ at  $\znot$  exceeds $1$, we
say that $\znot$ is a {\bf critical point of multiplicity $\mult_\varphi (\znot) := \deg_\znot (\varphi) -1$}.

Suppose that $\vphi(B) = B^\prime$ where $B$ is some subset of $\L$.
If there exists an integer $d_B \geq 1$ such that 
$$d_B = \sum_{\{ \z \in B \, / \, \varphi(\z) = \z^\prime \}}
\deg_\z (\varphi)$$
for all $\z^\prime \in B^\prime$,
then we say that $\varphi: B \rightarrow B^\prime$ has degree 
$d_B = \deg_B (\varphi)$. 

\medskip
Polynomials map balls onto balls (see~\cite{rivera-00} page 29):

\begin{proposition}
\label{balls-to-balls-p}
Let $\varphi(\z) \in \L [\z]$ be a polynomial of degree $\deg (\varphi)$.
Consider a closed (resp. open, irrational) ball $B \subset \L$. 
Then the following hold:
\begin{enumerate}
  \item[(i)]
$\varphi(B)$ is a closed (resp. open, irrational) ball.

\item[(ii)]
$\varphi : B \rightarrow \varphi(B)$ has a well defined degree
$\deg_B (\varphi)$.

\item[(iii)]
$\varphi^{-1}(B)$ is a disjoint union
of closed (resp. open, irrational) balls $B_1,  \dots,  B_k$ 
such that $$\sum \deg_{B_i} (\varphi) = deg(\varphi).$$

\item[(iv)]
$$\deg_{B} (\varphi) -1 = 
\sum_{\zeta \in B} (\deg_\zeta (\varphi) -1) = \sum_{\z \in \mathrm{Crit}(\varphi)\cap B} \mult (\z)$$
where  $\mathrm{Crit}(\varphi)$ is the set formed by the critical points of $\varphi$.
\end{enumerate}
\end{proposition}

Statement~(iv)  makes a substantial difference between dynamics over fields with characteristic zero residual 
fields  (e.g., $\L$) and dynamics over  fields with residual fields with  non-vanishing characteristic (e.g., $\CC_p$).

\medskip
\noindent
{\bf Sketch of the Proof.} 
Statements (i)--(iii) follow by inspection of the Newton polygon of $\varphi$.
We refer the reader to~\cite{cassels-85} for background on Newton  polygons 
and~\cite{rivera-00} for a proof of (i)--(iii) in the context of $p$-adic holomorphic functions 
that applies without modifications to our context.
Statement (iv) follows from a simple observation. Without loss of generality we may assume that: $B$ and $\varphi(B)$
are balls which contain the origin,  $\varphi(0) \neq 0$, and $\varphi^\prime (0) \neq 0$. Since natural numbers
have valuation $1$, the Newton polygon of $\varphi$ translated to the left by $1$ and restricted to the right half plane is the Newton
polygon of $\varphi^\prime$. 
It then follows that the number of zeros of $\varphi$ in $B$ minus $1$ coincides with the number of
zeros of $\varphi^\prime$ in $B$.
\hfill $\Box$

\medskip
We say that $A \subset \L$ is an {\bf annulus} if 
$$A = \{ \z \in \L \, / \, \log \lval \zeta - \znot \rval \in I \}$$
for some $\znot \in \L$ and some interval $I \subset (-\infty, \infty)$.
We say that $A$ is an open (resp. closed) annulus if $I$ is open (resp. closed) interval.
The length of $I$ is by definition the  {\bf modulus of $A$}, denoted $\mdl A$.
The next proposition describes how the modulus of an annulus changes
under the action of a polynomial $\varphi$.

\begin{proposition}
\label{annulus-p}
If $A, A^\prime$ are annuli and $\varphi(\z) \in \L [\z]$ is
such that $\varphi(A) = A^\prime$, then $\varphi: A \rightarrow A^\prime$
has a well defined degree $\deg_A ( \varphi)$ and 
$$\deg_A (\varphi) \cdot \mdl A = \mdl A^\prime.$$
\end{proposition}

The statement of Lemma 5.3 in~\cite{rivera-03}  is the same than the one of 
 the previous proposition but  in the context of holomorphic functions
in $\C_p$. Rivera's proof  applies to our setting as well.

\medskip
We will also need the following version of Schwarz's Lemma (see~\cite{rivera-00})

\begin{lemma}
\label{schwarz-l}
  Consider $\varphi(\z) \in \L [\z]$. Assume that
$\varphi (B_0 ) \subset B_1$ where $B_i$ is a ball of radius
$r_i$ for $i=0,1$.
Then, for all $\z_1, \z_2 \in B_0$:
\begin{eqnarray}
\lval \varphi(\zone) - \varphi(\z_2) \rval & \leq & 
\frac{r_1}{r_0} \lval \zone - \z_2 \rval  \label{sc1-e} \\
\lval \varphi^\prime (\zone) \rval & \leq &\frac{r_1}{r_0}.  \label{sc2-e}
\end{eqnarray}
Moreover, equality holds at some $\z_1, \z_2$ in (\ref{sc1-e}) or at some $\z_1$ in (\ref{sc2-e}) if and only if
equality holds for all $\z_1, \z_2$ in (\ref{sc1-e}) and all $\z_1$ in (\ref{sc2-e}).

\end{lemma}

The next lemma will be useful to count the number of fixed points inside a given closed ball.

\begin{lemma}
  \label{fixed-l}
Let $\varphi \in \L [\z]$. Let $B$ and $B^\prime$ be closed balls such that $B^\prime = \varphi(B) \supset B$.
Denote by $|\operatorname{Fix}_B (\varphi) |$ the number of fixed points of $\varphi$ in $B$ counting multiplicities.
If $\deg_B(\vphi) >1$ or $B \subsetneq B^\prime$, then
$$|\operatorname{Fix}_B (\varphi) | =  \deg_B (\varphi).$$
\end{lemma}

\noindent
{\bf Proof.}
Without loss of generality $B = B^+_1(0)$.

In the case that there exists $\znot \in B$ such that $\lval \varphi^\prime (\znot) \rval >1$, 
after conjugation  by $\z \mapsto \z - \znot$,  we may assume that $\znot =0$. 
It follows that  the Newton polygons for $\varphi(\z)$ and $\varphi(\z) - \z$ coincide and therefore
$|\operatorname{Fix}_B (\varphi) | =  \deg_B (\varphi)$.

For the case in which $\lval \varphi^\prime (\z) \rval \leq 1$ for all $\z \in B$  we  write
$$\varphi(\z) = \alpha_0 + \alpha_1 \z + \cdots + \alpha_n \z^n$$
and observe that $\varphi(B) =B$ and that $\lval \alpha_k \rval \leq 1$ for all $k$. Also,  the number of zeros of $\varphi$ in $B$ is 
$\deg_B (\varphi)$ and  coincides with  the maximal 
index $k$ for which $\lval \alpha_k \rval =1$. Since the coefficient of 
$\z^k $ in $\varphi(\z)-\z$ coincides with $\alpha_k$ for all $k \neq 1$,
if $\deg_B (\varphi) >1$, then $\vphi(\z)-\z$ has exactly  $\deg_B(\vphi)$ zeros in $B$ (counting multiplicities).
\hfill $\Box$

\subsection{Affine Partitions}
\label{affine-ss}
In the study of iterations of rational functions on $p$-adic fields it is useful to consider their action on projective systems
(see~\cite{rivera-00}).
For polynomials the situation is simpler and we will just need to consider affine partitions (compare with the ``classes'' 
of a ball in~\cite{escassut-03}). 

By definition, the {\bf canonical affine partition} $$\cP_c := \{ B_1 (c) \,\,/\, c \in \tL \}$$  is the collection of equivalence classes
of the ring  $\cO_\L = B^+_1(0)$ modulo the ideal $\cM_\L = B_1(0)$. 
The {\bf affine partition} $\cP_{B_0}$ associated to a  closed  ball $B_0$ is:
$$\{ h^{-1} (B) \,\, /\, B \in \cP_c \}$$
where $h: \L \rightarrow \L$ is an affine map such that $h(B_0) = B^+_1(0)$.
Affine partitions are parametrized by the residual field $\tL$ and the parametrization is unique up to 
$\tL$--affine maps. Therefore, affine partitions inherit the affine structure of $\A^1(\tL)$.

\begin{proposition}
\label{affine-p}
  Let $\varphi : \L \rightarrow \L$ be a polynomial. Given a closed ball $B_0 \subset \L$ let $B_1 = \varphi (B_0)$.
Denote by $\cP_0$ and $\cP_1$ the associated affine partitions. Then:
\begin{enumerate}
  \item[(i)]
 There is a well defined  induced action on the affine partitions given by:
$$    \begin{array}{rccl}
\varphi_* : & \cP_0 & \rightarrow & \cP_1 \\
                 &  B       & \mapsto      & \varphi(B)
    \end{array}
$$
Moreover, $\varphi_*$ is a polynomial from  the affine structure of $\cP_0$ to that of $\cP_1$.
\item[(ii)]
  $\deg (\varphi_*) = \deg_{B_0} (\varphi)$.
\item[(iii)]
  $ \deg_{B} (\varphi_*) = \deg_B(\varphi) $  for all $B \in \cP_0$.
\end{enumerate}
\end{proposition}

\noindent
{\bf Proof.}
We first apply  an affine change of coordinates in the domain and the range so that  $B_0=B_1 = B^+_1(0)$. Hence
$\varphi(\z) = \alpha_0 + \cdots + \alpha_n \z^n$ with $\lval \alpha_k \rval \leq 1$. Now let
$\pi : B^+_1 (0) \rightarrow \tL$ be the quotient map and for $\z \in B^+_1 (0) $ let  $\tz = \pi(\z)$.
It follows that $\tvarphi(\tz) = \talpha_0 + \cdots + \talpha_{d} \tz^d $ is such that $\pi \circ \varphi = \tvarphi \circ \pi$ where
 $d = \deg_{B_0} (\varphi)$. Thus $\varphi_{*}(\pi^{-1}(\tz)) = \pi^{-1} (\tvarphi(\tz))$ and $\varphi_*$, in these coordinates,
becomes $\tvarphi$. From where (i) and (ii) easily follow. 

For (iii), without loss of generality we may assume that $B=\varphi(B) = \pi^{-1}(0)$. Under this assumption
$\tvarphi(\tz) = \talpha_j \tz^j + O(\tz^{j+1})$. It follows that $j$ is the smallest index such that $\lval \alpha_j \rval =1$.
Looking at the Newton polygon of $\varphi$ we conclude that $j$ is the degree of $\varphi: B \rightarrow B$ and (iii) follows.
\hfill $\Box$

\subsection{Fatou and Julia Sets}
The {\bf chordal metric} on $\L$  is  defined by
$$\chordal (\znot, \zone) := \frac{ \lval \znot - \zone \rval }{\max( \lval \znot \rval, 1 ) \cdot \max( \lval \zone \rval, 1 )}$$
for $\znot, \zone \in \L$.

The {\bf Fatou set} $\fatou (\varphi)$ is  the set formed by all $\znot \in \L$ for which there exists a neighborhood
$U$ of $\znot$ such that $\{ \varphi^n : U \rightarrow \L \}_{n \geq 1}$ is an uniformly Lipschitz collection of functions
(with respect to the chordal metric).
The {\bf Julia set} $\julia (\varphi)$ is the complement of the Fatou set. That is, $\julia (\varphi) := \L \setminus \fatou (\varphi)$.

According to Hsia~\cite{hsia-00}, a sufficient condition for the collection 
$$\{ \varphi^n : U \rightarrow \L \}_{n \geq 1}$$ 
to be uniformly Lipschitz is that 
$$\L \setminus \underset{n \geq 1}{\cup} \varphi^n (U) \neq \emptyset.$$

Given $\varphi \in \L [ \z]$, in analogy with complex polynomial dynamics, the {\bf filled Julia set} is defined by
$$\filled (\varphi) := \{ \z \in \L \,\,/\, \lval \varphi^n (\z) \rval \not\rightarrow \infty \}.$$
That is, the filled Julia set is the complement of the basin of $\infty$.
Although $\julia(\varphi)$ might be empty (e.g., $J(\z^2)=\emptyset$)  the filled Julia set
 $\filled (\varphi)$  is always non-empty since $\filled(\varphi)$
 contains the  periodic points of $\varphi$.
According to Proposition 6.2 in~\cite{rivera-00} a polynomial Julia set can be characterized as follows:
\begin{eqnarray*}
\julia (\varphi) & = & \partial \filled (\varphi)  \\
& = & \{ \z \in \L \,\, / \, \underset{n \geq 1 }{\cup} \varphi^n (U) = \L \mbox{ for all open sets } U 
\mbox{ with  } \z \in U \}.
\end{eqnarray*}

\subsection{Dynamical balls and infraconnected components of a filled Julia set}
\label{infra-ss}
Throughout this subsection, let $\varphi$ be a degree $d > 1$ polynomial of the form:
$$\varphi (\z) = \alpha_0 + \alpha_1 \z + \cdots + \alpha_d \z^d \in \L[\z]$$
where $\alpha_d \neq 0$.
Following Section~6.1 of~\cite{rivera-00}, let 
$$R_\varphi := \max \left( \left| \frac{\alpha_i}{\alpha_d} \right|^{\frac{1}{d-i}}_o, \left| \frac{1}{\alpha_d} \right|^{\frac{1}{d-1}}_o \right).$$
Then it is easy to check that 
$K(\varphi) \subset \varphi^{-1}(\{ \lval \z \rval \leq R_\varphi \}) \subset \{ \lval \z \rval \leq R_\varphi \}$ 
and  
$$K(\varphi) = \{ \z \in \L \,\, /\, \lval \varphi^n (\z) \rval \leq R_\varphi \mbox{ for all } n \geq 1\}.$$

\begin{lemma}
  \label{diam-l}
  Given a polynomial $\psi \in \L [\z]$ there exists another polynomial $\varphi \in \L[\z]$ affine conjugate to $\psi$ such 
that $R_\varphi = \operatorname{diam} \filled(\varphi) = \sup \{ \lval \zone - \z_2 \rval \,\,/\, \zone, \z_2 \in \filled(\vphi) \}$.
\end{lemma}

\noindent
{\bf Proof.}
After an affine conjugacy $\psi$ becomes $\varphi (\z) = \alpha_1 \z + \cdots + \alpha_{d-1} \z^{d-1} + \z^d$.
Note that $$R_{\varphi} =\max ( \{\lval \alpha_j \rval^{\frac{1}{d-j}} \,\, / \, 1 \leq j <d \} \cup \{1 \}).$$ 
Hence, if $R_{\varphi}=1$, then $\filled (\varphi) = B^+_1(0)$. Otherwise, $R_{\varphi} >1$ and
from the Newton polygon of $\varphi$ we deduce that there exists $\znot$ such that $\lval \znot \rval = R_\varphi$
and $\varphi (\znot)=0 \in \filled (\varphi)$.
\hfill $\Box$ 

\begin{definition}
\label{dynaball-d}
 We say that $\db_0 = B^+_{R_\varphi} (0)$ is the {\bf dynamical ball of level $0$} of $\varphi$.
The set
$\varphi^{-n}(\db)$ is the union of finitely many pairwise disjoint closed balls which we call {\bf level $n$ dynamical balls}.
\end{definition}

Later we will introduce ``parameter'' balls of level $n$. Often, when clear from the context, a dynamical ball
will be  simply called a ball.

\medskip
Observe that each ball of level $n > 0$ is contained in exactly one of level $n-1$ and maps onto a level $n-1$ ball.

\begin{definition}
  \label{dynanest-d}
A {\bf dynamical end}
 $\dn$ is a sequence $\{ \db_n \}_{n \geq 0}$ such that $\db_n$ is a ball of level $n$ and $\db_{n+1} \subset \db_n$
for all $n$.
\end{definition}

The map $\varphi$ acts on ends. In fact, given an end  $E = \{ \db_n (E) \}$ let $\db_n (\varphi (E)) = \varphi (\db_{n+1} (E))$
for all $n \geq 0$. It follows that $\varphi (E) := \{ \db_n (\varphi(E)) \}$ is an end which we call 
{\bf the image of $E$ under $\varphi$}.

\medskip
Following Escassut~\cite{escassut-95} a subset $X$ of $\L$ is called {\bf infraconnected} if whenever
$ X \subset B_0 \cup B_1$ for some disjoint closed balls $B_0, B_1$, then $X \subset B_0$ or $X \subset B_1$.
An infraconnected component of $Y \subset \L$ is an equivalence class of the relation that identifies two 
points $\z_0, \z_1$ if there exists an infraconnected subset of $Y$ containing both $\z_0$ and $\z_1$.

\medskip
Proposition 6.8 in~\cite{rivera-00} reads as follows:

\begin{lemma}
  \label{infra-l}
  \begin{enumerate}
    \item[(i)]
If $\dn = \{ \db_n \}$ is an end, then $\cap \db_n$ is empty, or a singleton, or a closed ball, or an irrational ball.

\item[(ii)]
If $\z \in \filled (\varphi)$, then there exists a unique end $\dn (\z) = \{ \db_n (\z) \}$ such that $\z \in \cap \db_n (\z)$.
Moreover, the infraconnected component of $\filled (\varphi)$ which contains $\z$ is $\cap \db_n(\z)$.

\item[(iii)]
For any $\z \in \filled(\varphi)$, the infraconnected component of $\filled (\varphi)$ which contains $\z$ is
a singleton  if and only if $\z \in \julia(\varphi)$. 
  \end{enumerate}
\end{lemma}

A well known result in complex polynomial dynamics states that  the filled Julia set of a polynomial $f$ is connected
if and only if all the critical points of $f$ have bounded orbit (e.g., see Theorem 9.5 in~\cite{milnor-99}). We obtain 
a similar result for polynomial dynamics in $\L$:

\begin{corollary}
\label{connected-c}
  Let $\varphi \in \L [ \z]$ and denote by $\operatorname{Crit}(\varphi)$ the set of critical points of $\varphi$.
Then $\filled(\varphi)$ is infraconnected if and only if $\operatorname{Crit}(\varphi) \subset \filled (\varphi)$.
In this case, $\filled (\varphi)$ is a closed ball.
\end{corollary}

\noindent
{\bf Proof.} First suppose that $\operatorname{Crit}(\varphi) \subset \filled (\varphi)$. In view of Lemma~\ref{diam-l} we may 
assume the $R_\varphi = \diam \filled (\varphi)$. 
From Proposition~\ref{balls-to-balls-p} (iv) it follows that there
exists a unique level $1$ dynamical ball $\db_1$ which must coincide with the level $0$ ball $\db_0$ since
$R_\vphi = \diam \filled (\vphi)$ and $\filled (\vphi) \subset \db_1$.  Therefore $\filled (\vphi) = \db_0$.

If $\operatorname{Crit}(\varphi) \not\subset \filled (\varphi)$, then there 
exist a level with at least two disjoint balls, say $B_1$ and $B_2$. Each one of these balls
 $B_i$ contains a periodic point $\z_i$ because
there exists $k$ such that $\varphi^k(B_i) \supsetneq B_i$ (Lemma~\ref{fixed-l}). 
It follows that the infraconnected components $C(\z_0), C(\z_1)$ of $\filled(\varphi)$ 
containing $\z_0, \z_1$ (respectively) are distinct and therefore $\filled(\varphi)$ is not
infraconnected.
\hfill $\Box$

\medskip
Regarding compactness of $\julia (\varphi)$ we obtain the following result.(Compare with~\cite{bezivin-04}.) 

\begin{corollary}
  \label{compact-c}
Given $\varphi \in \L [\z]$ the following hold:

(i) If $\julia (\varphi)$ is compact and non-empty, then every infraconnected component of $\filled (\varphi)$ is a singleton.

(ii) If every infraconnected component of $\filled (\varphi)$ is a singleton, then all the cycles of $\varphi$ are repelling.
\end{corollary}

\noindent
{\bf Proof.}
For (i) we proceed by contradiction and suppose that $\julia(\varphi)$ is compact and non--empty and that there exists 
and end $E = \{ \db_n \}$ such that $C = \cap \db_n$ is a ball. Since $\julia (\varphi) = \partial \filled (\varphi)$, the Fatou set
contains $C$. Now let $\znot \in \julia (\varphi)$. For all $n \geq 0$ there exists $\z_n \in \db_n$ such that $\varphi^n (\z_n) = \z_0$ since
$\varphi^n(\db_n) = \db_0 \supset \julia (\varphi) \ni \znot$. Therefore, after passing to a convergent subsequence we obtain a limit point
 $\z \in J(\varphi) \cap C = \emptyset$ which is a contradiction.

For (ii), suppose that $\znot$ is a period $p$ periodic point. Then $\{ \znot \} = \cap \db_n (\znot)$ where $\db_n (\znot)$ 
is the level $n$ ball containing $\znot$. The orbit of $\znot$ does not contain critical points, for otherwise the infraconnected
component of $\znot$ in $\filled (\varphi)$ would contain points that are attracted to the cycle of $\znot$. Hence, for $n$ large
$\varphi^p : \db_n (\znot) \rightarrow \db_{n-p}(\znot)$ has degree $1$. By Schwarz Lemma~\ref{schwarz-l}, 
$\lval (\varphi^p)^\prime (\znot) \rval > 1$ and $\znot$ is repelling.
\hfill $\Box$

\subsection{Points and Annuli of level $n$}
Consider a polynomial  $\vphi \in \L [\z]$, an integer $n \in \N$ and a point $\z \in \vphi^{-n}(D_0)$ where $D_0 = B^+_{R_\varphi} (0)$ is the level $0$ ball of  $\vphi$. 
In this case we say that $\z$ is a {\bf level $n$ point}. Note that $\z$ is a level $k$ point for all $k \leq n$.
Also, $\z$ is contained in a unique level $n$ ball denoted $\db_n(\z)$. The radius of $\db_n (\z)$ will be denoted by $r_n (\z)$.
 
 Now let $\hr$ denote the radius of $\vphi(D_0)$. We say that 
$A_0 =B_{\hr}(0) \setminus D_0$ is  the {\bf level $0$ annulus of $\vphi$}.
For $n \in \N$, we say that the {\bf annulus of level $n$ around  $\z$} is $A_n(\z) = B_{r_{n-1} (\z)} (\z) \setminus D_n (\z)$
where $\z$ is a level $n$ point.
Note that:
$$\log \hr - \log r_n (\z) = \sum^n_{\ell=0} \mdl \da_\ell (\z).$$
Similarly if $E = \{ D_n \}$ is an end, then we denote by $D_n (E)$ the ball of level $n$ participating in $E$ and 
its radius will be denoted by $r_n (E)$. The {\bf level $n \geq 1$ annulus of $E$}
 is $A_n (E) = B_{r_{n-1}(E)}(\z) \setminus \db_n(E)$
where $\z$ is any point of $\db_n(E)$. 
Also, 
$$\log \hr - \log r_n (E) = \sum^n_{\ell=0} \mdl \da_\ell (E).$$

We omit the straightforward proof of the following result 
which shows the importance of studying the convergence of the sum of the moduli of annuli. 
\begin{lemma}
  \label{sum-l}
Let $\z \in \filled (\varphi)$ and $E$ be an end.
Then the following are equivalent:

(i) $r_n(\z) \rightarrow 0$ (resp. $r_n (E) \rightarrow 0$).

(ii) $\sum^{\infty}_{\ell =0} \mdl A_\ell (\z) = + \infty$ (resp. $\sum^{\infty}_{\ell =0} \mdl A_\ell (E) = + \infty$).

(iii) $\{ \z \} = \cap \db_n (\z)$ (resp.  $\cap \db_n (E)$ is a singleton).
\end{lemma}

\section{Polynomials with all critical points escaping}
\label{shift-s}
The  Julia set
of a degree $d>1$ polynomial $f : \CC \rightarrow \CC$ with all
its critical points escaping is a Cantor set. Moreover, the  dynamics over its Julia set $\julia (f)$
is topologically conjugated  to the one--sided shift on $d$ symbols and 
$f$ is uniformly expanding in a neighborhood of $\julia (f)$ (e.g., see Theorem 9.9 in~\cite{blanchard-84}).
The aim of this section is to prove the analogous result for polynomials acting on $\L$.

\begin{theorem}
\label{cantor-th}
  Let $\varphi : \L \rightarrow \L$ be a degree $d \geq 2$ polynomial
with all critical points escaping (i.e., $\omega \notin \filled (\varphi)$
for all critical points $\omega$). Then
   $\varphi : \filled (\varphi) \rightarrow \filled (\varphi)$
is topologically conjugated to the one--sided shift on $d$ symbols. Moreover, $\varphi$ is uniformly expanding
in a neighborhood of $\julia (\varphi)$. In particular, $\filled (\varphi)$ is a Cantor set and $\julia (\varphi) = \filled (\varphi)$. 
Furthermore, the intersection of every end is a singleton.
\end{theorem}

Before proving the theorem let us be more precise about the meaning of uniformly expanding maps.
(compare with Definition~3.1.~in~\cite{benedetto-01} and Definition 3 in~\cite{bezivin-04}).

\begin{definition}
\label{hyperbolicity-d}
We say that $\varphi$ is  {\bf uniformly expanding}  on a neighborhood $V$ of $\julia(\varphi)$
if there exist real numbers $0< c_1 < c_2$, a bounded continuous function $\tau : V \rightarrow [c_1, c_2]$ and
$\lambda >1$ such that
$$\tau(\varphi(\z)) \lval \varphi^\prime \rval \geq \lambda \tau(\z)$$
for all $z \in V$.
\end{definition}

\medskip
To prove the theorem we will need to label level $n$ balls:

\begin{lemma}
\label{label-l}
  Let $\varphi \in \L [\z]$ be a degree $d \geq 2$ polynomial.
  Then there exists a function $L$ that assigns to each level $n$
ball a subset of $\{ 1, \dots, d \}$ such that
if  $B$ and $B^\prime$ are dynamical balls of some level (not necessarily the same), then
the following hold:

\begin{enumerate}
\item[(i)]
  The cardinality of $L(B)$ is $\deg_B (\varphi)$.

\item[(ii)]
  If $B \subset B^\prime$, then $L(B) \subset L(B^\prime)$.

\item[(iii)]
  If $B \neq B^\prime$ and $\varphi(B) = \varphi ( B^\prime)$, then
$L(B) \cap L(B^\prime) = \emptyset$.
\end{enumerate}
\end{lemma}

\noindent
{\bf Proof.}
The  unique ball $\db_0$ of level $0$ is such that $\varphi^{-1}(\db_0) \subset
\db_0$, therefore  $\deg_{\db_0} (\varphi) = d$ and we are forced to 
label it $L(\db_0) = \{ 1, \dots, d \}$. 

We construct the labelling $L$ recursively.
Suppose that all level $n$ balls have been  labelled.
Consider a pair of level $n$ balls $\db_n, \db^\prime_n$ 
such that $\varphi(\db_n) \supset \db^\prime_n$.
We simultaneously label all the level $n+1$ balls $\db^1_{n+1}, \dots, \db^k_{n+1}$
contained in  $\db_n$ which map onto $\db^\prime_n$. In fact, since
 $$\sum \deg_{\db^j_{n+1}} (\varphi) = \deg_{\db_n} (\varphi)$$
it is sufficient to subdivide $L(\db_n)$ into $k$ sets 
$L(\db^1_{n+1}), \dots, L(\db^k_{n+1})$
of cardinalities $$\deg_{\db^1_{n+1}} (\varphi),\dots ,\deg_{\db^k_{n+1}} (\varphi).$$
Repeating this process for all pairs $\db_n, \db^\prime_n$ such that $\varphi(\db_n) \supset \db^\prime_n$ 
a definition of  $L$ on the level $n+1$ balls is achieved. 

Properties (i) through (iii) are easily checked for this labelling.
\hfill $\Box$

\medskip
\noindent
{\bf Proof of Theorem~\ref{cantor-th}.}
Consider a labelling as in the previous lemma. 
Let $N \geq 1$ be such that $\varphi^N(\omega) \notin \db_0$ for all
critical points $\omega$. That is, the level $N$ balls are  critical
point free and  therefore  each  level $N$ ball maps bijectively onto
one of  level $N-1$. In particular, $L(\db_N)$ is a singleton
for all level $N$ balls $\db_N$.

We first show that the intersection of every end $\dn = \{ \db_n \}$
is a singleton. For this we consider the metric
on $\varphi^{-(N-1)}(\db_0)$ defined by;
$$\rho(\z, \z^\prime) = \left\{ \begin{array}{ll} 
 {{\lval \z - \z^\prime \rval} \cdot {r_{N-1}(\z)}^{-1}} & \mbox{if } \db_{N-1} (\z)
= \db^\prime_{N-1}(\zp), \\
\lval \z  - \z^\prime \rval & \mbox{otherwise} \\
\end{array}\right. $$
where $r_{N}(\z)$ denotes  the radius of the level $N$ ball which contains $\z$.
Let 
$$\lambda = \min \left\{ \frac{r_{N-1} (\zeta)}{r_N(\zeta)} \, / \, \zeta \in \varphi^{-N} (\db_0) \right\} >1.$$
By Schwarz Lemma,
\begin{equation}
\label{expansion-e}
\rho(\varphi(\z), \varphi(\z^\prime)) \geq \lambda \rho (\z, \z^\prime)
\end{equation}
if $\db_N(\z) =\db_N(\z^\prime)$.
Moreover, 
$$r_{N-1}(\varphi(\z))^{-1} \lval \varphi^\prime (\z) \rval = r_N (\z)^{-1} \geq \lambda  r_{N-1}(\z)^{-1} $$
for all $\z \in \varphi^{-N}(\db_0)$.
In particular, $\varphi$ is uniformly expanding on the  neighborhood $\varphi^{-N}(\db_0)$ of $\filled(\varphi)$ 
taking  $\tau(\z)= r_{N-1}(\z)^{-1} $ in Definition~\ref{hyperbolicity-d}.

For $n \geq N-1$, 
let $$R_n = \max \{ \rho(\db_{n}) \, / \, \db_n \mbox{ ball of level } n \}$$ 
where $\rho(\db_n)$ is the $\rho$--radius of $\db_n$.
By (\ref{expansion-e}), $R_{N-1+k} \lambda^k \leq R_{N-1}$. 

It follows that if  $ \dn = \{ \db_n\}$ is an end, then 
$\rho(\db_n) \rightarrow 0$ as $n \rightarrow \infty$.
From the completeness of $\L$ we conclude that the  intersection of $\dn$ is
 a point. By Lemma~\ref{infra-l}, every infraconnected component
of $\filled(\varphi)$ is a point and $\julia (\varphi) = \filled (\varphi)$.

\smallskip
The labelling of the previous lemma determines  an itinerary for each end.
Namely, let $\dN$ denote the collection of all ends  and
$$\begin{array}{rccl}
\itin : & \dN & \rightarrow  & \{1, \dots, d\}^{\N \cup \{0 \}} \\
         & \{ \db_n \} & \rightarrow & (j_0,  j_1, \dots) \mbox{ if } \{ j_k \} = L(\varphi^k(\db_{N+k})).
\end{array} $$
It follows that the itinerary function is bijective. Moreover, 
for $\z \in \filled(\varphi)$, let $\dn (\zeta) = \{ \db_n (\zeta) \}$ be the end with intersection
$\{ \z \}$.
Then the map  $\z  \mapsto \itin (\dn (\zeta))$ gives the desired topological 
conjugacy between $\varphi: \filled (\varphi) \rightarrow \filled (\varphi)$ and 
the one--sided shift on $d$ symbols.
\hfill $\Box$

\section{Cubic polynomials: the dynamical space}
\label{dynamical-s}
From Corollary~\ref{connected-c} and Theorem~\ref{cantor-th} we conclude that
the filled Julia set of quadratic polynomials is either a closed ball or a Cantor set according to whether the unique critical 
point belongs to the filled Julia set or escapes to infinity.
For a cubic polynomial $\varphi \in \L[\z]$  we have three possibilities:

(i) All the critical points escape to infinity. In this case $\filled(\varphi)$ is a Cantor set (Theorem~\ref{cantor-th}).

(ii) All the critical points belong to $\filled (\varphi)$. Here $\filled (\varphi)$ is a closed ball (Corollary~\ref{connected-c}).

(iii) One critical point escapes to infinity and the other belongs to  $\filled(\varphi)$.

The aim of this section is to describe $\filled(\varphi)$ for polynomials as in (iii).

\subsection{Branner--Hubbard Tableaux}
\label{bh-ss}
{\bfseries\slshape 
  Our standing assumption for this subsection is that $\vphi$ is a cubic polynomial with two distinct critical points
$\omega^\pm$ such that $\omega^- \notin \filled (\vphi)$, $\vphi (\omep) \in \db_0$ and $R_\vphi = \diam \filled (\vphi)$ where $\db_0$ is the level $0$ ball of $\vphi$. The level $0$ annulus of $\vphi$ will be denoted $A_0$.}
(see Lemma~\ref{diam-l} and Definition~\ref{dynaball-d}).

\smallskip
According to Lemma~\ref{sum-l} to study the geometry of $\filled (\vphi)$ it is convenient to compute the moduli of
the annuli of level $n$, for all $n$.
The next pair of lemmas describe the behavior of level $n$ annuli under iterations:

\begin{lemma}
\label{level-zero-l}
Let $\varphi$ be a cubic polynomial with critical points $\omega^\pm$ such that $R_\varphi = \operatorname{diam} \filled(\varphi)$.
Suppose that $\varphi (\omega^+) \in D_0$  and $\omega^- \notin \filled(\varphi)$. 
Then the following hold:

(i) $\varphi(\omega^-) \notin D_0$.

(ii) There are exactly two level $1$ balls: $D_1 (\omega^+)$ and $D_1 (\gamma^+)$ where $\varphi(\gamma^+) = \varphi(\omega^+)$
and $\gamma^+ \neq \omega^+$.

(iii) The degree of $\varphi : D_1 (\omega^+) \rightarrow D_0$ is $2$ and the degree of $\varphi : D_1 (\gamma^+) \rightarrow D_0$
is $1$.

(iv) $\varphi (A_1 (\omega^+)) = A_0$ and $\varphi : A_1 (\omega^+) \rightarrow A_0$ is 
a degree $2$ map. Also,  $\varphi(A_1 (\gamma^+)) = A_0$ and 
$\varphi : A_1 (\gamma^+) \rightarrow A_0$ is a degree $1$ map.
\end{lemma}

Following Branner and Hubbard we say that $\gamma^+$ as in the lemma is the cocritical point of $\omega^+$.

\medskip
\noindent
{\bf Proof.}
For (i) we proceed by contradiction, if $\varphi(\omega^-) \in D_0$, then both critical points must be in the same level 
$1$ ball. Hence there would be a unique level $1$ ball which contains $\filled (\varphi)$ and has radius strictly smaller
than $R_\varphi$ which is a contradiction since $R_\varphi = \diam \filled(\varphi)$.

To prove statements (ii) and (iii) just observe that
from (i) it follows that  $\deg_ {D_1 (\omega^+)} (\varphi) =2$. 
Thus there exists another level $1$ ball which maps onto $D_0$ under $\varphi$ with degree $1$. 

For (iv),  note that $\varphi^{-1} (D_0) \cap B_{R_\varphi} (\omega^+) = D_1 (\omega^+)$ for otherwise 
$\filled (\varphi) \subset \varphi^{-1}(D_0) \subset B_{R_\varphi} (\omega^+)$ and $\diam \filled (\varphi) < R_\varphi$.
Now since $\varphi (B_{R_\varphi} (\omega^+)) = B_\hr (0)$ where $\hr$ is the radius of $\vphi(\db_0)$,
it follows that $\varphi(A_1(\omega^+)) =  A_0$
and the degree of $\varphi : A_1(\omega^+)  \rightarrow A_0$ is $2$.
The rest of (iv) follows along the same lines.
\hfill $\Box$

\begin{lemma}
  \label{modulus-action-l}
Let $\varphi$ be a cubic polynomial with critical points $\omega^\pm$ such that $R_\varphi = \operatorname{diam} \filled(\varphi)$.
Suppose that $\omega^- \notin \filled(\varphi)$. 
Consider $n \geq 1$ and assume that $\varphi^n (\omega^+) \in D_0$. 
Let $\znot$ be  a level $n$ point and let $E$ be an end.
Then the following hold:

(i) For any element $P$  of the affine partition associated to $D_{n-1}(\znot)$ or to $\db_{n-1} (E)$ there exists at most one ball of 
level $n$ contained in $P$ (see Subsection~\ref{affine-ss}).

(ii) $\varphi(A_n (\znot))  =  A_{n-1} (\varphi(\znot))$ and  $A_n (\znot) \subset \L \setminus \filled(\varphi)$.

(ii') $\varphi(A_n (E))  =  A_{n-1} (\varphi(E))$ and  $A_n (E) \subset \L \setminus \filled(\varphi)$.

(iii) $$ 
  \deg_{A_n(\znot)} (\varphi)  =  \left\{ \begin{array}{ll}
                                             1 & \mbox{ if } \omep \notin \db_n (\znot), \\
                                             2 & \mbox{ if }  \omep  \in \db_n (\znot). 
                                             \end{array}
                                   \right.
$$

(iii') $$ 
  \deg_{A_n (E)} (\varphi)  =  \left\{ \begin{array}{ll}
                                             1 & \mbox{ if } \omep \notin \db_n (E), \\
                                             2 & \mbox{ if } \omep  \in \db_n (E). 
                                             \end{array}
                                   \right.
$$
\end{lemma}

\noindent
{\bf Proof.}
We proceed by induction. For $n=1$ the previous lemma implies (i)--(iii).
Consider $n \geq 2$ and suppose that (i)--(iii) hold for $1, \dots, n-1$. 
We show that (i)--(iii) hold for $n$:

Note that $P \setminus \db_n (\znot)= A_n(\znot)$.

To prove (i) we proceed by contradiction and suppose that $P$ contains $D_n (\znot)$ and another level $n$ ball
$D_n (\z_1)$. By the inductive hypothesis,  the unique ball inside $\varphi (P)$ is $D_{n-1}(\varphi(\znot))$.
Therefore, $\deg_P (\varphi) = 2$ and 
$P$ contains the critical point $\omega^+$ which has to be outside $\varphi^{-1} (D_{n-1}(\varphi(\znot)))$.
Hence, $\varphi(\omega^+) \in A_{n-1} (\varphi(\znot)) \subset \L \setminus \filled (\varphi)$ which contradicts the hypothesis
of the lemma.

From (i) we have that $\varphi^{-1} (D_{n-1}(\varphi(\znot))) \cap P = D_{n}(\znot)$. Hence $\varphi (A_n (\znot))
= \varphi (P \setminus \db_n (\znot)) = \varphi (P) \setminus  \db_{n-1} (\varphi(\znot)) = A_{n-1}(\varphi(\znot))$
and (ii) follows.

For (iii) since $\varphi^{-1} (\da_{n-1}(\varphi(\znot))) \cap P = \da_{n}(\znot)$,  the degree
of $\varphi$ in $A_n(\znot)$ coincides with that of $\varphi$ in $P$.
The degree of $\varphi$ in $P$ is $1$ or $2$ according to whether $\omep \notin P$ or
$\omep \in P$.  From (i),  $\omep \in P$ if and only if $\omep \in D_n (\znot)$. Thus (iii) holds.

Choosing $\znot \in \db_n (E)$ parts (ii') and (iii') follow as well.
\hfill $\Box$

\medskip
Following Branner and Hubbard~\cite{branner-92} we will introduce marked grids and tableaux in order to keep track of 
the  moduli of annuli.

\begin{notation}{\em Let $\ell, k \geq 0$ be  integers.
Given $\z \in \L$ such that $\z_k = \varphi^k (\z)$ is a level $\ell$ point we denote by $A_{\ell, k} (\z)$ the level $\ell$ annulus
$A_\ell (\z_k)$ around $\z_k$. Similarly, if $E$ is an end, we denote by $A_{\ell,k}(E)$ the level $\ell$ annulus of $\varphi^k (E)$.}
\end{notation}

\begin{definition}
\label{tableaux-d}
  Let $\z$ be a level $n \in \N$ point. The {\bf level $n$ tableaux of $\z$}, denoted $\bT_n (\z)$ or simply $\bT (\z)$, 
is the two dimensional array:
$$ \bT ( \z) := (A_{\ell, k} (\z) )$$
where $0 \leq \ell, k$ and $\ell + k \leq n$. 
The {\bf level $n$ marked grid}, denoted $\bM_n (\z)$ or sometimes simply $\bM(\z)$ is the two dimensional array $(M_{\ell,k} (\z))$ where
$0 \leq \ell, k$ and $\ell + k \leq n$ and
$$M_{\ell, k}(\z) = \left\{ \begin{array}{ll}
1 & \mbox{if } A_{\ell,k}(\z) = A_{\ell}(\omep) \\
0 & \mbox{otherwise}.
\end{array} \right. $$
If $\z \in \filled (\vphi)$, then the {\bf tableaux of $\z$}
 is the infinite array  of annuli $\bT(\z) = ( A_{\ell,k}(\z))$ and the {\bf marked grid} is
$\bM(\z) = (M_{\ell, k} (\z))$ where $\ell, k \geq 0$.
Similarly, given an end $E$ we define the corresponding tableaux $\bT(E) = ( A_{\ell,k}(E))$  and
marked grid $\bM(E) = (M_{\ell, k} (E))$ where $M_{\ell, k} (E)$ is $0$ or $1$ according to whether 
$A_{\ell,k}(E)(\omep) \neq A_{\ell,k}(E)$ or $A_{\ell,k}(E)(\omep) = A_{\ell,k}(E)$.
\end{definition}

Marked grids are useful to compute the moduli of the annuli of the corresponding tableaux. In fact, from Lemma~\ref{modulus-action-l}, 
if $\z$ is a level $n$ point, let $$S_\ell = \sum^{\ell-1}_{i=0} M_{\ell-i,i} (\z).$$
Then $$\mdl A_\ell (\z) = 2^{-S_\ell} \mdl A_0$$ for all $\ell \leq n$. 

\smallskip
Marked grids satisfy four simple rules:

\begin{proposition}
\label{rules-p}
  Suppose that $\omep$ is a level $n$ point. Given a level $n$ point $\z$ (resp. an end $E$)
let $M_{\ell, k} = M_{\ell, k} (\z)$ (resp. $M_{\ell,k}=M_{\ell, k} (E)$) for $\ell + k \leq n$.  Then the following hold:

(Ma) If $\ell +k \leq n$ and $M_{\ell,k}$ is marked, then $M_{j,k}$ are marked
for all $j \leq \ell$.

(Mb) If  $\ell +k \leq n$ and $M_{\ell,k}$ is marked, then $M_{\ell-i,k+i} = M_{\ell-i,i}(\omep)$
for $0 \leq i \leq \ell$.

(Mc) If $\ell  + m +1 \leq n$ and $M_{\ell-i,i}(\omep)$ is not marked for all $0 < i <k$,
$M_{\ell+1-k,k}(\omep)$ is marked, $M_{\ell,m}$ is marked, and 
$M_{\ell+1,m}$ is not marked, then $M_{\ell+1-k,m+k}$ is
not marked.

(Md) If $\ell+k+1 \leq n$ and $M_{1,\ell}(\omep)$ is not marked, $M_{\ell,k}$ is marked,  $M_{\ell+1,k}$ 
is not marked and $M_{\ell-i,k+i}$ is not marked for $0 <i <\ell$, then $M_{1,k+\ell}$ is marked.
\end{proposition}

\begin{definition}
\label{admissable-marked-d}
A two dimensional array $\bM = (M_{\ell,k})_{\ell,k \geq 0}$ such that $M_{\ell, k} \in \{0,1\}$ is called an {\bf admissible
marked grid} if (Ma)--(Md) of Proposition~\ref{rules-p} hold for all $n$. If moreover $M_{\ell,0}$ is marked for all $\ell$, then
we say that $\bM$ is an {\bf admissible critical marked grid}. 
Similarly, an array $\bM_n =  (M_{\ell,k})$ where $\ell + k \leq n$ for which (Ma)--(Md) hold is called an {\bf admissible marked
grid of level $n$}. 
\end{definition}

\begin{remark}
  {\em The rule (Md) implies the fourth rule in~\cite{harris-99} but not conversely. In fact, 
let  $\bM = (M_{\ell, k})$ be such that  $1= M_{5,1} = M_{5,2} = M_{5,3}$ 
and $1 = M_{\ell,0} = M_{0,k}$ for all $\ell, k \geq 0$
but all the other positions are unmarked (i.e., $0$).
Then $\bM$ satisfies the first three rules but not the fourth. Such a grid is not the critical marked grid of a
cubic polynomial.}
\end{remark}

\noindent
{\bf Proof.}  We may assume that $M_{\ell, k} = M_{\ell, k} (\z)$ for some $\z \in \varphi^{-n}(D_0)$.
As usual, let $\z_k = \varphi^k (\z)$. 

(Ma) follows directly from the definitions.  

For (Mb) note that if $M_{\ell,k}$ is marked, then 
$A_{\ell, k}(\z) = A_{\ell, 0}(\omep)$. Therefore, $A_{\ell-i, k+i} = \varphi^i (A_{\ell}(\z_k)) = \varphi^i (A_{\ell} (\omep)) = 
A_{\ell-i, i} (\omep)$. 

Under the hypothesis of  (Mc) it follows that $\z_m \in \db_\ell (\omep) \setminus \db_{\ell+1} (\omep)$.
Since $\db_{\ell+1} (\omep)$ is the only preimage of $\db_\ell (\varphi(\omep))$ inside $\db_\ell(\omep)$,
we conclude that $\z_{m+1} \in \db_{\ell -1} (\varphi (\omep)) \setminus \db_\ell (\varphi(\omep))$.
Now $\varphi^{k-1}$ is one--to--one on $\db_{\ell -1} (\varphi (\omep))$, therefore 
$\z_{m+k} \in \db_{\ell-k} (\varphi^k(\omep)) \setminus \db_{\ell -k +1}(\varphi^k (\omep))$.
By assumption $M_{\ell -k +1, k} (\omep)$ is marked, thus  $\db_{\ell -k +1}(\varphi^k (\omep)) = \db_{\ell -k +1}(\omep)$.
Hence,  $\db_{\ell -k +1}(\z_{m+k}) \neq \db_{\ell -k +1}(\omep)$ and $M_{\ell-k+1,m+k} $ is unmarked.

Now under the hypothesis of (Md) we have that $\z_k \in \db_\ell(\omep) \setminus \db_{\ell+1} (\omep)$.
It follows that $\z_{k +\ell} \in \db_0 \setminus \db_1 (\varphi^\ell(\omep))$. By hypothesis, $\db_1 (\varphi^\ell(\omep))
= \db_1 (\gamma^+) \neq \db_1(\omep)$ where $\gamma^+$ is the cocritical point of $\omep$. Therefore,
$\z_{k+\ell} \in \db_1 (\omep)$ because there are only two level $1$ balls. Hence $M_{1, k +\ell}$ is marked.
\hfill $\Box$

The marked grid of the critical point plays a central role. If $\omep \in \filled (\varphi)$, then the 
critical marked grid $(M_{\ell, k}(\omep))$ is defined for all $\ell, k \geq 0$. In this case, the critical marked
grid is said to be {\bf periodic} of period $p>0$ if the $p$-th column is marked. That is, 
$M_{\ell, p} (\omep) =1$ for all $\ell \geq 0$ and $p>0$ is minimal with this property. 

For tableaux such that the corresponding grid satisfy (Ma)--(Mc) of Proposition~\ref{rules-p} and part (iii) of 
Lemma~\ref{modulus-action-l}, Branner and Hubbard (see Theorem 4.3 in~\cite{branner-92})  established the following:

\begin{theorem}[Branner and Hubbard]
\label{bh-t}
  Suppose that $\varphi$ is a cubic polynomial such that $R_\varphi = \diam \filled (\varphi)$, 
$\omep \in \filled (\varphi)$ and $\omega^- \notin \filled (\varphi)$. Then:

(i)  If the critical marked grid is not periodic, then 
$$\sum_{\ell \geq 0} \mdl A_{\ell} (E)$$
is divergent for all ends.

(ii) If the critical marked grid is periodic, then
$$\sum_{\ell \geq 0} \mdl A_{\ell} (E)$$
is convergent if and only if there exists $k \geq 0$ such that $A_{\ell, k}(E) = A_\ell (\omep)$ for all $\ell \geq 0$. 
\end{theorem}

From Lemma~\ref{sum-l} we obtain the following:
\begin{corollary}
  \label{bh-c}
 Suppose that $\varphi$ is a cubic polynomial  such that $R_\varphi = \diam \filled (\varphi)$, 
$\omep \in \filled (\varphi)$ and $\omega^- \notin \filled (\varphi)$. Then:

(i) The critical marked grid is not periodic if and only if the intersection of every end is a singleton.

(ii) The critical marked grid is periodic if and only if the infraconnected component of $\filled (\vphi)$ which contains
the critical point is a periodic ball.
\end{corollary}

\subsection{Dynamical space results.}
Now we are ready to prove a stronger version of Theorem~\ref{wandering-it}

\begin{theorem}
  \label{wandering-t}
Let $\varphi \in \L[\z]$ be a cubic polynomial. Then the following hold:

(i) Every end of $\varphi$ has non--empty intersection.

(ii) Every infraconnected component of $\filled (\varphi)$ is either a closed ball or a point.

(iii) An infraconnected component $C$ of $\filled(\varphi)$ is a closed ball 
if and only if $C$ eventually maps onto a periodic infraconnected component containing a critical point.
\end{theorem}

\noindent
{\bf Proof.}
By Theorem~\ref{cantor-th} and Corollary~\ref{connected-c}, we may assume that $\varphi$ has one critical point $\omega^-$
escaping to $\infty$ and another one $\omep$ in $\filled (\varphi)$. 
Moreover, we may also assume that $\vphi$ is normalized so that $R_\vphi = \diam \filled (\vphi)$ (Lemma~\ref{diam-l}). 
Therefore the definitions and results of Subsection~\ref{bh-ss} apply to $\vphi$.

Let $E$ be an end.  If $r_n(E) \rightarrow 0$ then
the intersection of $E$ is a point. If $r_n (E) \not\rightarrow 0$, then for some $k \geq 0$ we have
that  $\omep \in \varphi^k (\db_n (E))$ for all $n \geq 0$, by Theorem~\ref{bh-t}. In particular $\cap \db_n (E) \neq \emptyset$ 
and (i) follows.

Also note that $\filled (\varphi)$ has a non--trivial infraconnected component if and only if the critical marked grid of $\varphi$ is
periodic. In this case, every non--trivial infraconnected component eventually maps onto the periodic 
infraconnected component $C (\omep) = \cap \db_n (\omep)$.
Therefore, to finish the proof of the theorem, it suffices to show that $C(\omep)$ is a closed ball when the critical marked grid
is periodic. In fact, if $\bM(\omep)$ is periodic, say of period $p$, then there exists $\ell_0$ such that
$$\mdl A_{\ell+p} (\omep) = \frac{1}{2} \mdl A_\ell (\omep)$$
for all $\ell \geq \ell_0$. 
It follows that 
$$a= \frac{1}{\mdl A_0} \sum^\infty_{\ell =0} \mdl A_\ell (\omep)$$ 
is rational and therefore the radius of $C(\omep)$ is $\hr^{a-1} R^a_\varphi \in \valgL$ where 
$\hr$ is the radius of $\vphi(\db_0)$.
\hfill
$\Box$

\begin{corollary}
  \label{almost2-c}
  Let $\vphi \in \cubicl$ be a cubic polynomial. Then the following are equivalent:

(i) $\julia (\vphi)$ is compact.

(ii) $\julia (\vphi) = \filled (\vphi)$.

(iii) All the cycles of $\vphi$ are repelling.

(iv) $\vphi$ either has all its critical points escaping or it has exactly one critical point (counting multiplicities) 
which is non-escaping and the corresponding marked grid is not periodic. 
\end{corollary}

\medskip
\noindent
{\bf Proof.}
In view of Lemma~\ref{infra-l} and Corollary~\ref{compact-c} and we only have to prove that (iii) implies (iv)
and (iv) implies (i).
 
(iii) $\implies$ (iv):
If (iv) does not hold, then
 $\vphi$ is in the infraconnectedness locus or is in $\escapepml$ but the corresponding critical marked grid is periodic.
In both cases  there exists a periodic infraconnected component $B$ say 
of period $p$ which is a closed ball that contains a critical point (Corollary~\ref{bh-c}).
By Lemma~\ref{fixed-l}  , $\varphi^p : B \rightarrow B$ would have a fixed point, which by Schwarz Lemma~\ref{schwarz-l}
would be non--repelling. Hence (iii) does not hold.

(iv) $\implies$ (i): 
If we assume that $\vphi \in \shiftcubicl$ or $\vphi \in \escapepml$ and the corresponding critical marked grid is not periodic,
then the intersection of every end is a point.
Now let $\{ \z_k \} \subset \julia (\varphi)$ be a sequence. Then there exists an end $E$ and a subsequence $\{ \z_{k_i} \}$
such that for all $n$, there exists $i_0$ for which $\z_{k_i} \in \db_n(E)$ for all $i \geq i_0$. Since the intersection
of $E$ is a point, say $\z$, it follows that $\z_{k_i} \rightarrow \z \in \julia (\varphi)$.
\hfill $\Box$

\medskip
We end this section with a basic result  about the topological entropy of cubic polynomials:

\begin{proposition}
\label{top-ent-p}
  Let $\vphi : \L \rightarrow \L$ be a cubic polynomial. For a compact invariant subset $X$ of $\julia(\vphi)$ we denote
by $h_{top} (\vphi, X)$ the topological entropy of $\vphi: X \rightarrow X$. If $\julia  (\vphi) \neq \emptyset$ then
$$h_{top} (\vphi):= \sup_X h_{top}(\vphi, X) = \log 3.$$
\end{proposition}

\noindent
{\bf Proof.}
Suppose that $\vphi$ is normalized so that $R_\vphi = \diam \filled (\vphi)$.
We may assume that $\vphi$ has exactly one critical point $\omep$  in $\filled (\vphi)$, for otherwise $\julia (\vphi) = \emptyset$ or
$\vphi : \julia (\vphi) \rightarrow \julia (\vphi)$ is topologically conjugated to the one-sided shift on $3$ symbols. In the
latter case the topological entropy is clearly $\log 3$.

Let $\mathbf{Ends}$ be the set of all ends of $\vphi$ endowed with the metric defined by 
$\rho ( \{ \db_n \} , \{ \db^\prime_n \}) =1/(k+1)$ if $k$ is the largest integer such that $\db_k = \db^\prime_k$.
Denote by $\vphi_\#$ the action induced by $\vphi$ on $\mathbf{Ends}$.
For $\z \in \filled (\vphi)$,
the map $\pi: \z \mapsto \{ \db_n (\z) \}$  is a semiconjugacy between $\vphi: \filled (\vphi) \rightarrow \filled (\vphi)$ 
and $\vphi_\#:\mathbf{Ends} \rightarrow \mathbf{Ends}$. 
Since the number of dynamicals balls of level $n$ is $(3^n +1)/2$, it follows that the 
topological entropy of $\vphi_\#$ is exactly $\log 3$. 

If the marked grid of $\omep$ is not periodic, then $\julia (\vphi)$ is compact and  $\pi :\julia (\vphi) \rightarrow \mathbf{Ends}$ 
is a topological conjugacy. Hence, the claim of the proposition follows in this case.

In the case that the marked grid of $\omep$ is periodic denote by $E^*$ the critical end.
From Lemma~\ref{label-l} we obtain  consider a labelling $L$ of the
level $n$ balls  which, after switching symbols if necessary, is such that:
 $$\lim_{n \rightarrow \infty} L(\db_n(E^*) = \{ 1, 2 \}.$$
  Now consider the itinerary function
$$\itin : \mathbf{Ends} \rightarrow \{\{1\},\{2\},\{3\}, \{1,2\} \}^{\N \cup \{0\}}$$
defined by 
$$\itin (E) = (i_k(E))_{k \geq 0}=  (\lim_{n\rightarrow \infty} L(\db_n (\vphi^k_\#(E))))_{k \geq 0}.$$
Observe that $i_k(E) = \{1,2\}$ if and only if $\vphi^k_\# (E)$ is the critical end $E^*$.
Moreover the image of $\mathbf{Ends}$ is characterized as the sequences $(i_k)$ such that
if there exists $k \geq 0$ so that $i_{k + \ell} \subset i_\ell (E^*)$ for all $\ell \geq 0$, then 
$i_{k + \ell} = i_\ell (E^*)$.  Also, $\itin \circ \pi$ is
injective over the Julia set and $\itin \circ \pi (\julia (\vphi))$ is the set of all itineraries (in the image of $\itin$) with no
symbol equal to $\{ 1, 2\}$. So it is sufficient to construct compact subsets of $\itin \circ \pi (\julia (\vphi))$ invariant under
the one-sided shift $\sigma$ with topological entropy arbitrarily close to $\log 3$. 
For this, let $p$ denote the period of $E^*$ and for each $N >p$, consider
the set $Y_N$ of all symbol sequences $(i_k)$ with $i_k \neq \{1,2\}$ for all $k \geq 0$ and
such that for some $0 \leq j <N$ and all $\ell \geq 0$
$$i_{\ell N +j} = i_{\ell N +1 +j} = \cdots = i_{\ell N +p-1+j} = \{3\}.$$  
Since the topological entropy of $\sigma : Y_N \rightarrow Y_N$ is $(1- pN^{-1}) \log 3$, the proposition follows.
\hfill $\Box$ 
 
\medskip
It is worth to mention that
$\vphi_\# : \mathbf{Ends} \rightarrow \mathbf{Ends}$ is topologically conjugated to the dynamics of $\vphi$ over 
the Julia set of $\vphi$ in the Berkovich analytic space induced by  $\L$ (compare with~\cite{rivera-03,rivera-03a}). The complement
of $\pi(\julia(\vphi))$ consists of all ends which have empty or non-trivial intersection.
Favre and Rivera's
construction in~\cite{favre-04} produces an equilibrium measure supported in the Berkovich space Julia set of $\vphi$.
It is natural to expect that their measure corresponds to a measure of maximal entropy for $\vphi_\#$.

\section{Parameter space}
\label{parameter-s}
Recall that  we work in the parameter space $\cubicl$ of 
monic centered cubic polynomials with marked critical points. That is,
$$\cubicl:= \{ \varphiab \in \L[\z] \,\, / \, \varphiab(\z) = \z^3 - 3 \alpha^2 \z + \beta \mbox{ for some } 
(\ab) \in \L^2 \}$$
which is naturally identified with $\L^2$. Note that the critical points of $\vphiab$ are $\pm \alpha$.

The {\bf infraconnectedness locus $\coneccubicl$} is the subset of $\cubicl \equiv \L^2$
formed by all the cubic polynomials with infraconnected filled Julia set.
According to Corollary~\ref{connected-c}, $\varphiab \in \coneccubicl$ if and only if
$\pm \alpha \in \filledab$. 

The {\bf shift locus $\shiftcubicl$} is the subset of $\cubicl \equiv \L^2$ formed by all
the cubic polynomials with both critical points escaping.

The rest of parameter space splits into  two sets  $\cE^+(\L)$ and $\cE^-(\L)$
where $$\cE^+(\L) := \{ \varphiab \,\, /\,  +\alpha \in \filled (\varphi) \not\ni -\alpha \}$$
and
$$\cE^- (\L) := \{ \varphiab \,\, /\,  -\alpha \in \filled (\varphi) \not\ni +\alpha \}.$$

\smallskip
A quick computation leads to the following result.

\begin{proposition}
  \label{glance-p}
  \begin{eqnarray*}
    \coneccubicl & = & B^+_1 (0) \times B^+_1 (0).\\
    \shiftcubicl    & \supset & \{ \varphiab  \notin \coneccubicl \,\, /\,   \lval \alpha \rval^3 \neq \lval \beta \rval \}.\\
    \cE^\pm (\L)        & \subset &  \{ \varphiab \,\, /\, 1<  \lval \alpha \rval^3 = \lval \beta \rval \mbox{ \em{and} } \lval \mp 2  \alpha^3 + \beta  \rval \leq \lval  \alpha \rval \}.\\
  \end{eqnarray*}
\end{proposition}

Since conjugation of $\varphiab$ by $\z \mapsto -\z$ 
gives $\varphi_{-\alpha,-\beta} $, to describe how polynomials are organized in 
  $$\{ \varphiab \,\, /\, 1<  \lval \alpha \rval^3 = \lval \beta \rval \mbox{ and } \lval \mp 2  \alpha^3 + \beta  \rval \leq \lval  \alpha \rval \}$$
it is sufficient to understand the structure of
$$ \{ \varphiab \,\, /\, 1<  \lval \alpha \rval^3 = \lval \beta \rval \mbox{ and } \lval - 2 \alpha^3 + \beta  \rval \leq \lval \alpha \rval \}.$$
Our aim is  to study this set in detail. For this purpose:

\bigskip
{\bfseries\slshape 
Throughout  this section  we fix $\alpha \in \L$ such that $\lval \alpha \rval >1$ and let $\vphib = \vphiab$.}

\bigskip
Let $$\cEvphi_0 := 
\{ \beta \,\, / \,  \lval \alpha \rval^3 = \lval \beta \rval \mbox{ and } \lval - 2 \alpha^3 + \beta  \rval \leq \lval \alpha \rval \}.$$
To simplify notation, we identify $\cEvphi_0$ with $\{ \vphib \in \cubicl \,\, / \, \beta  \in \cEvphi_0 \}$.

\begin{remark}
  {\em We use the upper-script $\vphi$ to distinguish the sets associated to the family $\vphib$ 
from the corresponding sets associated to  another  family    $\psi_\nu$. The family $\psinu$ will be  introduced
in the next subsection.}
\end{remark}

\begin{lemma}
  Let $\varphi = \vphib \in \cEvphi_0$. Then the following hold:

(i) $R_\vphi = \lval \alpha \rval = \diam \filled (\vphi)$.

(ii) $\vphi (- \alpha ) \notin \dbb_0 \ni \vphi(\alpha)$ where $\dbb_0 = B^+_{R_\vphi} (0)$ is the level $0$ ball of $\vphi$.
\end{lemma}

\noindent
{\bf Proof.} From the definition of $R_\vphi$ it follows that $R_\vphi = \lval \alpha \rval$.
By inspection of the Newton polygon of $\vphi(\z) - \z$ it follows that
$\vphi$ has $3$ fixed points $\z_1, \z_2, \z_3$ in $\{ \lval \z \rval =  \lval \alpha \rval \}$.
Since $\z_1 + \z_2 + \z_3 =0$, at least two of the fixed points  are at distance $\lval \alpha \rval$ from each other.
Hence $ \lval \alpha \rval = \diam \filled (\vphi)$. The rest of the lemma is also straightforward.
\hfill $\Box$

\medskip
From the previous lemma, for all  $\vphib \in \cEvphi_0$ we have that 
$R_\vphib = \diam \filled (\vphib)$, $-\alpha \notin \filled (\vphib)$ and $\vphib (\alpha) \in \dbb_0$.
Thus, the assumptions and therefore the definitions 
and results contained in Subsection~\ref{bh-ss} apply to $\vphib \in \cEvphi_0$.

\begin{theorem}
\label{parameter-th}
  Consider an admissible  critical marked grid $\bM$ and
let $$ C_{\bM} := \{ \beta \in \cEvphi_0  \,\,/\, \alpha \in \filled(\vphi_\b) \mbox{ and } \bM = \bM^\b (\alpha) \}$$
where $\bM^\b (\alpha)$ is the marked grid of the critical point $\alpha$ under iterations of $\vphib$.
Then the following hold:

(i) If $\bM$ is periodic, then $C_{\bM}$ is a non-empty union of finitely many closed and pairwise disjoint balls.

(ii) If $\bM$ is not periodic, then $C_{\bM}$ is a non-empty compact set and  
$$C_{\bM} \subset  \partial \{ \beta \in \cEvphi_0 \,\,/\, \vphib \in \shiftcubicl \}. $$
\end{theorem}

We prove this theorem in Subsection~\ref{realization-ss}.
The proof relies on describing  how polynomials are organized in $\cEvphi_0$.
To describe $\cEvphi_0$, for $n \geq 0$ we introduce the sets
$$\cEvphi_n := \{ \vphib \in \cEvphi_0 \,\,/\, \vphi^{n+1}_\b (\alpha) \in \dbb_0 \}.$$
Note that  $\cEvphi_n$ is a finite disjoint union of closed balls.
Each of these balls  is called a {\bf $\vphi$--parameter ball of level $n$}.

\medskip
A level $n$ dynamical ball of $\vphib \in \cEvphi_0$  which contains $\z$ will be denoted by $\dbb_n (\z)$.
The level $n$ marked grid  of a point $\z$ by $\bM_n^\b (\z)$ and the corresponding  entries by $M^\b_{\ell,k}(\z)$.

\begin{definition}
  \label{center-d}
Let $n \in  \N$.
We say that $\vphib \in \cEvphi_0$ is a {\bf center of level $n$} if for some $p \geq 1$:

(i) $ \vphi^p_\b (\alpha)= \alpha$ and,

(ii) $\alpha \notin \vphi^k_\b (\dbb_{n+1} (\alpha))$ for $k=1, \dots, p-1.$

We say that $p$ is the {\bf period} of the center $\vphib$.
\end{definition}

The correspondence between level $n$ dynamical and parameter balls is stated in the next proposition.

\begin{proposition}
  \label{level-n-p}
  Let $\pb_n$ be a level $n$ parameter ball.
Then the following hold:

(i) $\pb_n = \dbb_n (\vphib(\alpha)) + 2 \alpha^3$ for all $\beta \in \pb_n$.

(ii) $\bM^\b_{n+1} (\alpha) = \bM^{\beta^\prime}_{n+1} (\alpha)$ for all $\beta, \beta^\prime \in \pb_n$.

(iii) There exists a unique center of level $n$ 
in $\pb_n$. The period of this center is $\min \{ k \geq 1  \,\, /\, M^\b_{n+1-k, k} (\alpha) =1 \}$ for any $\beta \in \pb_n$.
\end{proposition}

The proof of this proposition is in Subsection~\ref{paraballs-ss}.

In particular, the above proposition shows that the radius of $\pb_n$ is easily computed from 
$\bM^{\beta}_{n+1} (\alpha)$ for any  $\beta \in \pb_n$ and coincides with the radius of the level
$n$ dynamical ball around the critical value $\vphib (\alpha)$.  The proposition also says that  
if $\beta \in \pb_n$ and the critical point $\alpha$ is  periodic of period $q$ under $\vphib$, then
$q \geq p$ where $p$ is the period of the center of $\pb_n$.
Moreover, $p=q$  if and only if $\vphib$ is the unique level $n$ center in $\pb_n$. 

\smallskip
The next proposition  describes the correspondence between level $n+1$ parameter and dynamical 
balls.

\begin{proposition}
\label{level-n+1-p}
  Consider a level $n$ parameter ball  $\pb_n$ and let $P$ be an element of the affine partition associated to $\pb_n$. 
For any $\beta \in \pb_n$ we have the following:

There exists a level $n+1$ parameter ball contained in $P$ if and only if there exists a level $n+1$ dynamical ball 
$\dbb_{n+1} (\z)$ contained in  $P-2\alpha^3$.
In this case, the level $n+1$ parameter ball $\pb_{n+1}$ is unique and 
$$\bM^\b_{n+1} (\z) = \bM^{\beta^\prime}_{n+1} (\vphi_{\beta^\prime}(\alpha))$$
for all $\beta^\prime \in \pb_{n+1}$. In particular, the radii of  $\dbb_{n+1} (\z)$ and $\pb_{n+1}$ coincide.
\end{proposition}

The proof of this proposition is also given is Subsection~\ref{paraballs-ss}. 
The above propositions are easier to prove after a change of coordinates in the dynamical and parameter spaces.

\subsection{Change of coordinates}
\label{change-ss}
To prove our parameter space results is more comfortable to work with the family
$$\psinu (\z) = \psianu (\z) = \alpha^2 (\z -1)^2 (\z+2) + \nu$$
where $\nu  \in \L$. 
Note that the critical points of $\psinu$ are $\omega^\pm = \pm1$ and $\psinu (-2) = \psinu (\omep = +1) = \nu$.
Moreover, since $\alpha \neq 0$ the polynomial $\psinu$ is conjugate via $\z \mapsto \alpha \z$ to $\vphib$ where
$$\beta = \alpha \nu + 2 \alpha^3.$$
Observe that $\vphib \in \cEvphi_0$ if and only if  $\lval \nu \rval \leq 1$
so we let 
$$\cEpsi_0 := \{ \nu  \in \L  \,\, /\, \lval \nu \rval \leq 1 \} \equiv \{ \psinu \,\,/\, \lval \nu \rval \leq 1\}.$$

For all $\nu \in \cEpsi_0$, the level $0$ dynamical ball $\dbnu_0$ of $\psinu$
is $B^+_1(0)$ and $\psinu (\dbnu_0) = B^+_{\lval \alpha \rval^2} (0)$. 
Also, $\psinu (\omep) \in \dbnu_0$, $\psinu (\omega^-) \notin \dbnu_0$ and
$\diam \filled (\psinu)=1$. So we are under the assumptions of Subsection~\ref{bh-ss}.

 We let
$$\cEpsi_n := \{ \nu  \in \cEpsi_0 \,\,/\, \psinu^{n+1}(\omep) \in \dbnu_0 \}.$$
This set  is also the union of finitely many closed and disjoint balls which we call {\bf $\psi$-parameter balls of level $n$}.

\smallskip
Here we denote the level $n$ ball of $\psinu \in \cEpsi_0$ containing $\z$ by $\dbnu_n (\z)$, 
the radius of $\dbnu_n (\z)$
by $r^\nu_n (\z)$, the level $n$ annulus around $\z$ by $\danu_n (\z)$ and 
the level $n$ marked grid of $\z$ by $\bM^\nu_n (\z)$ with entries $M^\nu_{\ell,k}(\z)$.

\medskip
Similarly than in the $\vphi$-parameter space we say that $\psinu \in \cEpsi_0$ is a {\bf center of level $n$} if

(i) $\psi^p_\nu (\omep) = \omep$ and

(ii) $\omep \notin \psi^k_\nu (\dbnu_{n+1} (\omep))$ for $k=1, \dots, p-1.$

\smallskip
The next lemma is a straightforward consequence of the change of coordinates involved. We omit the proof.

\begin{lemma}
\label{coordinates-l}
  Let $\beta$ and $\nu$ be such that 
$\beta = \alpha \nu + 2 \alpha^3 \in \cEvphi_0$. Then the following hold:

(i) $\db^\nu_n$ is a level $n$
dynamical ball of $\psinu$ if and only if $\db^\b_n = \alpha \dbnu_n$ is a level $n$ dynamical ball of $\vphib$.

(ii) $\psinu \in \cEpsi_n$ if and only if $\vphib \in \cEvphi_n$.
In particular, 
$\pb^\psi_n$ is a $\psi$--parameter ball of level $n$ if and only if 
$\pb^\vphi_n= \alpha \pb^\psi_n + 2 \alpha^3$ is a $\vphi$--parameter ball of level $n$.

(iii) $\psinu$ is a center of level $n$ if and only if $\vphib$ is a center of level $n$.
\end{lemma}

\subsection{Thurston map}
Our next result is the key to prove propositions~\ref{level-n-p} and~\ref{level-n+1-p}. It shows that
given a polynomial $\psinu$ and a level $n+1$ dynamical ball $\db_{n+1}$ inside the critical value ball of level $n$ there exists
a parameter $\nup$ close to the level $n+1$ ball $\db_{n+1}$ such that the critical point of $\psinup$ is periodic with orbit close
to that of the points in $\db_{n+1}$. The precise statement is as follows:

\begin{proposition}
\label{fixed-p}
  Consider a parameter $\hnu \in \cEpsi_n$ and let $\hzone$ be a level $n+1$ point
such that $\dbhnu_{n+1}(\hzone) \subset \dbhnu_n(\hnu)$.
For $k \geq 0$, let $\hz_{k+1} = \psi^k_\hnu (\hzone)$ and 
$$\hp = \min \{ k \geq 1 \,\, /\, \omep \in \dbhnu_{n+2-k}(\hz_k) \}.$$
Then there exists a unique $\nup$ such that:

(i) $\lval \nup - \hzone \rval < r^\hnu_n (\hnu)$.

(ii) $\psinup^{\hp} (\omep) = \omep$.
\end{proposition}

This subsection is devoted to the proof of this  proposition so throughout  we consider $\hnu$,  $\hz_k$ and $\hp$ as above.
The parameter $\nup$ is obtained as the first coordinate of the  fixed point of an appropriate ``Thurston map'' which acts on:
$$\B := \{ (\z_1, \dots, \z_{\hp} = \omep) \,\, /\, \lval \z_k - \hz_k \rval < \rho_k \}$$
where $\rho_k = r^\hnu_{n+1-k}(\hz_k)$.

\smallskip
We start with two lemmas which apply to an arbitrary $\nu \in \cEpsi_n$:

\begin{lemma}
  \label{apriori-l}
  Let $n \geq 0$ be an integer and consider $\nu \in \cEpsi_n$. 
  For $k \geq 1$, let $\nu_k = \psi^k_\nu (\omep)$. Then

(i) 
 \begin{equation*}
    r^\nu_n(\nu)  <  r^\nu_{n+1-k}(\nu_k) \mbox{ for } k=2,\dots,n.
  \end{equation*}

(ii) Assume that $\zone$ is a level $n+1$ point in $\db^\nu_n(\nu)$. For $k \geq 1$, let $\z_{k+1} = \psi^k_\nu (\zone)$ and
$$p =  \min \{ k \geq 1 \,\, / \, \omep \in \db^\nu_{n+2-k} (\z_k) \}.$$
Then
\begin{equation*}
  r^\nu_{n+2 - (k+1)} (\z_{k+1})  =  \lval \alpha \rval^2 \cdot \lval \z_k - \omep \rval \cdot r^\nu_{n+2-k}(\z_k)
\end{equation*}
for all $k=1, \dots, p-1$.
\end{lemma}

\noindent
{\bf Proof.}
For each $\ell$ such that $0 \leq \ell \leq n+1-k$ let $\delta(\ell) \geq 1$ be the  integer such that
$\da^\nu_{\ell+\delta(\ell)} (\nu_{k-\delta(\ell)})$ is critical but $\danu_{\ell+i}(\nu_{k-i})$ is not critical for all
$0 < i < \delta(\ell)$. 
To find such an integer $\delta(\ell)$ start at $A_{\ell,k}(\omep)$ in the critical marked grid and follow the southwest 
diagonal until you hit a critical position. The number of columns that you moved to the left is $\delta(\ell)$. 
Note that $\ell + \delta(\ell)$  is a strictly increasing function of $\ell$ and $0 < \ell + \delta(\ell) \leq n+1$. 
Moreover, 
$$2 \mdl \danu_{\ell+\delta(\ell)} (\omep) = 2 \mdl \danu_{\ell+\delta(\ell)} (\nu_{k-\delta(\ell)}) = \mdl \danu_\ell(\nu_k).$$
Since $\psinu (\dbnu_0) = B^+_{\lval \alpha \rval^2}(0)$ and $\mdl A^\nu_0 = \log \lval \alpha \rval^2$,
\begin{eqnarray*}
  \log \lval \alpha \rval^2  - \log r^\nu_n(\nu) & = & \sum^n_{\ell=0} \mdl \danu_\ell (\nu) 
                                                         =  2 \cdot  \sum^{n+1}_{\ell=1} \mdl \danu_\ell (\omep) \\
                                                        & > & 2 \cdot \sum^{n+1-k}_{\ell =0}  \mdl \danu_{\ell+\delta(\ell)} (\omep)    
                                                         =   \sum^{n+1-k}_{\ell =0} \mdl \danu_\ell(\nu_k) \\
                                                        & = &   \log \lval \alpha \rval^2  - \log r^\nu_{n+1-k} (\nu_k).
\end{eqnarray*}
Hence (i) follows.

Now  we prove (ii). Fix $k \geq 1$ and let $\ell$ be such that $\danu_{\ell}(\z_k)$ is marked but $\danu_{\ell+1}(\z_k)$ is not.
It follows that $\dbnu_\ell(\omep) = \dbnu_{\ell}(\z_k)$ and $\dbnu_{\ell+1}(\omep) \neq \dbnu_{\ell+1}(\z_k)$.
By Lemma~\ref{modulus-action-l} (i),  $$r^\nu_\ell(\z_k) = \lval \z_k - \omep \rval.$$
Note that:
$$2 \mdl \danu_j (\z_k) = \left\{ \begin{array}{ll}
                                                   \log  \lval \alpha \rval^4   & \mbox{if } j=0 \\
                                                   \mdl \danu_{j-1}(\z_{k+1})  & \mbox{if }  1 \leq j \leq \ell 
                                                    \end{array} \right.$$
Also, $\mdl  \danu_j (\z_k) =  \mdl \danu_{j-1}(\z_{k+1}) $ if $j \geq \ell +1$.
Hence:
\begin{align*}
\begin{split}
  \log  \lval \alpha \rval^6 - \log r^\nu_{n+2-(k+1)} (\z_{k+1}) & =  
                                              \log \lval \alpha \rval^4 + \sum^{n+2-(k+1)}_{j=0} \mdl \danu_j (\z_{k+1}) \\
& = 2 \cdot \sum^\ell_{j=0} \mdl \danu_j (\z_k) + \sum^{n+2-k}_{j = \ell +1} \mdl \danu_j (\z_k) \\
&= \log  \lval \alpha \rval^2 - \log r^\nu_{\ell} (\z_k) +  \log  \lval \alpha \rval^2  \\
& \quad \quad \quad \quad \quad \quad \quad  \quad \quad \quad \quad - \log r^\nu_{n+2-k} (\z_k).
\end{split}
\end{align*}
Statement (ii) follows after replacing $r^\nu_\ell(\z_k) $ by  $\lval \z_k - \omep \rval.$
\hfill
$\Box$

\begin{lemma}
\label{apriori2-l}
Consider $n \geq 1$ and  $\nu \in \cEpsi_n$. Let $\nu^\prime \in \dbnu_n (\nu)$ and, 
for all $k \geq 0$, let $\nu_k = \psi^k_\nu (\omep)$ and $\nu^\prime_k=\psi^k_{\nup} (\omep)$.
Then  for all $k$ such that $0 \leq k \leq  n$ the following hold:

(i) $\dbnu_{n+1-k}(\nu_k) = \dbnup_{n+1-k} (\nup_k)$. In particular, $\bM^\nu_{n+1} (\omep) = \bM^{\nup}_{n+1}(\omep)$.

(ii) Let $\cP_{n+1-k}$ be the affine partition associated to $\dbnu_{n+1-k}(\nu_k)$, then:
$$\psi_{\nu *} = \psi_{\nup*} : \cP_{n+1-k} \rightarrow \cP_{n+1-(k+1)}.$$
\end{lemma}
                                                   
\noindent
{\bf Proof.}
Let $k$ be such that $0 \leq k \leq n$ and $\zp \in \dbnu_{n+1-k} (\nu_k)$.
Then 
\begin{eqnarray*}
  \lval \psinup (\zp) - \nu_{k+1} \rval & = & \lval \psinup (\zp) - \psinu(\zp) + \psinu( \zp) - \nu_{k+1} \rval\\
                                                         & \leq & \max \{  r^\nu_n(\nu), \lval \psinu( \zp) - \nu_{k+1} \rval\}\\
                                                         & <     & r^\nu_{n+1-(k+1)} (\nu_{k+1}). 
\end{eqnarray*}
Therefore, $\psinup (\zp) \in \dbnu_{n+1-(k+1)} (\nu_{k+1})$.
Hence $\psinup^{n+1-k} (\dbnu_{n+1-k} (\nu_k)) \subset \db_0$ and 
$\nup_k  \in \dbnu_{n+1-k} (\nup_k) \subset \dbnup_{n+1-k} (\nu_k) $.
In particular, $r^{\nup}_n (\nup) \geq r^\nu_n(\nu)$ and $\nu \in \dbnup_n (\nup)$. After switching $\nu$ for $\nup$
and repeating the above argument it follows that $\dbnup_{n+1-k} (\nup_k) \subset \dbnup_{n+1-k} (\nu_k) $
and (i) follows.

For (ii), let $P$ be an element of the partition $\cP_{n+1-k}$ and choose $\z \in P$.
Then 
$$\lval \psinu (\z) - \psinup(\z) \rval \leq r^\nu_n (\nu) < r^\nu_{n+1-(k+1)}.$$
Therefore, $\psinup(P) =\psinu(P) \in \cP_{n+1-(k+1)}$.
\hfill
$\Box$

\begin{lemma}
\label{thurston-l}
Let   $$\B := \{ (\z_1, \dots, \z_{\hp} = \omep) \,\, /\, \lval \z_k - \hz_k \rval < \rho_k \}$$
where $\rho_k = r^\hnu_{n+1-k}(\hz_k)$. 
Then for each $(\z_1, \dots, \z_{\hp}) \in \B$ there exists a unique $(\zp_1, \dots, \zp_{\hp}) \in \B$ such that
$\psi_{\z_1}(\zp_k) = \z_{k+1}$ for $1 \leq k < \hp$.
\end{lemma}

\noindent
{\bf Proof.}
Let  $\cP_{n+1-k}$ be the partition associated to $\dbhnu_{n+1-k} (\hz_k)$
and denote by $P(\z) = B_{\rho_k}(\z)$ the element of $\cP_{n+1-k}$ that contains $\z$.
Since $\psi_{\zone*} = \psi_{\hnu*} : \cP_{n+1-k} \rightarrow  \cP_{n+1-(k+1)}$ and $\psihnu: P(\hz_k) \rightarrow 
P(\hz_{k+1})$ is one-to-one for $1 \leq k < \hp$ we have that $\psi_{\zone} : P(\hz_k) \rightarrow 
P(\hz_{k+1})$ is also one-to-one. Therefore, there exists a unique $\zp_k \in  P(\hz_k) $ such that
$\psi_{\zone} (\zp_k) = \z_{k+1} \in P(\hz_{k+1})$.
\hfill
$\Box$

\medskip
We may now define a Thurston map as
$$\begin{array}{rccc}
T: & \B & \rightarrow & \B \\
& (\z_1, \dots, \z_{\hp}) & \mapsto & (\zp_1, \dots, \zp_{\hp})
\end{array}$$
if $\psi_{\z_1}(\zp_k) = \z_{k+1}$  for all $k=1, \dots, \hp-1$.

\begin{lemma}
  A parameter $\nup$ is such that (i) and (ii) of Proposition~\ref{fixed-p} hold if and only if 
$(\nup, \psinup^2 (\omep), \dots, \psinup^{\hp-1} (\omep), \omep)$ is a fixed point of $T$.
\end{lemma}

\noindent
{\bf Proof.} Given $\nup$ such that (i) of Proposition~\ref{fixed-p} holds, by Lemma~\ref{apriori2-l} (ii),
for $k=1, \dots, \hp-1$, we have that $\lval \psinup^k (\omep) - \hz_k \rval <  \rho_k$. Therefore, 
$(\nup, \psinup^2 (\omep), \dots, \psinup^{\hp-1} (\omep), \omep)$ belongs to $\B$ and clearly is a fixed point
of $T$. The converse is straightforward.
\hfill $\Box$

\smallskip
It follows that to prove the proposition is sufficient to show that $T$ has a unique fixed point.

\smallskip
In $\L^{\hp}$ we consider the $\sup$-norm:
$$\| \vec{\z} = (\z_1, \dots, \z_\hp) \|_\infty := \max \{ \lval \z_1 \rval, \dots, \lval \z_\hp \rval \}.$$

\begin{lemma}
  For all $\vec{\z}= (\z_1, \dots, \z_{\hp}) \in \B$ we have that  $T^n (\vec{\z})$ converges to a fixed point of $T$.
\end{lemma}

\noindent
{\bf Proof.}
Consider $(\zone^{(0)}, \dots, \z^{(0)}_\hp) \in \B$ and let
$$(\zone^{(n)}, \dots, \z^{(n)}_\hp)= T^n (\zone^{(0)}, \dots, \z^{(0)}_\hp).$$
For $k=1, \dots, \hp-1$, 
\begin{eqnarray*}
  \lval \z^{(n+1)}_k - \z^{(n+2)}_k \rval & \leq & \frac{1}{\lval \alpha \rval^2 \lval \hz_k - \omep \rval} \cdot \max \{  
\lval \z^{(n)}_{k+1} - \z^{(n+1)}_{k+1} \rval,  \lval \z^{(n)}_1 - \z^{(n+1)}_1 \rval \}\\
& = & \frac{\rho_k}{\rho_{k+1}} \max \{  
\lval \z^{(n)}_{k+1} - \z^{(n+1)}_{k+1} \rval,  \lval \z^{(n)}_1 - \z^{(n+1)}_1 \rval \}
\end{eqnarray*}
Hence, for all $n$, there exist  $k_1, \dots, k_j$ such that $1 \leq k_i \leq \hp -1$, $k_1 + \cdots + k_j =n$,  and
$$\lval \z^{(n+1)}_{1} - \z^{(n+2)}_{1} \rval \leq \frac{\rho_1}{\rho_{k_1 + 1}} \cdots \frac{\rho_1}{\rho_{k_j + 1}} \max \{ 
\lval \z^{(0)}_{k_j+1} - \z^{(1)}_{k_j+1}\rval , \lval \z^{(0)}_1 - \z^{(1)}_1 \rval \}.$$
Now let
$$\lambda = \max \left\{ \frac{\rho_1}{\rho_k} \,\, /\, k=2, \dots, \hp-1 \right\} <1.$$
Since, $j \geq \frac{n}{\hp-1}$, it follows that
$$\lval \z^{(n+1)}_{1} - \z^{(n+2)}_{1} \rval \leq \lambda^{ \frac{n}{\hp-1} } \| \vec{\z}^{(0)} - \vec{\z}^{(1)} \|_\infty.$$
Therefore, $\z^{(n)}_1$ converges to some $\z^{(\infty)}_1$ as $n \rightarrow \infty$ 
and $\z^{(n)}_{k+1} \rightarrow \z^{(\infty)}_{k+1} = \psi^k_{\z^{(\infty)}_1}(\z^{(\infty)}_1)$ for $k=1, \dots, \hp-1$.
It follows that $(\z^{(\infty)}_1, \dots, \z^{(\infty)}_\hp)$ is a fixed point for $T$.
\hfill
$\Box$

\begin{lemma}
$T$ has a unique fixed point in $\B$.  
\end{lemma}

\noindent
{\bf Proof.}
Suppose that $\vec{\z}= (\z_1, \dots, \z_\hp)$ and 
$\vec{\eta}=(\eta_1, \dots, \eta_\hp)$ are fixed points of $T$. For $k=1, \dots, \hp-1$
the polynomial $h: \z \mapsto \alpha^2 (\z -1)^2 (\z +2)$  maps $B_{\rho_k}(\z_k) = B_{\rho_k}(\eta_k)$ isomorphically onto its image
and $$\lval \frac{dh}{d\z} (\z_k) \rval = \lval \alpha \rval^2 \lval \z_k - \omep \rval.$$
Therefore, 
\begin{eqnarray*}
  \lval \z_{k+1} - \eta_{k+1} + \eta_1 - \z_1 \rval & =& \lval h(\z_k) - h(\eta_k) \rval \\
& = & \lval \z_k - \eta_k \rval  \lval \alpha \rval^2  \lval \z_k - \omep \rval.
\end{eqnarray*}
Hence, 
\begin{eqnarray*}
  \lval  \z_{k} - \eta_{k}  \rval & \leq & \frac{1}{\lval \alpha \rval^2  \lval \z_k - \omep \rval} \max \{ \lval  \z_{k+1} - \eta_{k+1} \rval, \lval  \eta_1 - \z_1\rval \} \\
& = & \frac{\rho_k}{\rho_{k+1}} \max \{ \lval  \z_{k+1} - \eta_{k+1} \rval, \lval  \eta_1 - \z_1\rval \}.
\end{eqnarray*}
It follows that for some $k$ such that $2 \leq k < \hp$:
$$\lval \z_1 - \eta_1 \rval \leq \frac{\rho_1}{\rho_k} \lval \eta_1 - \z_1 \rval.$$
Since $\rho_1 < \rho_k$ we have that $\z_1 = \eta_1$ and $\vec{\z}= \vec{\eta}$.
\hfill $\Box$

\subsection{Parameter balls}
\label{paraballs-ss}
Here  we prove simultaneously 
propositions~\ref{level-n-p} and~\ref{level-n+1-p} with their analogues for the family $\psinu$.

\begin{proposition}
\label{psi-level-n-p}
   Let $\pb^\psi_n$ be a $\psi$--parameter ball of level $n$.
Then the following hold:

(i) $\pbpsi_n = \dbnu_n (\nu) $ for all $\nu \in \pbpsi_n$.

(ii) $\bM^\nu_{n+1} (\omep) = \bM^{\nup}_{n+1} (\omep)$ for all $\nu, \nup \in \pbpsi_n$.

(iii) There exists a unique center of level $n$ 
in $\pbpsi_n$. The period of this center is $\min \{ k \geq 1  \,\, /\, M^\nu_{n+1-k, k} (\omep) =1 \}$ for any $\nu \in \pbpsi_n$.
\end{proposition}

\noindent
{\bf Proof of Propositions~\ref{level-n-p} and~\ref{psi-level-n-p}.} 
In view of Lemma~\ref{coordinates-l}, we just have to prove the above proposition concerning the family $\psinu$.
By Lemma~\ref{apriori2-l} (i) we have that  $\dbnu_n (\nu) \subset \pbpsi_n$
for all $\nu \in \pbpsi_n$. Moreover, for all $\nu, \nup \in \pbpsi_n$, the dynamical balls $\dbnu_n (\nu)$ and $\dbnup_n (\nup)$
are equal or disjoint. 
By Proposition~\ref{fixed-p}, for each $\nu \in \pbpsi_n$ the ball $\dbnu_n (\nu)$ contains at least one element of the finite set
$$\{ \nu \,\, / \, \psinu^k (\omep) = \omep \mbox{ for some } 1 \leq k \leq n+2 \}.$$
It follows that $\pbpsi_n = \cup_{\nu \in \pb_n} \dbnu_n (\nu)$ is a finite union of closed and pairwise disjoint balls.
This is only possible if $\pbpsi_n = \dbnu_n(\nu)$ for all $\nu \in \pbpsi_n$.  Hence we have proven statement (i).
Statement (ii) now follows from Lemma~\ref{apriori2-l} (i). 

For (iii), let $\hnu \in \pbpsi_n$ and note that $\hnu \in \dbhnu_n (\hnu) \subset \dbhnu_{n-1}(\hnu)$. Let $\hp = \min \{ k \geq 1 \,\, /\, 
\omep \in \dbhnu_{n+1-k} (\psihnu^k (\omep)) \}$. Proposition~\ref{fixed-p} says that there exists a unique $\nup$ 
such that $\psinup^\hp (\omep) = \omep$ and $\lval \nup - \hnu \rval < r^\hnu_{n-1} (\hnu)$. We must show that $\nup \in \pbpsi_n$.
In fact, since $\omep \in \psihnu^{\hp-1} (\dbhnu_{n} (\hnu))$ there exist a level $n+1$ ball $\dbhnu_{n+1}(\hz_1)$ in 
$\dbhnu_{n} (\hnu)$ that maps onto $\dbhnu_1 (\omep)$ under $\psihnu^{\hp-1}$. By Proposition~\ref{fixed-p}  we have that
$\lval \nup - \hz_1 \rval < r^\hnu_n (\hnu)$. Therefore $\nup \in \pbpsi_n$.
\hfill $\Box$

\begin{corollary}
\label{trivial-paraends-c}
  If $\omep \in \filled( \psinu)$ and 
$\bM^\nu (\omep)$ is not periodic and $\pbpsi_n$ is the level $n$ parameter ball containing $\nu$, 
then $\{ \nu \} = \cap \pbpsi_n$. 
\end{corollary}

\begin{proposition}
\label{psi-level-n+1-p}
  Consider a level $n$ parameter ball  $\pbpsi_n$ and let $P$ be an element of the affine partition associated to $\pbpsi_n$. 
For any $\nu \in \pbpsi_n$ we have the following:

There exists a level $n+1$ parameter ball contained in $P$ if and only if there exists a level $n+1$ dynamical ball 
$\dbnu_{n+1} (\z)$ contained in  $P$.
In this case, the level $n+1$ parameter ball $\pbpsi_{n+1}$ is unique and 
$$\bM^\nu_{n+1} (\z) = \bM^{\nup}_{n+1} (\nup)$$
for all $\nup \in \pbpsi_{n+1}$. In particular, the radii of  $\dbnu_{n+1} (\z)$ and $\pbpsi_{n+1}$ coincide.
\end{proposition}

\noindent
{\bf Proof of  Propositions~\ref{level-n+1-p} and ~\ref{psi-level-n+1-p}.}
As in the previous proof, we just have to prove the above proposition concerning the family $\psinu$.
By Proposition~\ref{fixed-p}, if an element $P$ of the partition $\cP_n$ associated to $\pbpsi_n$ contains
a level $n+1$ dynamical ball of some $\psinu$ with $\nu \in \pbpsi_n = \dbnu_n (\nu)$, 
then $P$ contains a level $n+1$ parameter ball.
Conversely, if $P$ contains a parameter ball $\pbpsi_{n+1}$ of level $n+1$ and $\nup \in \pbpsi_{n+1}$, then
$\psi^n_{\nup} (P) = \psi^n_\nu (P)$ for all $\nu \in \pbpsi_n$. Moreover, since $\psinup^n (\nup)$ is a level $1$ point
we have that $\psinup^n (P) = B_1 (\omep)$ or $B_1 (-2)$. In either case the preimage of $\db_1(\omep)$ or
$\db_1(-2)$ under $\psinu : P \rightarrow B_1(\omep)$ or $\psinu : P \rightarrow B_1(-2)$ is a level $n+1$ dynamical
ball contained in $P$. 

We must show that $P$ contains at most one parameter ball of level $n+1$.
Suppose that $\pb_{n+1}$ and ${\pb}^\prime_{n+1}$ are level  $n+1$  parameter balls contained in $P$. 
From Proposition~\ref{fixed-p}, 
it is sufficient to  show that the  periods $p$ and $p^\prime$  of their centers $\mu$ and $\mup$ coincide. 
By Lemma~\ref{apriori2-l} we have that $\psi^k_\mu (P) = \psi^k_\mup (P)$ for $k=1, \dots, n$.
Moreover, since there is at most one dynamical ball inside $\psi^{k-1}_\mu(P)$ and 
$ \db^\mu_{n+1}(\mu) = \pb_{n+1} \subset P$
we have that:
$\omep \in \psi^k_\mu(\db^\mu_{n+2} (\omep)) =  \psi^{k-1}_\mu (\db^\mu_{n+1}(\mu)) \subset \psi^{k-1}_\mu(P)$ if and only if
$\omep \in \psi^{k-1}_\mu(P)$. 
Similarly, $\omep \in \psi^k_\mup (\db^\mup_{n+2} (\omep))$ if and only if
$\omep \in \psi^{k-1}_\mup (P)$. 
Therefore, $p = p^\prime$.

Now let $P$ be an element of the affine partition of $\pbpsi_{n}$ and let $\nu \in \pbpsi_n$.
Suppose that there exists a level $n+1$
parameter ball $\pb_{n+1} \subset P$ and a level $n+1$ dynamical ball $\dbnu_{n+1} (\z) \subset P$.
To complete the proof of the proposition, given $\nup \in \pb_{n+1}$ we must show that 
$$\bM^\nu_{n+1} (\z) = \bM^{\nup}_{n+1} (\nup).$$
Since $\z \in \dbnu_n(\nu)$ we have that 
$\bM^\nu_n (\z) = \bM^\nu_n(\nu) = \bM^{\nup}_n (\nup)$.
Thus we just need to prove that $\omep \in \dbnu_{n+1-k} (\psinu^k(\z))$ if and only if 
$\omep \in \dbnup_{n+1-k}(\psinup^k (\nup))$, for all $k=0, \dots, n+1$.
Since $\nup \in \dbnu_n(\nu)$, Lemma~\ref{apriori2-l} implies that 
$\psi_{\nu*} = \psi_{\nup*}: \cP_{n+1-k} \rightarrow \cP_{n+1-(k+1)}$ for all $k$ such that  $1 \leq k \leq n$, where
$\cP_{n+1-k}$ is the affine partition associated to $\dbnu_{n+1-k}(\psinu^k (\omep))=\dbnup_{n+1-k}(\psinup^k (\omep))$.
Now $\omep \in  \dbnu_{n+1-k} (\psinu^k(\z))$ if and only if $\omep \in \psinu^k(P)$ 
and  $\omep \in  \psinup^k(P)$ if and only if   $\omep \in \dbnup_{n+1-k}(\psinup^k (\nup))$. 
Therefore the proposition follows from the fact that $\psinu^k(P)=\psinup^k(P)$.
\hfill $\Box$

\subsection{Realization}
\label{realization-ss}
In order to prove Theorem~\ref{parameter-th} we first have to show that 
 every admissible critical marked grid of level $n$  is realized by a cubic polynomial.

\begin{proposition}
\label{realization-p}
  Let $n \geq 0$ and  $\bM_{n+1}$ be an admissible critical  marked grid of level $n+1$.
Then there exists $\nu \in \cEpsi_n$ such that  $\bM^\nu_{n+1} (\omep)= \bM_{n+1}$.
\end{proposition}

\noindent
{\bf Proof.}
Since the proposition is clearly true for $n=0$ we proceed by induction. That is we suppose that $\bM_{n+2}$ is an admissible 
critical marked grid of level $n+2$ and 
$\nu \in \cEpsi_n$ is  such that $\bM^\nu_{n+1} (\omep)$ coincides with $\bM_{n+2}$ 
in all the positions $(\ell,k)$ with  $\ell + k \leq n+1$. By Proposition~\ref{fixed-p}, it suffices to show that
there exists a level $n+1$ point $\z$ contained in $\dbnu_{n}(\nu)$ such that $\bM^\nu_{n+1}(\z)$ coincides with the
grid obtained from $\bM_{n+2}$ after erasing its first column. For this purpose let $\nu_j = \psinu^j(\omep)$, 
$p$ be the minimal $k\geq1$ such that the $(n+2-k,k)$ position of $\bM_{n+2}$ is marked and
$k_1, \dots, k_m$ be such that $0< k_1 < \cdots < k_m < p$ and $(n+1-j,j)$ is marked in $\bM_{n+2}$ if and only if
$j = k_i$ for some $i$. 

Our task boils down to find a dynamical ball $\db_{n+1}$ of level $n+1$ contained in $\dbnu_{n}(\nu)$ such that
$\omep \notin \db_{n+1} \cup \cdots \cup \psinu^{p-2} (\db_{n+1})$ and $\psinu^{p-1}(\db_{n+1}) = \db_{n+2-p}(\omep)$.

We first claim that there exist $2$ level $n+2 - k_m$ dynamical 
balls $B^0_m$ and $B^1_m$  contained in $\dbnu_{n+1-k_m}(\nu_{k_m})$
such that for $i=0,1$ we have that $\omep \notin B^i_m \cup \cdots \cup \psi^{p-1-k_m}( B^i_m)$ and
$\omep \in \psi^{p-k_m}(B^i_m)$. 
There are two cases according to whether $p < n+2$ or $p=n+2$.

In the case that $p < n+2$, by (Mc) of Proposition~\ref{rules-p}, $g_m = \psinu^{p-k_m}: \dbnu_{n+1-k_m}(\nu_{k_m}) \rightarrow 
\dbnu_{n+1-p}(\nu_p)$ has degree $2$ and $g_m(\omep) \notin \dbnu_{n+2-p}(\omep)$. Hence, in this case we let $B^0_m$ and 
$B^1_m$ be the two preimages  of $\dbnu_{n+2-p}(\omep)$ under $g_m$.

In the case that $p= n+2$, by (Md) of Proposition~\ref{rules-p}, the position $(1,n+1-k_m)$ in $\bM_{n+2}$ is marked
(otherwise, taking $\ell =n+1-k_m$ and $k=k_m$ in (Md) we would have that $M_{1,n+1}$ would be marked and therefore
$p=n+1$).
Hence, $g_m = \psinu^{n+1-k_m} : \dbnu_{n+1-k_m} (\nu_{k_m}) \rightarrow \db_0$ has degree $2$ and $g_m(\omep) \in \db_1(\omep)$.
Therefore, in this case, we let $B^0_m$ and $B^1_m$ be the two preimages of $\db_1(-2)$ under $g_m$.

We now claim that for all $j=1, \dots, m$ there exist at least $2$ level $n+2-k_j$ dynamical
balls $B^0_j$ and $B^1_j$  contained in 
$\dbnu_{n+1-k_j}(\nu_{k_j})$ such that for $i=0,1$ we have that $\omep \notin B^i_j \cup \cdots \cup \psinu^{p-1-k_j}( B^i_j)$ and
$\omep \in \psinu^{p-k_j}(B^i_j)$. In fact, for $j=m$ we have already established this so we may assume the above true 
for $j+1$ and prove it for $j$. Since 
$g_j = \psinu^{k_{j+1}-k_j}: \dbnu_{n+1-k_j}(\nu_{k_j}) \rightarrow \dbnu_{n+1-k_{j+1}}(\nu_{k_{j+1}})$
has degree $2$, it follows that  at least one of the two balls $B^0_{j+1}$, $B^1_{j+1}$ does not contain $g_j (\omep)$, say 
$B^0_{j+1}$, and we may let $B^0_j$ and $B^1_j$ be the preimages of $B^0_{j+1}$ under $g_j$.

Finally, note that $g_0= \psinu^{k_1-1} : \dbnu_{n}(\nu) \rightarrow \db_{n+1-k_1}(\nu_{k_1})$ is one-to-one. The preimage
of $B^0_1$ under $g_0$ is a level $n+1$ dynamical ball $\db_{n+1}$ contained in $\dbnu_{n}(\nu)$ and such that
$\omep \notin \db_{n+1} \cup \cdots \cup \psinu^{p-2}(\db_{n+1})$ and $\omep \in \psinu^{p-1}(\db_{n+1})$. 
\hfill $\Box$ 

\medskip
Again we simultaneously proof Theorem~\ref{parameter-th} and the corresponding version for
the family $\psinu$.

\begin{theorem}
  \label{psi-parameter-th}
  Consider an admissible  critical marked grid $\bM$ and
let $$ C^\psi_{\bM} := \{ \nu \in \cEpsi_0  \,\,/\, \omep \in \filled(\psinu) \mbox{ and }  \bM = \bM^\nu (\omep) \}$$
Then the following hold:

(i) If $\bM$ is periodic, then $C^\psi_{\bM}$ is a non-empty union of finitely many closed and pairwise disjoint balls.

(ii) If $\bM$ is not periodic, then $C^\psi_{\bM}$ is a non-empty compact set and 
$$C^\psi_{\bM} \subset \partial \{ \nu \in \cEpsi_0 \,\, /\, \{ \omega^\pm = \pm 1 \} \subset \L \setminus \filled (\psinu) \}.$$
\end{theorem}

\noindent
{\bf Proof of Theorems~\ref{parameter-th} and~\ref{psi-parameter-th}.} 
In view of Lemma~\ref{coordinates-l} it is sufficient to prove Theorem~\ref{psi-parameter-th}.
Let $M_{\ell, k}$ denote the $(\ell,k)$ entry of $\bM$. 
For $n \geq 0$, let $\bM_n$ denote the level $n$ grid with entries $M_{\ell,k}$ where $\ell +k \leq n$.

(i) Suppose that $\bM$ is periodic of period $p$ . Let $n_0$ be such that $\bM_{n, n-k}$ is unmarked for all
$1 <k <p$ and all $n > n_0$. Denote by $\nu_1, \dots, \nu_m$ the level $n_0$ centers with critical marked
grid of level $n_0 +1$ that coincides with $\bM_{n_0 +1}$. It follows that $\bM = \bM^\nu_i (\omep)$ for all $i=1, \dots, m$.
For each $i$, let $X_i$ be the infraconnected component of $\filled (\psi_{\nu_i})$ that contains $\nu_i$.
Recall that $X_i$ is a closed ball (Theorem~\ref{wandering-t}). 
Moreover, if $\pb_n (\nu_i)$ denotes the level $n$ parameter ball which contains
$\nu_i$, then $X_i = \cap \pb_n (\nu_i)$. Hence, $\bM^\nu (\omep)= \bM^{\nu_i}(\omep) = \bM$ for all $\nu \in X_i$.
That is, $X_1 \cup \cdots \cup X_m \subset C^\psi_\bM$.  Given $\nu \in C^\psi_\bM$, to complete the proof of (i),
 it is sufficient to show that $\nu \in X_i$ for some $i$. In fact, for $n > n_0$,
the center $\nu_c$ of the parameter ball $\pb_n (\nu)$ of level $n$ containing $\nu$ is such that
$\bM^{\nu_c}_{n+1} (\omep) = \bM^\nu_{n+1} (\omep) = \bM_{n+1}$. Therefore $\nu_c = \nu_i$ for some $i$
and $\nu \in X_i$.

\smallskip
(ii) Let $Y_n$ be the union of all level $n$ parameter balls $\pb_n$ such that $\bM_{n+1} = \bM^\nu_{n+1}(\omep)$
for all $\nu \in \pb_n$. It follows that $C^\psi_\bM = \cap_{n \geq 0} Y_n$. The radius $r_n$ of the parameter balls that
participate in $Y_n$ coincide since it only depends on $\bM_{n+1}$. Moreover $r_n \rightarrow 0$ as $n \rightarrow \infty$.
Therefore $C^\psi_\bM$ is a compact non-empty set. For all $\nup \in C^\psi_\bM$ the level $n$ parameter ball $\pb_n (\nup)$
which contains $\nup$ has non-empty intersection with $\cEpsi_n \setminus \cEpsi_{n+1}$. If 
$\nu \in \cEpsi_n \setminus \cEpsi_{n+1}$, then both critical points of $\psinu$ escape to infinity. It follows that
$\nup \in \partial \{ \nu \in \cEpsi_0 \,\, /\, \{ \omega^\pm = \pm 1 \} \subset \L \setminus \filled (\psinu) \},$ since 
$\{ \nup \} = \cap \pb_n (\nup)$.
\hfill $\Box$

\subsection{Proof of Theorem~\ref{compact-it} and a corollary}
\label{proof23-ss}
We will need the following result:

\begin{lemma}
  \label{open-l}
 If $\bM$ is  a periodic admissible critical marked grid, then
$$S_\bM = \{ \vphiab \in \cE^+(\L) \,\,/\, \bM^{\ab}(\alpha) = \bM \}.$$
is open in $\cubicl \equiv \L^2$.
\end{lemma}

\noindent
{\bf Proof.}
Let $\vphi_{\alpha_0, \beta_0} \in S_\bM$ and denote by $B$ the infraconnected component of $\filled (\vphi_{\alpha_0, \beta_0} )$
which contains $\alpha$. Hence, $B$ is periodic under $\vphi_{\alpha_0, \beta_0}$, say of period $p$.
Let $r$ be the radius of $B$. Then, for a sufficiently small neighborhood $V$ of $(\alpha_0, \beta_0) $,
$$\lval \vphi^p_{\alpha_0, \beta_0} (\zeta) - \vphi^p_{\alpha, \beta} (\z) \rval < r$$
for all $\zeta$ in a closed ball   containing $B$.  It follows that $\vphiab (B) = B$ 
 and therefore
$\bM^\ab (\alpha) $ is periodic for all $(\ab) \in V$.
Similarly, for any fixed level $n$, the  dynamical balls of level $n$ of $\vphiab$ are locally constant.
It follows that  there exists $n_0$ such that, for all $(\alpha, \beta) \in V$, the periodic critical marked grid 
$\bM^\ab(\alpha)$ is uniquely determined by which dynamical balls of level $n_0$ contain the critical point $\alpha$.
Therefore, after shrinking $V$, if necessary, $\bM^\ab (\alpha) = \bM$ for all $(\ab) \in V$ 
\hfill $\Box$

\medskip
\noindent
{\bf Proof of Theorem~\ref{compact-it}.}
By Corollary~\ref{almost2-c}, it is sufficient to show that $\vphiab \in \partial \shiftcubicl$ if and only if
$\vphiab \in \escapepml$ and the corresponding critical marked grid is not periodic.
From the previous lemma we obtain that if $\vphiab \in \partial \shiftcubicl$, then 
$\vphiab \in \escapepml$ and the corresponding critical marked grid is not periodic.
Conversely, Theorem~\ref{parameter-th} (ii) says that polynomials $\vphiab \in \escapepml$ with 
aperiodic critical marked grid are in $\partial \shiftcubicl$.
\hfill $\Box$

\begin{corollary}
\label{hyperbolic-c}
Let $\cH \subset \cubicl$ be the set formed by all the cubic polynomials $\vphiab$ such that the Julia set 
$\julia (\vphiab)$ is critical point free (i.e., $\{ + \alpha, -\alpha \} \cap \julia (\vphiab) = \emptyset$). 
Then $\cH$ is open and dense in $\cubicl$.
\end{corollary}

\medskip
\noindent
{\bf Proof of Corollary~\ref{hyperbolic-c}.}
First note that the polynomials in $\escapepml$ with aperiodic critical marked grid
are in $\partial \shiftcubicl$ (Theorem~\ref{parameter-th} (ii)). Hence, $\cH$ which is  the union
 of $\coneccubicl$, $\shiftcubicl$ and
the polynomials in $\escapepml$ that have periodic critical marked grid, is open and dense.
\hfill $\Box$

\subsection{Dynamics over finite extensions of $\laut$}

The aim of this subsection is to discuss the dynamical behavior of cubic polynomials with coefficients 
in a finite extension of $\laut$. In particular, we
show that some cubic polynomials with coefficients in a finite extension of $\laut$
have a non-periodic recurrent critical point.
That is, a critical point which is an accumulation point of its orbit but it  is not periodic.
Examples of non-archimedean dynamical systems over finite extensions of $\Q_p$ with wild recurrent critical points were 
recently given by Rivera in~\cite{rivera-04}. We emphasize that Rivera shows the existence of a {\em wild} recurrent critical
point (i.e., a critical point where the local degree is a multiple of $p$).

If $\psianu$ is such that $\omep \in \filled(\psianu)$ and $\omega^-$ escapes to infinity, then
$\omep$ is recurrent if and only if the critical marked grid of $\psianu$ is not periodic and has marked columns of arbitrarily
long depth. There is a stronger notion of critical recurrence associated to marked grids called {\em persistent recurrence}
(e.g., see~\cite{harris-99,milnor-00}).

As a corollary of our description of the parameter space of cubic polynomials we will be able to use the examples of recurrent and 
persistently recurrent critical marked grids due to Harris~\cite{harris-99} to prove the following:

\begin{corollary}
  \label{recurrent-c}
Let $\psi_{ \nu} (\z) = t^{-2} (\z-1)^2 (\z+2) + \nu$.  Then the following hold:

a) There exists $\nu_a \in \laut$ such that the critical point $\omep = +1$ is recurrent and not periodic under iterations of 
$\psi_{\nu_a}$.

b) There exists $\nu_b \in \Q^a((t^{1/2}))$ such that critical point $\omep = +1$ is  persitently 
recurrent and not periodic under iterations of  $\psi_{\nu_b}$.
\end{corollary}

The corollary will follow from two lemmas and the work of Harris cited above. 

\medskip
In order to simplify notation, for $\z \in \pui$ we say that the {\em algebraic degree of $\z$} is
$$\delta(\z) := [ \laut (\z) : \laut].$$
When $\pui$ is regarded as the inductive limit of $\{ \lautm; m \in \N \}$ the algebraic degree
of $\z$ coincides with the smallest $m$ such that $\z \in \lautm$. 

\smallskip
Let us fix $\alpha \in \pui$ and consider as before  $\psi_\nu(\z) =\psi_{\alpha, \nu} (\z) = \alpha^{-2} (\z -1 )^2 (\z +2) + \nu$
with $\nu \in \S$. 
We define the {\em algebraic degree of $\psi_\nu$} as
$$\delta(\psinu) := \max \{ \delta(\alpha), \delta(\nu) \} \quad \quad \mbox{if } \nu \in \pui$$
and $\infty$ otherwise. 
The {\em algebraic degree of a ball $B \subset \S$} is defined as:
$$\delta (B) := \min \{ \delta (\psi_\nu) \,\, /\, \nu \in B \}$$

We will be interested in computing the algebraic degree of the $\psi$-parameter balls in the $\alpha$-slice of 
the $\psi$-parameter space (see Subsection~\ref{change-ss}). 
Clearly $\delta (\pb_0) =\delta(\alpha)$. Our next result shows that the center of a parameter ball minimizes
the algebraic degree of the elements of the ball. More precisely:

\begin{lemma}
  \label{center-degree-l}
  Let $\pb_{n+1} \subset \pb_{n}$ be parameter balls of levels $n+1$ and $n$ respectively with $n \geq 0$.
Denote by $\nu_{n+1}$ the center of $\pb_{n+1}$ and by $P_n$ the element of the affine partition of $\pb_n$ that
 contains $\pb_{n+1}$.  Then:
$$ \delta(\psi_{\nu_{n+1}})= \delta(\pb_{\nu_{n+1}}) = \delta(P_n).$$
\end{lemma}

Before proving the lemma let us remark that if  $\vphi$ is a polynomial with coefficients in $\lautm$, 
and $B, B^\prime  \subset \S$ are balls such that  $\vphi : B \rightarrow B^\prime$ is bijective. Then 
$\delta (\z) \leq \max \{m , \delta(\vphi(\z)), \delta(\nu) \}$ for all $\nu \in B$. 
This easily follows from the Newton polygon of $\vphi(\cdot - \nu) - \vphi(\z)$. 

\medskip
\noindent
{\bf Proof.}  Let $\hnu \in P_n$ be such that $\delta(\psihnu) = \delta(P_n)$.
Let $\hp$ be the period of the center $\psi_{\nu_{n+1}}$. By Proposition~\ref{psi-level-n+1-p} there exists 
a level $n+1$ dynamical ball $\dbhnu_{n+1}$ contained in $P_n$ such that $\psi^{\hp-1}_{\hnu} : \dbhnu_{n+1} \rightarrow 
\dbhnu_{n+2-\hp}(\omep)$ is one--to--one. Hence there exists 
a unique $\hz_1 \in \dbhnu_{n+1}$ such that $\psi^{\hp-1}_{\hnu}(\hz_1) = \omep$. It follows that $\delta(\hz_1) \leq \delta(\psihnu)$.
Let $\hz_k = \psihnu^{k-1}(\hz_1)$ and consider $\B$ as in Lemma~\ref{thurston-l}. From this lemma we conclude that there is 
a well defined Thurston map $T: \B \rightarrow \B$. If $T(\z_1, \dots , \z_\hp) = (\zp_1, \dots, \zp_\hp)$, then
$\delta(\zp_i) \leq \max\{ \delta(\zp_{i+1}), \delta(\psi_{\z_1}) \}$ since $\zp_i$ is obtained as a preimage of $\z_{i+1}$ 
under the restriction of $\psi_{\z_1}$ to a ball where this polynomial is injective. The first coordinates of 
$T^k (\hz_1, \dots \hz_\hp)$ converge to $\nu_{n+1}$ as $k \rightarrow \infty$. Therefore, 
$$\delta(\psi_{\nu_{n+1}}) \leq \delta(\psi_{\hz_1}) \leq \delta(\psihnu) = \delta (P_n).$$
The lemma easily follows.
\hfill $\Box$

\medskip
Next we show that $\delta(\pb_n)$ is in fact computable from the information contained in the level $n+1$ critical marked
grid of the parameters in $\pb_n$:

\begin{lemma}
  \label{delta-comp-l}
  Let $\pb_{n+1} \subset \pb_{n}$ be parameter balls of levels $n+1$ and $n$ with centers $\nu_{n+1}$ and $\nu_n$ respectively.
Denote by $r_n$ the radius of $\pb_n$ and let $s_n = \min \{ s \in \N \,\,/ \, s | \log r_n | \in \N \}$. Then the following hold:

a) If $\nu_{n+1} = \nu_n$, then $\delta(\pb_{n+1}) = \delta(\pb_{n})$.

b) If $\nu_{n+1} \neq \nu_n$, then $\delta(\pb_{n+1}) = \max \{ s_n , \delta(\pb_n) \}$.
\end{lemma}

\noindent
{\bf Proof.}
Part a) follows immediately from the previous lemma. Suppose that $\nu_{n+1} \neq \nu_n$ and observe that
by Proposition~\ref{psi-level-n+1-p}, $\lval \nu_{n+1} - \nu_n \rval = r_n$. 
Hence $\order (\nu_{n+1} -\nu_n) = -\log r_n$. It follows that 
$$s_n \leq \max\{\delta(\nu_{n+1}), \delta(\nu_n) \} 
           \leq \max\{ \delta(\psi_{\nu_{n+1}}), \delta(\psi_{\nu_n}) \}
           =    \delta(\psi_{\nu_{n+1}}).$$
Therefore, $\max\{ s_n, \delta(\psi_{\nu_n}) \} \leq \delta(\psi_{\nu_{n+1}})$.

Now let  $P_n$ be the ball of the affine partition of $\pb_n$ that contains $\pb_{n+1}$.
It follows that $P_n$ contains a series of the form $\nu_n + a t^{q/s_n}$ for some $0 \neq a \in \Q^a$ and some $q \in \N$
relatively prime with $s_n$. 
So $$\delta(\psi_{\nu_{n+1}}) = \delta(P_n) \leq \max\{ \delta(\psi_{\nu_n}), s_n \}$$
which, in view of the previous lemma, finishes the proof.
\hfill $\Box$

\medskip
Let us now illustrate how the above lemmas may be used to compute the algebraic degree of some parameters.
For simplicity we restrict to the case in which $\alpha = t^{-1}$ and let $\psi_\nu$ be a polynomial with critical marked grid
$\mathbf{M}= (M_{\ell,k})$. For all $n$, we denote by $\pb_n$ the parameter ball of level $n$ containing $\nu$,
by $r_n$ its radius and by $\nu_n$ its center. 
Furthermore we suppose that 
$\omep$ is periodic of period $p$ under $\psi_\nu$ if the critical marked grid is periodic of period $p$.
Now the smallest integer $s_k$ such that 
$$s_k | \log r_n | = s_k \cdot 2 \cdot \sum^{k+1}_{\ell=1 } 2^{\sum^{\ell-1}_{i=0} {{M}_{\ell-i,i}}} \in \N$$
is clearly computable from $\mathbf{M}$. 
The previous lemma implies that
$$\delta(\pb_{n+1}) = \max \{ s_k \,\, / \nu_{k+1} \neq \nu_k, \, 0 \leq k \leq n \}.$$
Moreover, $\nu$ is algebraic over $\laut$ if and only if
$$\delta(\mathbf{M}) = \sup \{ s_k \,\, /\, \nu_{k+1} \neq \nu_k \} < \infty$$
and in this case the algebraic degree of $\psi_\nu$ coincides with $\delta(\mathbf{M})$.

The above formula for $\delta(\psi_\nu)$ coincides  with Branner and Hubbard's formula for the ''length'' of a ''turning curve''
passing through a complex cubic polynomial with critical marked grid $\mathbf{M}$. 
In~\cite{harris-99},  Harris shows the existence of critical marked grids satisfying rules (Ma) through (Md)
which are (resp. persistently) critically recurrent and aperiodic such that $\delta(\mathbf{M})=1$  (resp. $\delta(\mathbf{M})=2$).
Corollary~\ref{recurrent-c} now follows.

\section{Complex cubic polynomials}
\label{complex-s}
Recall that $\cubicc$ denotes the space of monic centered critically marked cubic polynomials.
That is, polynomials of the form:
$$g_{a,b} (z) = z^3 - 3 a^2 z +b$$
where $(a,b) \in \C^2$.
To prove Theorem~\ref{complex-it} is convenient to change coordinates and work in the space $\hcubicc$ of polynomials
of the form
$$\fav (z) = z^3 - 3a^2z + 2 a^3 +v$$
where $(a,v) \in \C^2$. Thus we identify $\hcubicc$ with $\C^2$. 
The critical points of $\fav$ are also  $\pm a$ but here we have the advantage
that  $v = \fav(+a)$ is a critical value (compare with~\cite{milnor-91}).
Observe that $$\fav = g_{a, 2a^3 +v}.$$
Moreover,  instead of working with the family $\vphib (\z) = \z^3 - 3 t^{-2} \z + \beta \in \S[\z]$ with $\beta \in \S$ it 
is easier to work with 
$$\psinu (\z) = t^{-2} (\z -1)^2 (\z +2) + \nu$$
where $\nu \in \S$. Now for $\beta = 2 t^{-3} + t^{-1} \nu$ we have that  $\psinu(\z) = t \vphib (t^{-1} \z)$.
The critical points of $\psinu$ are at $\omega^\pm= \pm 1$. 

The sets $A_\C$ and $A_\S$ of the introduction correspond to 
$$\hA_\C := \{ (a,v) \in \C^2 \,\,/\, +a \in \filled (\fav) \mbox{ and } C_{a,v} (+a) \mbox{ is not periodic} \}.$$
where $C_{a,v}(+a)$ denotes the connected component of $\filled (\fav)$ which contains $+a$. Also
$$\hA_\S := \{ \nu \in \S \,\,/\, \omep \in \filled (\psinu) \mbox{ and } IC_\nu (\omep)  \mbox{ is not periodic} \}.$$
where $IC_\nu(\omep)$ is the infraconnected component of $\filled (\psinu)$ which contains $\omep$.

It is easy to verify that Theorem~\ref{complex-it} in this new coordinates is equivalent to:

\begin{theorem}
  \label{complex-th}
  There exists $\epsilon >0$ such that for all $\nu = \sum_{\lambda \geq 0} a_\lambda t^\lambda$
in $\hA_\S$ the series $$ \sum_{\lambda \geq 0} a_\lambda  \mathe^{2 \pi i T \lambda}$$
is the Fourier series of an analytic almost periodic function $\bv_\nu : \hepsilon \rightarrow \C$.
Moreover, 
$$  \begin{array}{rccc}
  \tilde{\Psi} : & \hepsilon \times \hA_\S & \rightarrow &  \hA_\C \cap \{ | a | > 1/\epsilon \} \\
                      & (T, \nu)                           & \mapsto      & (\mathe^{-2 \pi i T}, \mathe^{-2 \pi i T} \bv_\nu(T))
  \end{array}$$
is a well defined and onto  map which is continuous in $(T,\nu)$ and holomorphic in $T$.
Furthermore, $\tilde{\Psi}$ projects to a homeomorphism:
$$  
  {\Psi} :  \hepsilon \times A^\psi_\S  / \sim   \rightarrow   \hA_\C \cap \{ | a | > 1/\epsilon \} 
$$
where $\sim$ is the smallest equivalence relation that identifies $(T-1, \nu)$ with $(T, \sigma (\nu))$.
\end{theorem}

The parameter space $\hcubicc$ is also stratified according to how many critical points escape to $\infty$.
The connectedness locus $\hconeccubicc$ is the set of polynomials $\fav$ with connected Julia set.
Here $\hconeccubicc$ is also compact, connected and cellular (see~\cite{branner-88}). 
The shift locus $\hshiftcubicc$ is the
set of polynomial $\fav$ with all its critical points escaping and 
$$\hescapepmc = \{ \fav \in \hcubicc \,\,/\, \pm a \in \filled (\fav) \ni \mp a \}.$$

\subsection{The combinatorics of complex cubic polynomials}
Important tools to study the dynamics of complex polynomials are the Green function and the \bottcher  map.
Given a degree $d$  monic polynomial $f : \C \rightarrow \C$ the Green function 
$$\begin{array}{rccl}
G_f :& \C & \rightarrow & \R_{\geq 0} \\
         & z  & \mapsto     & \lim \frac{\log_+ |f^n(z)|}{d^n}
\end{array}$$
is a well defined continuous function which vanishes in $\filled (f)$ and is harmonic in $\C \setminus  \filled(f)$.
The \bottcher  map $\bot_f : \basinstar(f) \rightarrow \CDC$ is a conformal isomorphism from the basin of infinity $\basinstar$
under the gradient flow $\nabla G_f$ into $\CDC$ which conjugates $f$ with $z \mapsto z^d$ 
(i.e. $\bot_f ( f (z)) = \bot_f (z)^d$ for all $z \in \basinstar(f)$) and is asymptotic to the identity at infinity
(i.e. $\bot_f (z) = z + o(z)$ as $|z| \rightarrow \infty$).

For general background on polynomial dynamics see Section 18 in~\cite{milnor-99}. We now specialize on the set $V$ 
cubic polynomials where
$$V= \{ \fav \,/\,\, -a \notin \filled(\fav), G_\fav (+a) < G_\fav (-a) \}.$$

Following Branner and Hubbard~\cite{branner-88,branner-92} we now summarize the basic combinatorial structure 
of the dynamical plane of polynomials in $V$.
Consider  $f=\fav \in V$ then $$\db^{f}_0 = \{ z \, /\,\, G_{f} (z) < 3 \cdot G_f (-a) \}$$
is a topological open disk which we call
the {\bf dynamical disk of level $0$ of $f$}.  The set $ \{ z \, /\,\, G_{f} (z) < 3^{-n+1} \cdot G_f (-a) \}$ is a finite disjoint 
union of open topological disks that we call {\bf dynamical disks of level $n$}. 
Equivalently, a dynamical disk of level $n$
is a connected component of $f^{-n} (\db^{f}_0)$.   A point $z$ is a {\bf level $n$ point} if $G_{f} (z) < 3^{-n+1} \cdot G_f (-a)$
and the disk of level $n$ containing $z$ is denoted $\db^f_n (z)$. The {\bf level $0$ annulus} is 
$$\da^f_0 = \{ z \, /\,\, G_f (-a) < G_{f} (z) < 3 \cdot G_f (-a) \}$$ and has modulus 
$\mdl A = \pi^{-1} \cdot G_f (-a)$ since it is conformally 
isomorphic to $\{ z \,/ \, \, \mathe^{G_f(-a)} < |z| < \mathe^{3 \cdot G_f(-a)} \}$.  
The {\bf level $n$ annulus} around a level $n$ point $z$
is $$A^f_n (z) = \db^f_n (z) \setminus \{ w \,/\,\, G_{f} (w) \leq 3^{-n} \cdot G_f (-a) \}.$$
A level $n$ annulus $A^f_n$ is said to be {\bf critical}  if $A^f_n=A^f_n (+a)$. 

Here we have that if $z$ is a level $n+1$ point, then $f(\db^f_{n+1}(z)) = \db^f_n (f(z))$ and
$f: \db^f_{n+1}(z) \rightarrow \db^f_n (f(z))$ is a branched covering of degree $2$ if $+a \in \db^f_{n+1}(z) $ and of
degree $1$ otherwise. Now if $+a$ is a level $n+1$ point, then
$f: \da^f_{n+1}(z) \rightarrow \da^f_n (f(z))$ is a covering map of degree $2$ when
$\da^f_{n+1}(z)$ is critical and of  degree $1$ otherwise. 

\medskip
The {\bf tableaux and marked grids for $f$} are defined similarly than in Definition~\ref{tableaux-d}.
Marked grids satisfy the rules stated in Proposition~\ref{rules-p}. 

\smallskip
According to~\cite{branner-88} I.1 the \bottcher map $\bot_\fav (z)$ depends holomorphically in
$(a, v , z)$ for $z$ near infinity. 

\begin{lemma}
  \label{locally-constant-l}
  Suppose that $f_{\ha,\hv} \in V$ and $+\ha$ is a level $n+1$ point under $f_{\ha,\hv}$. 
Then there exists a neighborhood $U$ of $f_{\ha,\hv}$ such that the critical marked grid of level 
$n+1$ is constant in $U$. 
\end{lemma}

\noindent
{\bf Proof.}
Let $f = f_{\ha,\hv}$.
The condition of the critical point $+a$ being a level $n+1$ point is clearly open.
Fix $\ell \geq 0$ and $k \geq 0$ such that $\ell + k \leq n+1$. Let $\rho >0$ be such that 
$3^{-\ell+1} G_f (-\ha) > \rho > 3^{-\ell}G_f (-\ha)$.
The portion of the equipotential $G_f = \rho$ contained in $\db^f_\ell (f^k(\ha))$ is a smooth Jordan curve which varies
smoothly in a neighborhood of $f$. Hence, for $\fav$ in a sufficiently small neighborhood $U$ of $f$, 
the annulus $A^\fav_\ell (f^k_{a,v} (a))$ is critical if and only if 
$+a$ is enclosed by this Jordan curve. Therefore, $A^f_\ell (f^k(\ha))$ is critical if and only if $A^\fav_\ell (f^k_{a,v})$
for all $\fav \in U$ and the lemma follows.
\hfill $\Box$

\subsection{Periodic curves, ends and Puiseux series dynamics}
The key to pass from the parameter space of $\psinu (\z) \in \S [\z]$ to 
the parameter space of cubic polynomials $\hcubicc$ are the Puiseux series of the ends of periodic curves.
The {\bf periodic curve}  $\curve_p$ (of period $p$)
is  the set formed by all $(a,v) \in \hcubicc$ such that $+a$ is periodic of period exactly $p$ under
$\fav$. 

It follows that $\curve_p$ is an algebraic curve in $\C^2 \equiv \hcubicc$. 
For example, $\curve_1 = \{ (a,v) \, /\,\, v-a =0 \}$ and $\curve_2$ is the zero set of
$$(v-a)(v+2a) +1 = \frac{(v-a)^2 (v+2a) + v-a}{v-a}.$$
In general $\curve_p$ is the algebraic curve determined by the polynomial
\begin{equation*}
F_p (a,v) = \frac{f^n_{(a,v)} (+a) - a}{
  \prod_{d | n \;  d \neq n}  f^n_{(a,v)}(+a) -a}  \in \C [a,v].
\end{equation*}
According to Milnor~\cite{milnor-91} 
(compare with Rees~\cite{rees-03}) the curve $\curve_p$ is a smoothly embedded (possibly disconnected) Riemann 
surface in $\C^2$. It is an open problem to determine whether $\curve_p$ is connected (i.e. irreducible) for all $p \geq 1$.
For more about periodic curves we refer the reader to~\cite{milnor-91,rees-03}. 

\medskip
We compactify $\hcubicc \equiv \C^2$ by adding a line $\cL_\infty$ at infinity to obtain $\CP2 \equiv \cubicc \cup \cL_\infty$
and  denote the closure of $\curve_p$ in $\CP2$ by $\overline{\curve_p}$. Since the highest order term of the polynomial 
that defines  $\curve_p$ is of the form $(v-a)^j (v+2a)^k$ for some positive integers $j,k$, we have that 
$\overline{\curve_p} \cap \cL_\infty = \{ [1:1:0], [1:-2:0] \}$. 

\medskip
Our main interest is in $\curve_p \cap \hescapepc$. A connected component $\cF$ of $\curve_p \cap \hescapepc$ 
is called an {\bf end} of $\curve_p$.

\begin{lemma}
  \label{marked-constant-c}
  Let $\cF$ be an end of a periodic curve. Then the marked grid $\bM^\cF$ of $+a$ under $\fav$
is independent of $\fav \in \cF$.
\end{lemma}

\noindent
{\bf Proof.}
Since $+a \in K (\fav)$ for all $\fav \in \cF$, by Lemma~\ref{locally-constant-l}, the subset of $\cF$
where the $(\ell, k)$ position of the critical marked grid is marked (resp.  unmarked) is open, and therefore
closed in $\cF$. \hfill $\Box$

\begin{lemma}
\label{bot-coor-l}
  Let $\cF$ be an end of $\curve_p$. Consider the map
$$\begin{array}{rlll}
\bot_\cF: &\cF & \longrightarrow & \CDC \\
  & (a,v) &\mapsto & \bot_{\fav} (2a).
\end{array}$$
Then $\bot_\cF$ is a holomorphic covering map of finite degree.
\end{lemma}
We say the the degree of $\bot_\cF$ 
is the {\bf multiplicity of the end $\cF$.}

\smallskip
Note that $2a$ is the  ``cocritical point'' of $-a$. That is, $\fav(2a) = \fav (-a) = 4a^3 + v$.
Therefore $G_\fav (2a) = G_\fav(-a)$. It follows that $2a \in \basinstar(\fav)$ for all $\fav \in \curve_p \cap \hescapepc$.

\medskip
\noindent
{\bf Proof.}
From Branner and Hubbard's wringing construction~\cite{branner-88} it follows that $\bot_\cF$ is a local homeomorphism.
Note  that $\bot^{-1}_\cF (w_0)$ is closed and bounded. To prove that this set is finite it suffices to show that every point is 
isolated. In fact, if $(a_0,v_0) \in \bot^{-1}_\cF (w_0)$, then $\bot_\cF$ extends holomorphically to a map $\bot$ on a 
neighborhood $U_0$ of $(a_0,v_0)$ in $\C^2$ and therefore $\bot^{-1} (w_0) \cap \curve_p$ is discrete in $U_0$.
\hfill $\Box$

\begin{corollary}
  If $\cF$ is an end of $\curve_p$, then  $\overline{\cF} \cap \cL_\infty$ consists of exactly one point $x$. Moreover, 
the germ of the analytic set $\cF \cup \{ x \}$ at $x$  is irreducible. 
\end{corollary}

\noindent
{\bf Proof.} 
Since $\bot_\cF$ is a finite degree covering of a punctured disk, it follows that $\cF$ is conformally isomorphic to a 
punctured disk and $\overline{\cF} \cap \cL_\infty$ consists of exactly one point $x$. Moreover, the germ of
$\cF \cup \{ x \}$ at $x$ is irreducible since a fundamental system of (punctured) neighborhoods of $x$ in $\cF$ 
is given by the connected sets $\bot^{-1}_\cF (\{ |w| >k\})$ where $k \in \N$. \hfill $\Box$

\medskip
To study $\curve_p$ near $\cL_\infty$ in $\CP2 = \{ [a:v:w] \,/\,\, (a,v,w) \in \C^3 \setminus \{ 0\} \}$ we consider the
open set $U_1 = \{ [a:v:w] \,/\,\, a \neq 0 \}$ which we identify  via $[1:\bv:\ba] \mapsto (\ba,\bv)$
 with a copy $\C^2_{(\ba,\bv)}$ of $\C^2$.
Note that $\cL_\infty \cap \C^2_{(\ba,\bv)}$ is
the $\bv$-axes (i.e.  $\ba=0$). Moreover, every end $\cF$ is contained in $\C^2_{(\ba,\bv)}$.
The series
$$\nu = 
\sum_{j \geq 0} a_j t^{j/m} \in \pui$$
is called a {\bf Puiseux series of an end $\cF$} if there exists $ \epsilon >0$
such that $$\sum_{j \geq 0} a_j s^j$$
converges in $|s| < \epsilon^{1/m}$ and
$$\begin{array}{clc}
  \{ 0 < |s| < \epsilon^{1/m} \}  & \rightarrow & \cF \cap 
  \{ (\ba, \bv) \in \C^2_{(\ba,\bv)} / 0 < |\ba| < \epsilon \} \\
  s &\mapsto & (s^m , \sum_{j \geq 0} a_j s^j)
\end{array}$$
is a conformal isomorphism (e.g., see~\cite{brieskorn-86,casas-00,fischer-94}). 

Puiseux series always exist and they are unique modulo an automorphism of
$\lautm$ over $\laut$. That is, an automorphism
which  sends $t^{1/m}$ to $e^{2 \pi i k/m} t^{1/m}$ for some $0 \leq k < m$.

\medskip
A sufficient condition for the existence of a Puiseux series which ``converges'' in an $\epsilon$-disk follows immediately
from standard covering space theory:

\begin{lemma}
  \label{epsilon-pui-l}
Let $\cF$ be an end of a curve $\curve_p$.
  Suppose $\epsilon >0$ is such that the projection:
 $$\begin{array}{clc}
\cF \cap 
  \{ (\ba, \bv) \in \C^2_{(\ba,\bv)} / 0 < |\ba| < \epsilon \} & \rightarrow & \{ 0 < |\ba| < \epsilon \}  \\
  (\ba,\bv) &\mapsto & \ba
\end{array}$$
is a covering of degree $m$.
Then there exists a holomorphic map $g(s) = \sum_{j \geq j_0} a_j s^j$ defined on
$|s| < \epsilon^{1/m}$ so that $s \mapsto (s^m, g(s))$ is a conformal isomorphism from
$\{ 0< |s| < \epsilon^{1/m} \}$ onto $\cF \cap   \{ (\ba, \bv) \in \C^2_{(\ba,\bv)} / 0 < |\ba| < \epsilon \}$.
Moreover, if $\tilde{g}$ is another holomorphic function such that the above holds, then
$g(s) = g(\mathe^{2 \pi i k/m} s)$ for some integer $k$. 
\end{lemma}

We refer the reader to~\cite{brieskorn-86} or~\cite{fischer-94}  
for an elementary exposition of the previous results.
For future reference we state several equivalent conditions in the following lemma.
The proof is straightforward and we omit it. 

\begin{lemma}
\label{cover-l}
  Consider a real number  $\epsilon >0$, an end $\cF$ of $\curve_p$ and a series
$\nu  = \sum_{j \geq 0} a_j t^{j/m} \in \pui$. Then the following are equivalent:

(i) $$\begin{array}{clc}
  \{ 0 < |s| < \epsilon^{1/m} \}  & \rightarrow & \cF \cap 
  \{ (\ba, \bv) \in \C^2_{(\ba,\bv)} / 0 < |\ba| < \epsilon \} \\
  s &\mapsto & (s^m , \sum_{j \geq 0} a_j s^j)
\end{array}$$
is a conformal isomorphism.

(ii) $$\begin{array}{rclc}
 \bav_\nu: & \hepsilon  & \rightarrow & \cF \cap 
  \{ (\ba, \bv) \in \C^2_{(\ba,\bv)} / 0 < |\ba| < \epsilon \} \\
 & T &\mapsto & (\mathe^{2 \pi i T} , \sum_{j \geq 0} a_j \mathe^{2 \pi i T j/m})
\end{array}$$
is a well defined conformal map which is the universal covering of its image.

(iii) $$\begin{array}{rclc}
 \av_\nu:  & \hepsilon  & \rightarrow & \cF \cap 
  \{ (a, v) \in \cubicc \, / \,\, |a| >  \epsilon^{-1} \} \\
  & T &\mapsto & (\mathe^{-2 \pi i T} , \mathe^{-2 \pi i T} \sum_{j \geq j_0} a_j \mathe^{2 \pi i T j/m})
\end{array}$$
is a well defined conformal map which is the universal covering of its image.
\end{lemma}

To characterize Puiseux series from an algebraic and dynamical viewpoint we observe that
the polynomial $F_p (a,v)$ which defines $\curve_p$ 
 becomes:
$$\overline{F_p} (\ba,\bv) = \ba^{\deg F_p} F_p \left(\frac{1}{\ba}, \frac{\bv}{\ba}\right) \in \C [\ba,\bv]$$
in $(\ba,\bv)$ coordinates.
When $\overline{F_p} (t, \bv)$ is regarded as an element of $\laut [\bv]$, then 
$\nu \in \pui$ is a Puiseux series of an end of $\curve_p$
if and only if $\nu$ is a root of $\overline{F_p} (t, \bv) \in \laut [\bv]$. 

\begin{lemma}
\label{dyn-char-l}
  For $\nu \in \sol$ let 
  $$\psi_\nu (\z) = t^{-2} (\z -1)^2 (\z +2) + \nu.$$
A series $\nu \in \pui$ is the Puiseux series of an end of $\curve_p$ if and only if $\omep =+1$ is  periodic of period
exactly $p$ under $\psi_\nu : \sol \rightarrow \sol$. 
\end{lemma}

\noindent
{\bf Proof.}
From the algebraic characterization of Puiseux series described above,
 $\nu \in \pui$ is the Puiseux series of an end of $\curve_p$
 if and only if $\nu \in \pui$ is such that $\overline{F_p}(t,\nu) =0 \in \pui$.
Since $\pui$ is algebraically closed, this is equivalent to $\nu$ being an element of 
$\sol$ is such that $\overline{F_p}(t,\nu) =0 \in \sol$.
The latter equality is equivalent to:
\begin{eqnarray*}
  f^p_{(1/t, \nu/t)} (1/t) - 1/t & = & 0 \quad  \mbox{ and}, \\
  f^k_{(1/t,\nu/t)}(1/t)-1/t & \neq & 0 \quad  \mbox{ for } k < p.
\end{eqnarray*}
After changing coordinates to $\z = t z$ the map $f_{(1/t, \nu/t)}(z)$ becomes $\psi_\nu (\z)$ and the lemma follows.
\hfill $\Box$

\begin{remark}{\em
\label{automorphism-r}
Recall that $\sigma: \S \rightarrow \S$ is the unique automorphism such that $\sigma(t^{1/m}) = \mathe^{2 \pi i /m} t^{1/m}$ for all
$m \in \N$. 
For all $\nu$ and $\z$ in $\sol$,
$$\psi_{\sigma(\nu)} (\z) = \sigma \circ \psinu \circ \sigma^{-1}(\z).$$
In particular, 
$\omep \in \filled (\psinu)$ if and only if $\omep = \sigma(\omep) \in \filled (\psi_{\sigma(\nu)})$.
In this case, $\bM^\nu (\omep) = \bM^{\sigma (\nu)} (\omep)$.}
\end{remark}

\subsection{Convergent series}
\label{convergent-ss}
Here we remark a trivial but useful fact related to ``convergent series'' in $\S$.
We say that $$\z = \sum_{\lambda \geq \lambda_0} a_\lambda t^\lambda \in \S$$
{\bf converges in $\hepsilon$} to $T \mapsto \z (\mathe^{2 \pi i T})$ if the series
$$\sum_{\lambda \geq \lambda_0} a_\lambda \mathe^{2 \pi i T \lambda}$$
converges for all $T \in \hepsilon$.
For convergent series, evaluation and iteration commute:
\begin{lemma}
  \label{evaluation-l}
  Assume that $\nu, \z \in \S$ are convergent in $\hepsilon$. Let $a(T) = \mathe^{-2 \pi iT}$,
$v(T) = \mathe^{-2 \pi iT}\nu (\mathe^{2 \pi i T})$ and $z(T) = \zeta (\mathe^{2 \pi iT})$.
Then  $\psinu^k(\z)$ converges in $\hepsilon$ and 
$$\mathe^{2 \pi iT} f^k_{a(T),v(T)} (\mathe^{-2 \pi iT} z(T)) = \left( \psinu^k(\z) \right) (\mathe^{2 \pi iT})$$
for all $k \geq 1$.
\end{lemma}
The proof of the lemma is straightforward and we omit it.

\subsection{Uniform radius of convergence}
Our aim now is to show that the ``radii of convergence''  of any  Puiseux series of any end of any curve $\curve_p$ 
is uniformly bounded below. More precisely:

\begin{proposition}
  \label{uniform-p}
  There exist $C>0$ and $\epsilon > 0$ such that
for any  Puiseux series 
$$\nu = \sum_{\lambda \geq \lambda_0} a_\lambda t^\lambda$$
of any end of any periodic curve $\curve_p$ the following hold:

(i) $\bv_\nu (T) = \sum_{\lambda \geq 0} a_\lambda \mathe^{2 \pi i T \lambda}$
converges in $\hepsilon$.

(ii) $|\bv_\nu (T)| \leq C$ for all $T \in \hepsilon$.
\end{proposition}

Note that $\epsilon$ and $C$ in the above proposition are independent of $p \geq 1$.

The rest of this subsection is devoted to the proof of this proposition which relies on a series of lemmas.
The main idea is to use the \bottcher  function to change coordinates near $\cL_\infty$.

We will need the following elementary estimates:

\begin{lemma}
  \label{estimates0-l}
  Suppose that $(a,v) \in \C^2$ is such that $|a| >2$ and
either $| v +2a| <|a|$ or $|v -a| <|a|$. Then:

(i) $|f^n (-a)| \geq |a|^{3^n} \left( \sqrt{2} \right)^{-{3^{n-1} -1}}$.

(ii) $\log |a| - \frac{1}{6} \log 2 \leq G_{\fav} (-a) \leq \log |a| + \frac{1}{2} \log \frac{3}{2}.$

(iii) $G_{\fav} (+a) \leq \frac{1}{3} \log |a| + \frac{1}{18} \log 25 \cdot 32$.
\end{lemma}

\noindent
{\bf Proof.}
Let $f = \fav$.
First note that 
\begin{equation}
  \label{est1-e}
  |z| > 3 |a| \implies |f(z)| > 3 |a|.
\end{equation}
In fact, using that $|v| < 3 |a|$ we obtain $|f(z)| \geq |z|^3 ( 1 - |\frac{a}{z}|^2 - 2|\frac{a}{z}| - |\frac{v}{z^3}|) \geq 12 |z| \geq 3 |a|$.
Similarly we have:
\begin{equation}
  \label{est3-e}
 |z| > 3 |a| \implies \frac{1}{2} |z|^3 \leq |f(z)| \leq \frac{3}{2} |z|^3.
\end{equation}
Also,
\begin{equation}
  \label{est2-e}
  |z| < 3 |a| \implies |f(z)| < 40 |a|^3,
\end{equation}
since  $|f(z)| \leq  |a|^3 ( 2 + \frac{|v|}{|a|^3} + 3 | \frac{z}{a}| +  | \frac{z}{a}|^3) \leq 40 |a|^3$.
Inductively applying equation (\ref{est3-e}) we have:
\begin{equation}
  \label{est5-e}
3 |a| < |z| < 40 |a|^3 \implies |f^n (z)| < 40^{3^n} |a|^{3^{n+1}}\left(\frac{3}{2}\right)^{\frac{3^{n-1}-1}{2}}
\end{equation}
 
\smallskip
Now we prove (i). Since  $|f(-a)| = |4a^3 +v | \geq |a|^3 > 3 |-a|$, by induction, from equation (\ref{est3-e})  we conclude
that for all $n$:
$$|f^n (-a)| \geq |a|^{3^n} \left(\frac{1}{2}\right)^{\frac{3^{n-1} -1}{2}}.$$
From here, taking logarithms, dividing by $3^n$ and passing to the limit
we obtain the lower bound of (iii). 
For the upper bound of (iii) note that $3|-a|< |f(-a)| \leq 7 |a|^3 < 40 |a|^3$. Hence, applying equation (\ref{est5-e})
the desired upper bound follows.

For (iii), either $G_f(a) =0$ or there exists $n_0 \geq 1$ such that $|f^{n}(a)| \leq 3|a|$ for $n \leq n_0$
and $3|a| < |f^{n_0+1}(a)| < 40 |a|^3$. In the latter case, let $w = f^{n_0}(a)$ and from equation (\ref{est5-e}) conclude
that $G_f(w) \leq 3 \log |a| + \frac{1}{2} \log{25 \cdot 32}$. Since $G_f(w) = 3^{n_0+1} G_f (a) \geq 9 G_f (a)$ the upper
bound of (iii) follows.
\hfill $\Box$

\medskip
The complex plane $\{ (a,a) \,\, /\, a \in \C \}$ is the period $1$ curve $\curve_1$ and 
$$\{(a,-2a) \,\,/\, a \in \C\}$$ 
is the set of parameters for which $+a$ is a prefixed critical point. The asymptotic behavior of the \bottcher map along
this complex planes is easily computed (compare with~\cite{milnor-91}).
 
\begin{lemma}
  \label{bot-l}
  For $q =-2$ or $1$ the following holds:
$$\lim_{|a| \rightarrow \infty} \frac{\bot_{f_{(a, qa)}}(2a)}{a} = 2^{2/3}.$$
\end{lemma}

\noindent
{\bf Proof.}
Let $f_a = f_{(a,qa)}$ and $\bot_a = \bot_{f_a}$. 
Then
$$ \bot_a (2a) = \lim_{n \rightarrow \infty} 2a \left( \frac{f_a(2a)}{(2a)^3} \right)^{1/3} \cdot  
\left( \frac{f^2_a(2a)}{(f_a (2a))^3} \right)^{1/3^2} \cdots \left( \frac{f^n_a(2a)}{(f^{n-1}_a (2a))^{3^{n-1}}} \right)^{1/3^{n-1}}$$
which, for $|a|$ sufficiently large, is a uniform limit. Since $f(-a) = f(2a)$ we have that
\begin{eqnarray*}
\log \left| \frac{f^i (2a)  }{f^{i-1}(2a)^{3^{i-1}}} \right|^{1/3^i} & = & \frac{1}{3^i} \log| 1 - 3 \left( \frac{a}{f^i (2a)} \right)^2 + 2 
 \left( \frac{a}{f^i (2a)} \right)^3 + \frac{q a}{(f^i(2a))^3} \\
& \leq & \frac{K}{3^i}
\end{eqnarray*}
for some $K \geq 0$.
The lemma follows after  switching the order of  the limits.
\hfill $\Box$

\begin{lemma}
  \label{change-l}
  There exists a neighborhood $U$ of $\{ [1:2:0], [1:1:0] \}$ in $\CP2 \equiv \cubicc \cup \cL_\infty$ such that
$$  \begin{array}{rlll}
      \Phi : & U & \rightarrow & \C^2 \\
               &  x & \mapsto      & \left\{ \begin{array}{ll}
                                                    \left(\frac{1}{\bot_{\fav}(2a)}, \frac{v}{a} \right) & \mbox{if } x = \fav \in \cubicc \\
                                                    (0,\frac{v}{a})                                                      & \mbox{if } x = [a:v:0]
                                                           \end{array} \right. \\
\end{array}$$
is a well defined holomorphic function which is locally biholomorphic at $x =  [1:2:0]$ and $x= [1:1:0]$.
\end{lemma}

\noindent
{\bf Proof.}
Let $U_0 = \{ [1:\bv:\ba] \,/\,\, |\ba|< 1/3 \mbox{ and } \min\{|\bv +2|,|\bv-1|\}<1 \}$.
By Lemma~\ref{estimates0-l} (ii) and (iii), 
$\Phi : U_0 \rightarrow \C^2$ is well defined. Moreover,
$\Phi$ is  holomorphic in $U_0 \cap \cubicc = U_0 \setminus \cL_\infty$
and continuous in $U_0$. It follows that $\Phi$ is holomorphic in $U_0$. 
Using $(\ba,\bv)$ coordinates in $U_0$, by the previous lemma, for $q=1$ or $-2$ the derivative $D \Phi_{[1:q:0]}$ has the form:
$$\left[ \begin{array}{cc}
2^{-2/3} & \star  \\
0            &  1    \\
  \end{array} \right]$$
Hence $\Phi$ is locally biholomorphic at $[1:q:0]$.
\hfill $\Box$

\begin{lemma}
  \label{third-l}
  Let $f = \fav$ and suppose that 

(i) $G_f (-a) \geq 3 G_f(v)$,

(ii) $\db^f_1 (v) = \db^f_1 (q a)$ for $q=1$ or $-2$.

Then
$$\left| \frac{v}{a} -q \right| \leq \frac{16}{\mathe^{|q| G_f (-a)} -16}.$$
\end{lemma}

\noindent
{\bf Proof.}
Note that $2 \pi \mdl A^f_1 (v) =  |q| G_f (-a)$ and $A^f_1 (v)$ separates $qa$ and $v$ from $-a$ and $\infty$.
Let $\Gamma(z) = -a (z+2a)^{-1}$ and note that $\Gamma(A^f_1(v))$ separates $-1$ and $0$ from $\infty$
and $\Gamma(v) = (\frac{v}{a} - q)^{-1}$. According to Chapter III in~\cite{ahlfors-66}, for some function $\Psi$:
$$\frac{1}{2 \pi} |q| G_f(-a) \leq \frac{1}{2 \pi} \Psi (|\Gamma(v)|) \leq \frac{1}{2 \pi} \log 16(|\Gamma(v)|+1).$$
From where the desired inequality immediately follows.
\hfill $\Box$

\begin{lemma}
  \label{contained-l}
Let $q=-2$ or $1$. For $\rho >0$ let $E_\rho$ be the set formed by all $(a,v) \in \hcubicc$ such that

(i) $G_{\fav}(-a) \geq 3 G_{\fav}(v)$.

(ii) $\db^{\fav}_1 (v) = \db^{\fav}_1 (qa)$.

(iii) $G_{\fav}(-a) > \rho$.

For any neighborhood $W$ of $[1:q:0]$ in $\CP2 \equiv \cubicc \cup \cL_\infty$, there exists $\hat{\rho}$
such that $E_{\hat{\rho}} \subset W$.
\end{lemma}

\noindent
{\bf Proof.}
Take $\rho_0 >0$ such that (i), (ii) and $G_{\fav}(-a) > \rho_0$ imply that $|v-qa| < |a|$ (see Lemma~\ref{third-l}).
Let $\rho_1 \geq \rho_0$ be such that $G_{\fav} (-a) \leq \rho_1$ for all $(a,v)$ such that $|a| \leq 2$ and $|v|\leq 6$.
It follows that $G_{\fav}(-a) > \rho_1$ implies that $|a|>2$ and $|v-qa| < |a|$. Now let $\delta >0$ be such that
$\{ [1: \bv: \ba] \, / \,\, |\ba| < \delta, |\bv - q| < \delta \} \subset V$. By Lemma~\ref{estimates0-l} (ii) and Lemma~\ref{third-l}, 
there exists $\hat{\rho}$ such that $G_{\fav} (-a) > \hat{\rho}$ implies that $|a|^{-1} = |\ba| < \delta$ and $|v/a -q| = |\bv -q| < \delta$.
\hfill $\Box$

\medskip
\noindent
{\bf Proof of Proposition~\ref{uniform-p}.}
Let $U$ be a neighborhood of $\{[1:1:0],[1:-2:0]\}$ so that $\Phi : U \rightarrow \Phi(U) \subset \C^2$
is biholomorphic. In the range we will use coordinates $(\bb, \bv)$.
Let $\rho >0$ be such that $\{ \fav \in \curve_p \,/\,\, G_{\fav}(-a) \geq \rho \} \subset U$.

Since $\bot: \curve_p \rightarrow \CDC$ given by $\bot (a,v) = \bot_{(a,v)}(2a)$ is a local homeomorphism, for all $p \geq 1$,
the line tangent to $\Phi(\curve_p \cap U)$ at any $(\bb_0, \bv_0)$ is not vertical. Therefore, it has equation
$$\bv -\bv_0 = M(\bb_0,\bv_0)(\bb-\bb_0)$$
for some $M(\bb_0,\bv_0) \in \C$. 

We claim that there exists $C_\rho >0$ such that
$$|M(\bb_0,\bv_0)| < C_\rho$$
for all $(\bb_0,\bv_0) \in \cup_{p\geq 1} \Phi(\curve_p \cap U) \cap \{| \bb_0 | = \mathe^{-\rho}\}$.
For otherwise, for all $n\geq 1$ there would exist $(\bb_n,\bv_n) \in \curve_{p_n}$ with $|M(\bb_n,\bv )| \geq n$.
Without loss of generality we may assume that $\bb_n \rightarrow \bb_\infty$ and $\bv_n \rightarrow \bv_\infty$.
Since the projections of $\Phi(\curve_{p_n} \cap U)$ to the $\bb$--axis are covering maps, there exist a neighborhood $W$ 
of $\bb_\infty$ and holomorphic functions $h_n : W \rightarrow \C$ with graphs contained in $\curve_{p_n}$ such that
$h_n (\bb_n) = \bv_n$. By normality of $\{ h_n \}$ we may suppose that $h_n \rightarrow h_\infty$ locally uniformly in $W$.
Hence $M(\bb_n, \bv_n) = h^\prime_n (\bb_n) \rightarrow h^\prime_\infty (\bb_\infty)$ which contradicts 
$|M(\bb_n,\bv_n )| \geq n$.

\smallskip
We will now apply the Maximum Principle to conclude that $|M(\bb_0,\bv_0)| < C_\rho$ for all 
$(\bb_0,\bv_0) \in \cup_{p\geq 1} \Phi(\curve_p \cap U) \cap \{| \bb_0 | \leq \mathe^{-\rho}\}$.
According to our description of the level $1$ parameter balls the Puiseux series $\nu$ of any end $\cF$
is either $\mathe^{-1}$ close to $1 \in \sol $ or $\mathe^{-2}$ close to $-2 \in \sol$. In particular,
$$\nu = a_0 + a_1 t + \mbox{ higher order terms} \in \lautm$$
for some $m \in \N$.
If follows that the germ of $\Phi(\overline{\cF} \cap U)$ at $\Phi(\overline{\cF} \cap \cL_\infty)$ has Puiseux series
of the form
$$\mu = b_0 + b_1 u  +  \mbox{ higher order terms} \in \C((u^{1/m})).$$
Since the $\bb$--projection $\Pi: \Phi({\cF} \cap U) \rightarrow \C$ is a covering map and 
$\cF \cap U \supset \cF \cap \{ (a,v) \,/\,\, |\bot_{\fav} (2a)| \geq \mathe^\rho \}$, we have that
for all $s$ such that $|s| < \mathe^{\rho/m}$, the series $\mu(s^m)$ is convergent and 
$s \mapsto (s^m, \mu(s^m))$ parametrizes $\Phi(\cF)$. It follows that 
$$|M(s_0^m, \mu(s_0^m))| = \left|\left(\frac{d\mu(s^m)}{ds} (s_0)\right)/(ms^{m-1}_0)\right| = O(1)$$
as $s_0 \rightarrow 0$.
Hence, by the Maximum Principle, the holomorphic function $M(s^m,\mu(s^m))$ is bounded by $C_\rho$ in
$|s| \leq \mathe^{-\rho/m}$. Therefore, $|m(\bb_0,\bv_0)| \leq C_\rho$ for $|\bb_0| \leq \mathe^{-\rho}$. 

\smallskip
The proposition follows if we show that for $\epsilon >0$ sufficiently small the $\ba$--projection
of $\curve_p \cap \{[1:\bv:\ba] \, /\,\, |\ba| < \epsilon \}$ is a covering map. Equivalently,
the tangent line to points in $\curve_p \cap \{[1:\bv:\ba] \, /\,\, |\ba| < \epsilon \}$ is not vertical.
In $U$ we work with $(\ba,\bv)$ coordinates and observe that $D\Phi_x (0,1) = (0,1)$ for
$x = (0,1)$ or $(0,-2)$. Hence, shrinking $U$ if necessary, we have that 
$|d/c| > C_\rho$ where $(c,d) = D\Phi_{(\ba,\bv)} (0,1)$ for all $(\ba,\bv) \in U$.
It follows that a tangent line to $\curve_p \cap U$ is never vertical. 
To conclude the proof of the proposition take $\epsilon >0$ small so that 
for all $(a,v) \in \curve_p$ such that $|\ba| < \epsilon$ we have that $G_{\fav}(-a) > \hat{\rho}$ 
where $\hat{\rho}$ is such that $(a,v) \in \curve_p$ and $G_{\fav}(-a) > \hat{\rho}$ implies that
$(a,v) \in U$. 
\hfill $\Box$

\subsection{Marked grid correspondence}

For Puiseux series of ends of periodic curves the marked grids of the dynamics over $\S$ and $\C$ coincide:

\begin{proposition}
  \label{correspondence-p}
Consider an end $\cF$ and let $\bM^\cF$ be the critical marked grid of all the polynomials in $\cF$.
If $\nu \in \S$ is a Puiseux series of $\cF$, then $\bM^\nu = \bM^\cF$ where $\bM^\nu$ is the critical marked
grid of $\psinu : \S \rightarrow \S$.
\end{proposition}

The rest of this subsection is devoted to the proof of the above proposition. 
For this we let $\epsilon >0$ be such that $\nu$ converges in $\hepsilon$.

To simplify notation we let $f_T = f_{a(T), v(T)}$ where $a  (T) =\mathe^{-2 \pi iT}$
and $ v  (T) = \mathe^{-2 \pi iT} \nu(\mathe^{2 \pi iT}) $.
The corresponding Green function will be denoted by $G_T$ 
and the level $n$ dynamical disks by $\db^T_n(z)$. Similarly the level $n$ annuli are written as $\da^T_n (z)$. 

Now for $k \geq 0$ let $\nu_k = \psi^k_\nu(\omep)$ and $f^k_T (a (T)) = v_k(T)$. 
From Lemma~\ref{evaluation-l} we have that
$v_k(T)  = \mathe^{-2 \pi i T} \cdot \nu_k ( \mathe^{2 \pi i T})$.

\begin{lemma}
  \label{zero-l}
The following hold:

(i) $$\frac{\pi \mdl A^T_0}{\log | a(T)|} = \frac{G_T(-a(T))}{\log | a(T)| } \rightarrow 1 \mbox{  as  } \im T \rightarrow +\infty.$$

(ii) Let $q=1$ or $-2$ and 
assume that $z_i : \hepsilon \rightarrow \C$ are functions such that, for all $T \in \hepsilon$,
$z_1(T)$ is a level $2$ point under $f_T$, $z_1(T) \in \db^T_1 (q a(T))$ and there exists an annulus 
$A^T$ which separates $z_1(T)$ and $z_2(T)$ from $-a(T)$ and $\infty$. 
If there exists $S >0$ such that 
$$\frac{\mdl A^T}{\mdl A^T_0} \geq S, $$ then
$$\left| \frac{z_1(T)}{a (T)} - \frac{z_2(T)}{a(T)} \right| = O( | a (T) |^{-2S})$$
as $\im T \rightarrow +\infty$.
\end{lemma}

\noindent
{\bf Proof.}
Since (i) is a direct consequence of Lemma~\ref{estimates0-l} we proceed to prove (ii).
For this consider the M\"oebius transformation $\Gamma_T (z) = (a(T) + z_1(T))(z -z_1(T))^{-1}$.
The annulus $\Gamma_T(A^T)$ separates $-1$ and $0$ from $\infty$ and 
$$\left( 1 + \frac{z_1(T)}{a(T)} \right) \left( \frac{z_1(T)}{a (T)} - \frac{z_2(T)}{a(T)} \right)^{-1}.$$
From Chapter III in~\cite{ahlfors-66} we have that
$$S \mdl A^T_0 \leq \mdl A^T \leq \frac{1}{2 \pi} \log 16(\Gamma_T(z_1(T)) +1).$$
From (i) we conclude that for $\im T$ sufficiently large:
\begin{equation}
\label{extremal-e}
\left| \frac{z_1(T)}{a (T)} - \frac{z_2(T)}{a(T)} \right| \leq \frac{16}{a(T)^{2S} - 16} \left( 1 + \frac{z_1(T)}{a(T)} \right)
\end{equation}
Now replace in the previous equation $z_1(T)$ by $ q a(T)$, 
$z_2(T)$ by $z_1(T)$, $A^T$ by the annulus $A^T_1(qa(T))$ of modulus $|q+1|^{-1} \mdl A^T_0$ and
conclude that, as $\im T \rightarrow \infty$,  
$$1 + \frac{z_1(T)}{a(T)} \rightarrow q+1.$$
Combining this with equation (\ref{extremal-e}) part (ii) of the lemma follows.
\hfill $\Box$

\medskip
For any $j, k \geq 0$, the level $n$ disks $\db^T_n(v_j(T))$ and  $\db^T_n(v_k(T))$
are either equal for all $T \in \hepsilon$ or distinct for all $T \in \hepsilon$.
If we denote by $M^\cF_{\ell, k}$ the entries of $\bM^\cF$, then
$$\mdl A^T_\ell (v_k (T)) = 2^{-S_\ell} \mdl A^T$$
where $S_\ell = \sum^{\ell-1}_{k =0} M^\cF_{\ell-k,j+k}$.
By the Gr\"otzsch  inequality  (see~\cite{ahlfors-66}) it follows that 
$$\mdl \db^T_1 (v_k(T)) \setminus \overline{\db^T_n (v_k(T))} \geq 
\sum^{n}_{\ell =1}  \mdl A^T_\ell (v_k(T)).$$
Taking $A^T =  \db^T_1 (v_k(T)) \setminus \overline{\db^T_n (v_k(T))}$ in the previous lemma
we immediately obtain the following:

\begin{corollary}
  \label{order-c}
  If for some $j,k \geq 0$ and $n \geq 1$ we have that  $\db^T_n (v_j(T)) = \db^T_n (v_k(T))$, then
$$\order(\nu_j - \nu_k) \geq \frac{2}{\mdl A^T_0}  \sum^n_{\ell =1} \mdl A^T_\ell (v_j(T)).$$
\end{corollary}

Proposition~\ref{correspondence-p} is an immediate consequence of the following:

\begin{lemma}
  \label{n-l}
For all $n \geq 1$ and $j, k \geq 0$ we have that
$\db^\nu_n (\nu_j) = \db^\nu_n (\nu_k)$ if and only if $\db^T_n (v_j(T)) = \db^T_n(v_k(T))$.
\end{lemma}

\noindent
{\bf Proof.} We proceed by induction on $n$.

For $n=1$, if $D^T_1 (v_j(T)) = D^T_1 (a(T))$, then $\order (\nu_j - \omep) \geq 1$. Hence $\nu_j \in \dbnu_1 (\omep)$.
If $D^T_1 (v_j(T)) = D^T_1 (-2 a(T))$, then $\order (\nu_j - (-2)) \geq 2$. Therefore, $\nu_j \in \dbnu_1 (-2)$. 
It follows that $D^T_1 (v_j(T)) =D^T_1 (v_k(T))$ if and only if $\dbnu_1 (\nu_j) = \dbnu_1(\nu_k)$.

Now suppose that the lemma is true for $n$. Note that this implies that
$$\mdl \da^\nu_\ell (\nu_j) = \frac{2 \pi}{G_T(-a(T))} \mdl A^T_\ell (v_j(T))$$
for all $\ell \leq n$. Therefore, if $\db^T_{n+1}(v_j (T)) =\db^T_{n+1}(v_k (T))$ then
\begin{eqnarray*}
\order (\nu_j - \nu_k ) & \geq &\sum^{n+1}_{\ell=1} \mdl \da^T_\ell (v_j(T))  > \sum^{n}_{\ell=1} \mdl \da^T_\ell (v_j(T)) \\
                                   & =       & \sum^{n}_{\ell=1} \mdl \da^\nu_\ell (\nu_j) = -\log r^\nu_n (\nu_j). 
\end{eqnarray*}
Hence, $\lval \nu_j - \nu_k \rval < r^\nu_n (\nu_j)$ and $\nu_j, \nu_k$ belong to the same element of the affine
partition associated to $\db^\nu_n (\nu_j)$. By Lemma~\ref{modulus-action-l} (i), $\nu_k \in \db^\nu_{n+1}(\nu_j)$.

\smallskip
To finish the proof it is sufficient to show that if $\db^T_{n}(v_j (T)) = \db^T_{n}(v_k (T))$ and
$\db^T_{n+1}(v_j (T)) \neq  \db^T_{n+1}(v_k (T))  $, then $\db^\nu_{n+1}(\nu_j) \neq \db^\nu_{n+1}(\nu_k) $.
There are two cases.

\smallskip
\noindent
{\it Case 1. $\db^T_{n}(v_{j+1} (T)) \neq \db^T_{n}(v_{k+1} (T))$:} By the inductive hypothesis, 
$\db^\nu_{n}(\nu_{j+1}) \neq \db^\nu_{n}(\nu_{k+1})$. Thus, $\db^\nu_{n+1}(\nu_j) \neq \db^\nu_{n+1}(\nu_k) $.

\smallskip
\noindent
{\it Case 2. $\db^T_{n}(v_{j+1} (T)) = \db^T_{n}(v_{k+1} (T))$:} In this case $j,k >0$ and 
$f_T : \db^T_n (v_j(T)) \rightarrow \db^T_{n-1}(v_{j+1}(T))$ has degree $2$. Thus, there exist
$v^\prime_j (T)$ and $v^\prime_k (T)$ in $ \db^T_n (v_j(T)) $ distinct from $v_j (T)$ and $v_k (T)$,
respectively, such that $f_T (v^\prime_j (T)) = v_{j+1}(T)$ and $f_T (v^\prime_k(T)) = v_{k+1}(T)$. 
Similarly, $\psinu: \dbnu_n (\nu_j) \rightarrow \dbnu_{n-1}(\nu_{j+1})$ has degree $2$ and there
exist $\nup_j$ and $\nup_k$ in $\dbnu_n(\nu_j)$ distinct form $\nu_j$ and $\nu_k$, respectively,
such that $\psinu(\nup_j)= \nu_{j+1}$ and $\psinup(\nup_k) = \nup_{k+1}$. It follows that
$\nup_j (\mathe^{2 \pi i T}) = \mathe^{2 \pi i T} v^\prime_j (T)$ and $\nup_k (\mathe^{2 \pi i T}) = \mathe^{2 \pi i T} v^\prime_k (T)$.
We claim that $$\order(\nup_j - \nu_k) > \order (\nu_j - \nu_k).$$
In fact, let $$\Gamma_T (z) = \frac{z- v_k (T)}{v_j (T) - v_k (T)}.$$ 
Since  $v^\prime_j (T) \in \db^T_{n+1}(v_k(T))$ the annulus $\da^T_{n+1}(v_k(T))$ separates $v^\prime_j$ and $v_k$
from $v_j$ and $\infty$. Therefore, $\Gamma(v^\prime_j) \rightarrow 0$ as $\im T \rightarrow \infty$ because
$\mdl \da^T_{n+1}(v_k(T)) \rightarrow \infty$. Thus $\order(\nup_j - \nu_k) > \order (\nu_j - \nu_k)$.
From here we conclude that $\db^\nu_{n+1}(\nu_j) \neq \db^\nu_{n+1}(\nu_k) $. For otherwise $\dbnu_{n+1}(\nup_j) 
\neq \dbnu_{n+1}(\nu_k)$ and
$$\order (\nu_j - \nu_k) > \sum^n_{\ell =1} \danu_\ell(\nu_j) = \order (\nup_j - \nu_k).$$
\hfill $\Box$

\subsection{Almost periodic functions}
\label{almost-ss}
Let us briefly summarize some facts about analytic almost periodic functions (in the sense of Bohr).
For a more detailed discussion we refer the reader to~\cite{besicovitch-54}.
A function $f : \R \rightarrow \C$ is {\bf almost periodic} if the family $\{f (\cdot + h) \,\, /\, h \in \R \}$ is a normal family.
That is, every sequence in this family has a subsequence which converges uniformly in $\R$ (i.e., in the $\sup$-norm).
An analytic function $f: \hepsilon \rightarrow \C$ is {\bf almost periodic} if $f( \cdot + i y) :\R \rightarrow \C$ 
is almost periodic  for all $ y > - (2 \pi)^{-1} \log \epsilon$. If $f: \hepsilon \rightarrow \C$ is an analytic almost periodic
function, then 
$$a_\lambda = \lim_{h \rightarrow \infty} \frac{1}{2h} \int^h_{-h} \mathe^{- 2 \pi i \lambda x} f(x + i y) dx$$
exists for all $\lambda \in \R$ and $a_\lambda$  is independent of $ y > - (2 \pi)^{-1} \log \epsilon$.
Moreover, $a_\lambda$ is non-zero for at most countably many $\lambda \in \R$. Furthermore, the 
{\bf Fourier series of $f$}:
$$\sum_{\lambda \in \R}  a_\lambda \mathe^{- 2 \pi i \lambda T}$$
converges to $f$ in the norm $\| \cdot \|_M$ defined in the space of analytic almost periodic functions with
domain $\hepsilon$ by:
$$\| g \|^2_M =  \lim_{h \rightarrow \infty} \frac{1}{2h} \int^h_{-h}  | g (x + iy) |^2 dx$$
for any  $y > - (2 \pi)^{-1} \log \epsilon$. 

\medskip
Let $\per  = \{ \nu \,\, /\, \omep \mbox{ is periodic  under } \psinu \}$. 
According to Proposition~\ref{uniform-p} we may consider $\epsilon >0$ such that
$\nu$ converges in $\H_{2 \epsilon}$ for all $\nu \in \per$. More precisely, 
if $\nu = \sum a_\lambda t^\lambda \in \per$, then  the sum
$$\sum_{\lambda \geq 0} a_\lambda \mathe^{2 \pi i \lambda T}$$ 
converges uniformly in $\H_{2 \epsilon}$ to a periodic analytic function $\bv_\nu (T)$. 

\begin{proposition}
\label{convergence-p}
Suppose that $\hnu \in \sol$ is such that $\omep \in \filled(\psi_\hnu)$ and $\bM^{\hnu} (\omep)$ is not periodic.
Then there exists $\bv_\hnu: \hepsilon \rightarrow \C$ such that if $\{\nu_n\} \subset \per$ is any sequence converging to
$\hnu$ in $\sol$, 
then $\bv_{\nu_n} \rightarrow \bv_\hnu$ uniformly in $\hepsilon$.  
\end{proposition}

Since the center of the $\psi$-parameter balls of level $n$ that contain $\hnu$ converge to $\hnu$ 
(Corollary~\ref{trivial-paraends-c}) the previous proposition
imposes a non-empty condition on $\bv_\hnu$.

\medskip
To prove the above proposition we start by showing that $\bv_{\nu_n}$ converges locally uniformly to some function
$\bv_\hnu$. For this we combine the Fourier analysis lemma below with the fact that,
according to Proposition~\ref{uniform-p}, the collection  $\{ \bv_\nu \,\,/\, \nu \in \per \}$ is a normal family.

\begin{lemma}
\label{fourier-l}
Let $\{f_n : \R \rightarrow \C \}$ be a uniformly bounded sequence of functions that converges locally uniformly
to a function $f$. 
  For all $n \geq 1$, suppose $$f_n (x) = \sum_{ \Lambda_n \ni \lambda \geq M_n} a^{(n)}_\lambda \mathe^{2 \pi i \lambda x}$$
where $\Lambda_n$ is a discrete subset of  $\R$, the sum converges uniformly and
$$M_n \rightarrow \infty \mbox{  as  } n \rightarrow \infty.$$
 Then $f(x) =0$ for all $x \in \R$. 
\end{lemma}

\noindent
{\bf Proof.}
Denote by $S(\R)$ the Schwarz space (see VI.4.1 in~\cite{katnelzon-76}).
To show that $f \equiv 0$ it is sufficient to prove that $\int f g =0$ for all $g \in Z$
where $$Z = \{ g \in S(\R) \,\,/\, \hat{g} \mbox{ has compact support} \}$$
and  $\hat{g}$ denotes the Fourier transform of $g$. In fact, $Z$ is dense in $S(\R)$ (see~\cite{gelfand-64} II.1.6) and
$g \mapsto \int f g$ is a continuous functional.

Observe that if $$h(x) = \sum_{\lambda \in \Lambda} a_\lambda \mathe^{2 \pi i \lambda x}$$
where $\Lambda$ is a discrete and bounded below subset of $\R$ and the sum converges uniformly,
then
$$\int_\R h(x) g(x) dx = \sum_{\lambda \in \Lambda} a_\lambda \hat{g}(\lambda)$$
for all $g \in S(\R)$.

Given $g \in Z$ the support of $\hat{g}$ is contained in $[-R,R]$ for some $R >0$.
Hence, there exists $N$ such that $M_n > R$ for all $n >N$. 
It follows that 
$$\int_\R f(x) \cdot g(x) dx = \lim_{n \rightarrow \infty} \int f_n (x) \cdot g(x) dx$$
$$= \lim_{n \rightarrow \infty} \sum_{\lambda \in \Lambda_n} a^{(n)}_\lambda \hat{g} (\lambda)  =0$$
since $ \sum_{\lambda \in \Lambda_n} a^{(n)}_\lambda \hat{g} (\lambda) =0$ for all $n >N$. 
\hfill $\Box$

\begin{corollary}
\label{convergence-c}
Suppose that $\hnu \in \sol$ is such that $\omep \in \filled(\psi_\hnu)$ and 
the critical marked grid $\bM^{\hnu} (\omep)$ is not periodic.
Then there exists $\bv_\hnu: \hepsilon \rightarrow \C$ such that if $\{\nu_n\} \subset \per$ is any sequence converging to
$\hnu$ in $\sol$, 
then $\bv_{\nu_n}$ converge $ \bv_\hnu$ locally uniformly in $\H_{2 \epsilon}$.  
\end{corollary}

\noindent
{\bf Proof.} 
Since $\{ \bv_{\nu_n} \}$ is a normal family we may suppose that a subsequence $\{ \bv_{\nu_{n_k}} \}$
converges locally uniformly to some function
$\bv_\nu$. If $\{ \bv_{\nu_n} \}$ does not converge to $\bv_\nu$, then there exists another subsequence 
$\{ \bv_{\nu_{m_k}} \}$ converging to some function $\tilde{v}_\nu$. In $\sol$ we have that
$\nu_{m_k} - \nu_{n_k} \rightarrow 0$. Therefore, by the previous lemma, 
$\bv_{\nu_{m_k}}(x - i (2 \pi)^{-1} \log \epsilon) - \bv_{\nu_{n_k}}( x - i (2 \pi)^{-1} \log \epsilon)$ 
converges locally uniformly to $0$ as
$k \rightarrow \infty$. Hence $\bv_\nu = \tilde{v}_\nu$ and the lemma follows.
\hfill $\Box$

\medskip
\noindent
{\bf Proof of Proposition~\ref{convergence-p}.}
By the Maximum Principle it is enough to show that the convergence is uniform in $\im T = - (2 \pi)^{-1} \log \epsilon$.
We proceed by contradiction and suppose that there exists $\delta >0$, $n_k \rightarrow \infty$,
$x_k \in \R$ such that 
$$|\bv_{\nu_{n_k}} (x_k - i \log \epsilon / 2 \pi) - \bv_{\nu}(x_k - i \log \epsilon / 2 \pi) | > \delta.$$
By passing to a subsequence we may also assume that there exists $x_\infty \in \R$, $M_k \in \Z$
and $|\Delta_k| <1/2$ such that
$$x_k = x_\infty +M_k + \Delta_k$$
and $\Delta_k \rightarrow 0$.
From the compactness given by Theorem~\ref{psi-parameter-th} and by Remark~\ref{automorphism-r}, we may
pass to a subsequence and  assume that
$\sigma^{M_k} (\nu) \rightarrow \hnu$. Thus,  $\sigma^{M_k} (\nu_k) \rightarrow \hnu$. 
It follows that
$$ \bv_{\nu_{n_k}} (x_k - i \log \epsilon / 2 \pi) = \bv_{\sigma^{M_k} \nu_k}(x_\infty + \Delta_k - i \log \epsilon / 2 \pi) \rightarrow
\bv_\hnu (x_\infty - i \log \epsilon / 2 \pi)$$
and similarly
$$ \bv_{\nu} (x_k - i \log \epsilon / 2 \pi) = \bv_{\sigma^{M_k} \nu}(x_\infty + \Delta_k - i \log \epsilon / 2 \pi) \rightarrow
\bv_\hnu (x_\infty - i \log \epsilon / 2 \pi)$$
which is a contradiction.
\hfill $\Box$

\medskip
Since the uniform limit of periodic functions is an almost periodic function and the corresponding Fourier coefficients also
converge we obtain the following (see~\cite{besicovitch-54}):

\begin{corollary}
If $\hnu = \sum_{\lambda \geq 0} a_\lambda t^\lambda \in \sol$ is such that $\omep \in \filled (\psihnu)$ and 
$\bM^\hnu (\omep)$ is not periodic, then $\bv_{\hnu} : \hepsilon \rightarrow \C$ is an almost periodic function 
with Fourier series $\sum_{\lambda \geq 0} a_\lambda \mathe^{2 \pi i  \lambda T} $.
\end{corollary}

Now the map $\tilde{\Psi} (T,\nu) = (\mathe^{-2 \pi i  \lambda T},\mathe^{-2 \pi i  \lambda T} \bv_\nu (T))$ is
defined and continuous for all $(T,\nu) \in \hepsilon \times \hA_\S$.

\begin{lemma}
 There exists $k \in \Z$ such that $T - T^\prime =k$ and
$\sigma^k(\nu) = \nu^\prime$ if and only if $\tilde{\Psi} (T, \nu) = \tilde{\Psi}(T^\prime, \nu^\prime)$.
\end{lemma}

\noindent
{\bf Proof.}
It is easily verified that if $\sigma^k(\nu) = \nu^\prime$ then $\bv_\nu (T+k) = \bv_{\nu^\prime} (T)$. Hence,
$  \tilde{\Psi} (T, \nu)  = \tilde{\Psi}(T -k , \nu^\prime)$.

Looking at the first coordinates we have that
 if $  \tilde{\Psi} (T, \nu) = \tilde{\Psi}(T^\prime, \nu^\prime)$, then $T - T^\prime =k $ for some $k \in \Z$
and $\bv_{\sigma^k \nu} (T^\prime) = \bv_\nu(T^\prime +k) = \bv_{\nu^\prime} (T^\prime)$.  
Let $\{ \nu_n \}$ and $\{ \nup_n \}$ be sequences in $\per $ converging to $\nu$ and $\nup$ respectively.
For $n$ large there exists $T^{''}$ such that $\bv_{\sigma^k \nu_n} (T^{''}) = \bv_{ \nup_n} (T^{''})$.
Thus $\sigma^k \nu_n = \nup_n$, for $n$ large, and therefore $\sigma^k \nu = \nup$.
\hfill $\Box$

\begin{lemma}
\label{image-l}
If $(T, \nu) \in  \hepsilon \times \hA_\S$,  $\bM$ is the critical marked grid of $\psinu$ and $(a,v) =\tilde{\Psi}(T,\nu)$,
then $(a,v) \in \hat{A}_\C \setminus \operatorname{int}\hA_\C  $ and  the critical marked grid of $\fav$ is $\bM$.
Moreover, the image of $\tilde{\Psi}$ is $\hat{A}_\C \setminus \operatorname{int}\hA_\C  $.
\end{lemma}

\noindent
{\bf Proof.}
Let  $\{ \nu_k \} \subset \per$ be a sequence converging to $\nu$.
It follows that $\{(a, v_k(T))=(\mathe^{- 2 \pi i T}, \mathe^{- 2 \pi i T} \bv_{\nu_k} (T))\}$ converges to $(a,v)$.
It suffices to show that  the critical marked grid $\bM^f$ of $\fav$ coincides with $\bM$.
Fix $n \geq 0$ and denote by $\bM^f_{n+1}$ the level $n+1$ critical marked grid of $\fav$.
There exists an integer $k$ sufficiently large so that $\nu_k$ and $\nu$ are in the same $\psi$-parameter ball of level $n$.
Therefore $\bM^{(a, v_k(T))}_{n+1} (a) = \bM_{n+1}$ (Lemma~\ref{locally-constant-l}).
By Proposition~\ref{psi-level-n-p} $\bM^\nu_{n+1} (\omep) = \bM^{\nu_k}_{n+1} (\omep)$ and
by Proposition~\ref{correspondence-p} $\bM^{\nu_k}_{n+1} (\omep) = \bM^f_{n+1}$.
Since this occurs for any $n \geq 0$, it follows that $\bM^f = \bM $.

Now we have to show that any $(a,v) \in \hA_\C \setminus \operatorname{int}\hA_\C$ with $|a| > 1/ \epsilon$  is in the image of
$\tilde{\Psi}$. In fact, by Theorem~4.2 in~\cite{mcmullen-94}, there exists a sequence $\{ (a_k,v_k) \}$ converging to
$(a,v)$ such that $a_k$ is periodic under $f_{a_k, v_k}$. It follows that 
$(a_k, v_k ) = (\mathe^{- 2 \pi i T_k}, \mathe^{- 2 \pi i T_k} \bv_{\nu_k} (T_k))\}$ for some $\nu_k \in \per$.
Without loss of generality we may suppose that $\{ T_k \}$ converges to some value $T$.
Let $\bM$ be critical marked grid of $\fav$. Note that $\bM$ is a not periodic.
Let $\bM_{n+1}$ denote the level $n+1$ critical marked grid of $\fav$.
From Proposition~\ref{psi-level-n-p} it follows that if  $r_n$ denotes the radius of the $\psi$-parameter balls of level $n$
that contain parameters with level $n+1$ grid $\bM_{n+1}$, then $r_n \rightarrow 0$. 
By Lemma~\ref{locally-constant-l} and Proposition~\ref{correspondence-p}, there exists $k_0$ such 
that $\bM_{n+1} = \bM^{\nu_k} (\omep)$ for all $k \geq k_0$. In particular, if we consider a nested sequence $\{ \pb_n \}$
of $\psi$-parameter  balls so that each contains infinitely many elements of $\{ \nu_k \}$, then $\cap \pb_n$ is
a singleton, say $\{ \nu \} \subset \hA_\S$, 
and there is a subsequence of $\{ \nu_k \}$ converging to $\nu$. It follows that 
$\tilde{\Psi}(T, \nu) = (a,v)$.
\hfill $\Box$

\begin{corollary}
   $\operatorname{int}\hA_\C \cap \{ |a| > 1/\epsilon \} = \emptyset$.
\end{corollary}

\noindent
{\bf Proof.}
If $\operatorname{int}\hA_\C \neq \emptyset$, then there exists $T \in \hepsilon$ and an aperiodic critical marked grid
$\bM$ such that the closed and bounded set 
$S_\bM = \{ v \,\,/\, a = \mathe^{- 2 \pi i T} \mbox{  and  } \bM^{a,v} (a) = \bM \} \subset \C$ 
has non-empty interior. According to Theorem~\ref{psi-parameter-th}, the set $C_\bM$ of all parameters $\nu$ such that 
$\psinu$ has critical marked grid $\bM$ is compact, non-empty and totally disconnected. 
By the previous results, $\tilde{\Psi}(T, \cdot) : C_\bM \rightarrow  S_\bM$ is continuous, one--to--one and
contains $\partial C_\bM$. But since $C_\bM$ is compact, $\tilde{\Psi}(T, \cdot)$ is a homeomorphism between
the totally disconnected set $C_\bM$ and its image. Therefore, $\partial S_\bM$ is totally disconnected and hence
$S_\bM$ is totally disconnected. 
\hfill $\Box$

\begin{lemma}
  The function $\Psi$ is a topological  embedding.
\end{lemma}

\noindent
{\bf Proof.}
It suffices to show that $\Psi$ is a closed map.
For this let $X \subset \hepsilon \times \hA_\sol$ be the preimage of a closed set under the quotient map. 
We must show that $\tilde{\Psi}(X)$ is closed. In fact, consider a sequence $\{ (T_k, \nu_k ) \}\subset X$ such that
$\tilde{\Psi} (T_k, \nu_k ) \rightarrow (a,v) \in \hA_\C $. 
Since $X$ is the preimage of a set under the quotient map, by Remark~\ref{automorphism-r},
we may assume that $\{T_k \}$ converges, say to $T$. Then, as in the proof of Lemma~\ref{image-l}, by passing to a subsequence one
may suppose that $\{ \nu_k \}$ converges to some $\nu \in \hA_\sol$. 
It follows that $(T, \nu) \in X$ and $(a,v) \in \tilde{\Psi}(X)$.
\hfill $\Box$

\bibliographystyle{plain}

\end{document}